\documentclass[15pt,oneside,reqno]{amsart}

\usepackage[margin=20truemm]{geometry}
\usepackage{amsmath,amssymb,amsthm}
\usepackage{float}
\usepackage{shuffle}
\usepackage{multicol}
\usepackage{ascmac}
\usepackage{mathtools}
\usepackage{bm}
\usepackage{amscd}
\usepackage{comment}
\usepackage{mathrsfs}
\usepackage[all]{xy}
\usepackage{color}
\usepackage[dvipsnames]{xcolor}
\usepackage{setspace}
\usepackage{cases}
\usepackage{paralist}
\usepackage{autobreak}
\usepackage{enumitem}
\usepackage{stackrel}

\definecolor{Mycolor}{cmyk}{ 0.90, 0.00, 0.34, 0.35} 

\usepackage{iftex}
\ifptex
  \usepackage[dvipdfmx,bookmarksnumbered,colorlinks]{hyperref}
  \usepackage{pxjahyper}
\else
  \usepackage[bookmarksnumbered,colorlinks]{hyperref} 
\fi

\hypersetup{%
 setpagesize=false,%
 bookmarks=true,%
 bookmarksdepth=tocdepth,%
 bookmarksnumbered=true,%
 colorlinks=true,
 linkcolor=Mycolor,   
 citecolor=blue,   
 pdftitle={},%
 pdfsubject={},%
 pdfauthor={Takumi Anzawa},%
 pdfkeywords={SMZVs}}

\title[A Lie algebra associated with AMZVs]{A Lie algebra associated with adjoint multiple zeta values}
\author{Takumi Anzawa}
\address{Center for Promotion of Higher Education, Kogakuin University, 2665-1 Nakano-machi, Hachioji-shi, Tokyo 192-0015, Japan}
\email{atakumi86@gmail.com}

\newcounter{dfcounter}

\makeatletter
\@addtoreset{dfcounter}{chapter}    
\@addtoreset{dfcounter}{section} 
\makeatother

\newtheorem{df}{Definition}[dfcounter]
\newtheorem{thm}[df]{Theorem}
\newtheorem{cor}[df]{Corollary}

\newtheorem{conj}[df]{Conjecture}
\newtheorem{qu}[df]{Question}
\newtheorem{mt}{Main Theorem}

\theoremstyle{remark}{
\newtheorem{lem}[df]{Lemma}
\newtheorem{prop}[df]{Proposition}
\newtheorem{ex}[df]{Example}
\newtheorem{rmk}[df]{Remark}
}

\DeclareMathOperator{\Hom}{Hom}

\DeclareMathOperator{\Span}{Span}

\newcommand{\ide}[1]{\mathfrak{#1}}

\newcommand{\del}{\partial}

\newcommand{\engMonth}{
	\ifcase\month
		\or January 
		\or February 
		\or March 
		\or April
		\or May 
		\or June
		\or July 
		\or August 
		\or September
		\or October 
		\or November 
		\or December\fi
} 
\newcommand{\fullday}{
	\engMonth \, \the\day, \the\year	
} 

\newcommand{\Z}{\mathbb{Z}}
\newcommand{\Q}{\mathbb{Q}}

\newcommand{\C}{\mathbb{C}}

\let\Im\relax
\newcommand{\Im}{\textrm{Im}}

\newcommand{\DMR}{\mathrm{DMR}}
\newcommand{\dmr}{\mathfrak{dmr}}

\newcommand{\AdDMR}{\mathrm{AdDMR}}
\newcommand{\addmr}{\mathfrak{addmr}}

\newcommand{\der}{\mathfrak{der}}

\newcommand{\Ad}{\mathrm{Ad}}
\newcommand{\ad}{\mathfrak{ad}}

\newcommand{\tm}{\mathfrak{tm}}
\newcommand{\Fad}{\mathfrak{F}_2^{\ad(x_1)}}
\newcommand{\FAd}{\mathfrak{F}_2^{\Ad(x_1)}}

\DeclareMathOperator{\inv}{inv}

\newcommand{\cotimes}{\widehat{\otimes}}

\DeclareMathOperator{\EDS}{EDS}

\let\S\relax
\newcommand{\S}{\mathcal{S}}

\newcommand{\Admzv}[2]{\zeta_{\mathrm{Ad}}(#2;#1)}

\newcommand{\sAdmzv}[2]{\zeta_{\mathrm{Ad}}^{\shuffle}(#2;#1)}
\newcommand{\hAdmzv}[2]{\zeta_{\mathrm{Ad}}^{*}(#1;#2)}

\DeclareMathOperator{\sft}{sft}

\DeclareMathOperator{\mi}{mi}
\DeclareMathOperator{\ma}{ma}
\DeclareMathOperator{\vimo}{vimo}

\let\x\relax
\newcommand{\x}{\mathrm{x}}
\newcommand{\y}{\mathrm{y}}
\newcommand{\z}{\mathrm{z}}
\let\Y\relax
\newcommand{\X}{\mathrm{X}}
\let\Y\relax
\newcommand{\Y}{\mathrm{Y}}
\let\V\relax
\newcommand{\V}{\mathrm{V}}

\let\u\relax
\newcommand{\u}{\mathrm{u}}
\let\v\relax
\newcommand{\v}{\mathrm{v}}
\let\w\relax
\newcommand{\w}{\mathrm{w}}

\newcommand{\strprty}{\mathrm{str.prty}}

\def\utilde#1{\mathord{\vtop{\ialign{##\crcr
$\hfil\displaystyle{#1}\hfil$\crcr\noalign{\kern1.5pt\nointerlineskip}
$\hfil\tilde{}\hfil$\crcr\noalign{\kern1.5pt}}}}}

\newcommand{\Qalg}{\mathbb{Q}\text{-}\mathrm{\mathbf{Alg}}}

\DeclareMathOperator{\dep}{dep}

\DeclareMathOperator{\wt}{wt}

\allowdisplaybreaks[0]

\makeatletter
\newcommand{\nsubsection}[1]{%
  \par\bigskip
  \noindent\textbf{#1}\par
  \medskip
}
\makeatother

\begin{document}
\maketitle

\begin{abstract}
Jarossay (arXiv math.NT1412.5099) introduced adjoint multiple zeta values and, by using Racinet's dual formulation of the generating series of multiple zeta values, found $\mathbb{Q}$-algebraic relations among them, referred to as the \textit{adjoint double shuffle relations}.
Additionally, Jarossay defined the affine scheme $\mathrm{AdDMR}_0$ determined by the adjoint double shuffle relations and posed a question whether $\mathrm{AdDMR}_0$ is isomorphic to Racinet's double shuffle group $\mathrm{DMR}_0$ (Publ. Math. Inst. Hautes \'{E}tudes Sci. (2002), no. 95). 
In this paper, we refine Jarossay's question by introducing the condition referred to as the adjoint conditions, and, based on this refinement, we study the corresponding Lie algebraic aspect.
Within this framework, we construct the Lie algebra associated with the adjoint double shuffle relations by imposing Hirose's parity results.
\end{abstract}

\tableofcontents

\section{Introduction}

\subsection{Multiple zeta values}

For a tuple of integers $(k_1,\ldots,k_r)\in\mathbb{Z}_{>0}^r$ with $k_r>1$, a \emph{multiple zeta value} is defined by

\[
\zeta(k_1,\ldots,k_r) := \sum_{0<m_1<\cdots <m_r} \frac{1}{m_1^{k_1}\cdots m_r^{k_r}}.
\]
The condition $k_r>1$ ensures the convergence of the above multiple series.
We set $\zeta(\emptyset) := 1$.
Let $\mathcal{Z}$ be the $\mathbb{Q}$-linear space spanned by MZVs.
When $r=1$, $\zeta(k_1)$ coincides with the special values of the Riemann zeta function.
MZVs were first studied by Euler \cite{Euler}, who investigated double zeta values (the case $r=2$) and obtained the $\Q$-linear relations among them, called the sum formula.
Around 1990, Hoffman \cite{Hoffman-92} and Zagier \cite{Zagier-92} revisited MZVs and their applications.

The extended double shuffle relations are a family of $\mathbb{Q}$-linear relations among MZVs that arise from the combination of two kinds of products on $\mathcal{Z}$. The first is the \textit{shuffle product} $\shuffle$, which comes from the iterated integral expressions, and the second is the \emph{harmonic product} $*$, which comes from the multiple series expressions.

In \cite{Ihara-Kaneko-Zagier}, Ihara, Kaneko, and Zagier constructed two regularizations of MZVs, the $\shuffle$-regularization $\zeta^{\shuffle}(k_1,\ldots,k_r)$ and the $*$-regularization $\zeta^{*}(k_1,\ldots,k_r)$, for $k_1,\ldots,k_r\in\mathbb{Z}_{>0}$.
These regularizations extend the definition of MZVs to indices with $k_r=1$ and are compatible with the shuffle and harmonic products, respectively.
Using these two regularizations of MZVs, the authors of \cite{Ihara-Kaneko-Zagier} formulated the class of $\mathbb{Q}$-linear relations among MZVs, called the \emph{extended double shuffle relations}.

By $\mathcal{Z}^f$, we denote a unital $\Q$-algebra generated by formal symbols $\zeta^f(k_1,\ldots,k_r)$ ($r\in\mathbb{Z}_{>0}$, $k_1,\ldots,k_r\in\mathbb{Z}_{>0}$) which satisfy the extended double shuffle relations. The $\Q$-algebra $\mathcal{Z}^f$ introduced in \cite{Ihara-Kaneko-Zagier} is called the formal multiple zeta space, and its generators $\zeta^f(k_1,\ldots,k_r)$ are called formal multiple zeta values.
Conjecturally, the map $\mathcal{Z}^f\to\mathcal{Z}$ sending $\zeta^f(k_1,\ldots,k_r)$ to $\zeta^{\shuffle}(k_1,\ldots,k_r)$ (for $r\ge1$, $(k_1,\ldots,k_r)\in\mathbb{Z}_{>0}^r$) is a $\mathbb{Q}$-algebra isomorphism.
In particular, this implies that the extended double shuffle relations generate all $\mathbb{Q}$-linear relations among MZVs.

Let $x_0$, $x_1$ be letters and, for $\Q$-algebra $R$, $R\langle \langle x_0,x_1 \rangle\rangle$ denotes the algebra of non-commutative power series generated by $x_0$, and $x_1$ over $R$.
Independent of \cite{Ihara-Kaneko-Zagier}, Racinet \cite{Racinet} also formulated the extended double shuffle relations via generating series of $\shuffle$-regularized MZVs:
\[
\Phi_{\shuffle} := 1 + \zeta^{\shuffle}(2) x_0x_1 - \zeta^{\shuffle}(2) x_1x_0
+ \zeta^{\shuffle}(3)x_0x_0x_1 + \zeta^{\shuffle}(1,2)x_1x_0x_1 + \cdots \in \mathcal{Z}\langle\langle x_0,x_1 \rangle\rangle.
\]
(More precise definitions are provided in Section \ref{subsect:DMR}).
Using this series, Racinet \cite{Racinet} provided a dual framework for the extended double shuffle relations.
Moreover, for a $\Q$-algebra $R$, he introduced the set $\DMR_0(R) \subset R\langle \langle x_0,x_1\rangle\rangle$ defined by the extended double shuffle relations together with the additional relation ``$\pi i = 0$".
Then, $\DMR_0(R)$ forms a group endowed with $\circledast$.

\subsection{Adjoint multiple zeta values}

For $(k_1,\ldots,k_r) \in \mathbb{Z}_{>0}^r$ and $l\in \mathbb{Z}_{\ge 0}$, we define the \emph{adjoint multiple zeta value} (AdMZVs) by
\begin{align*}
\Admzv{l}{k_1,\ldots,k_r} :=& \sum_{i = 0}^r (-1)^{k_{i+1}+\cdots+k_r + l} \zeta^{\shuffle}(k_1,\ldots,k_i)\zeta^{\shuffle}_l(k_r,\ldots,k_{i+1}) \mod \zeta(2),
\end{align*}
where
\[
\zeta_{l}^{\shuffle}(k_r,\ldots,k_{i+1}) := (-1)^l \sum_{\substack{l_{i+1}+\cdots+l_r = l \\ l_{i+1},\ldots,l_r\ge 0}}
\left(\prod_{j = i+1}^r \binom{k_j + l_j-1}{l_j}\right)
\zeta^{\shuffle}(k_r+l_r,\ldots,k_{i+1}+l_{i+1}).
\]
Note that, for $k_1,\ldots,k_r \in \mathbb{Z}_{>0}$, $\zeta_{\S}^{\shuffle}(k_1,\ldots,k_r) = \sAdmzv{0}{k_1,\ldots,k_r}$ holds.
It is known that adjoint multiple zeta values are a generalization of symmetric multiple zeta values. Indeed, the definition of a symmetric multiple zeta value is
\[
\zeta_{\S}(k_1,\ldots,k_r) := \Admzv{0}{k_1,\ldots,k_r}
\]
for $(k_1,\ldots,k_r)$ (\cite{Kaneko-Zagier}).
The motivation for this work is that AdMZVs arise in the coefficients of $\Phi_{\shuffle}^{-1}x_1 \Phi_{\shuffle}$.
More precisely, for $k_1,\ldots,k_r \in \mathbb{Z}_{>0}$ and $l\in \mathbb{Z}_{\ge 0}$, the following
\[
\Admzv{l}{k_1,\ldots,k_r} = \langle \Phi_{\shuffle}^{-1}x_1\Phi_{\shuffle} \mid x_0^lx_1x_0^{k_r-1}x_1\cdots x_0^{k_1-1}x_1 \rangle \mod \zeta(2)\mathcal{Z}
\]
holds.
From this generating series, Jarossay derived $\mathbb{Q}$-algebraic relations among AdMZVs, called the \emph{adjoint double shuffle relations}, and defined the set
\[
\AdDMR_0(R)\subset R\langle\!\langle x_0,x_1\rangle\!\rangle .
\]
It is known that there is an injective map from Racinet's $\DMR_0(R)$ to Jarossay's $\AdDMR_0(R)$:
\[
\Ad^R(x_1) : \DMR_0(R) \rightarrow \AdDMR_0(R) : \phi \mapsto \phi^{-1}x_1\phi
\]
(See Proposition \ref{prop:t-adic}).
Related to this embedding, Jarossay posed the following question:
\begin{qu}[ = Question~\ref{qu:Jarossay}]\label{qu:Jarossay-qu}
Is the embedding
\[
\Ad^R(x_1) : \DMR_0(R) \rightarrow \AdDMR_0(R) : \phi \mapsto \phi^{-1}x_1\phi
\]
surjective?
\end{qu}

\subsection{Main Theorems}

Let $\ide{F}_2$ be the degree completion of the free Lie algebra over $\Q$ generated by $x_0$ and $x_1$. For $k \in \mathbb{Z}_{\ge 1}$, let $\ide{F}_2^{\ge k}$ denote the Lie subalgebra of $\ide{F}_2$ consisting of elements with no terms of degree less than $k$.
To study Question \ref{qu:Jarossay-qu}, we consider the tangent spaces of $\AdDMR_0(\Q)$ and $\DMR_0(\Q)$.
Let $\dmr_0$ be the tangent space of $\DMR_0(\Q)$ at $1$, and $\addmr_0$ the tangent space of $\AdDMR_0(\Q)$ at $\Ad^{\Q}(x_1)(1)=x_1$.
At the level of the tangent space, Question~\ref{qu:Jarossay-qu} becomes the following statement:
\begin{qu}[ = Question~\ref{qu:Jarossay-Lie}]\label{qu:Jarossay-qu-tan}
Is the embedding
\[
\ad(x_1) : \dmr_0 \rightarrow \addmr_0 : \psi \mapsto [x_1,\psi]
\]
surjective? Here, for $f$, $g \in \Q\langle\langle x_0,x_1 \rangle\rangle$, we define $[f,g] := fg - gf$.
\end{qu}
In fact, the answer is negative.
\begin{prop}\label{prop:omega}
The $\Q$-linear map stated in Question \ref{qu:Jarossay-qu-tan} is not $\Q$-linear isomorphism.
More precisely, put
\[
\omega :=
-\bigl[[x_0,x_1],[x_0,[x_0,x_1]]\bigr]
+\bigl[[x_0,x_1],[x_1,[x_0,x_1]]\bigr]
-\bigl[x_1,[x_0,[x_0,[x_0,x_1]]]\bigr]
+\bigl[x_1,[x_1,[x_0,[x_0,x_1]]]\bigr].
\]
Then $\omega\in \addmr_0$, but $\omega\notin \mathrm{Im}\,\ad(x_1)(\dmr_0)$.
\end{prop}
As a consequence of Proposition \ref{prop:omega}, we obtain the following:
\begin{mt}\label{mt:Jarossay-negative}
The answer to Question \ref{qu:Jarossay-qu} is also negative.
In particular, if we take $R=\Q[\varepsilon]/(\varepsilon^2)$, then the map
\[
\Ad^R(x_1)\colon \DMR_0(R)\longrightarrow \AdDMR_0(R),\qquad
\phi \longmapsto \phi^{-1}x_1\phi
\]
is not surjective.
\end{mt}
Due to Main Theorem \ref{mt:Jarossay-negative}, we refine Questions~\ref{qu:Jarossay-qu} and \ref{qu:Jarossay-qu-tan}.
Let $\ide{F}_2$ be the degree completion of the free Lie algebra over $\Q$ generated by $x_0$ and $x_1$.
For $k \in \mathbb{Z}_{>0}$, let $\ide{F}_2^{\ge k}$ denote the Lie subalgebra of $\ide{F}_2$ consisting of elements whose coefficients of words of degree less than $ k $ vanish.
We define $\Q$-linear space $\Fad$ by
\[
\Fad := \{\Psi \in \ide{F}_2^{\ge 3} \mid {}^\exists\, \psi \in \ide{F}_2^{\ge 2}\ \text{such that}\ \Psi = [x_1,\psi]\}.
\]
\begin{qu}[ = Question~\ref{qu:addmr-Lie}]\label{qu:Jarossay-Lie-intro}
Is the $\Q$-linear map
\[
\ad(x_1)\colon \dmr_0 \longrightarrow \addmr_0 \cap \Fad,\qquad \Psi \longmapsto [x_1,\Psi]
\]
an isomorphism?
\end{qu}
For each unital commutative associative $\Q$-algebra $R$, we define
\[
\FAd(R)
:=
\left\{
\phi^{-1}x_1\phi
\ \middle|\
\phi\in 1+R\langle\langle x_0,x_1\rangle\rangle_{\ge 2},
\quad
\Delta_{\shuffle}(\phi)=\phi\otimes\phi
\right\}.
\]
and consider $\AdDMR_0(R) \cap \FAd(R)$, where $\Delta_{\shuffle} : \Q\langle\langle x_0,x_1 \rangle\rangle \rightarrow \Q\langle\langle x_0,x_1 \rangle\rangle^{\cotimes 2}$ is a shuffle coproduct (see Section \ref{subsect:dual-of-h}).
Then, we have the following embedding:
\[
\Ad^R(x_1)\colon \DMR_0(R) \longrightarrow \AdDMR_0(R) \cap \FAd(R),\qquad \phi \longmapsto \phi^{-1} x_1 \phi.
\]
We note that the tangent space of $\AdDMR_0(\Q)\cap \FAd(\Q)$ at $x_1$ is
\[
T_{x_1}\bigl(\AdDMR_0(\Q)\cap \FAd(\Q)\bigr)=\addmr_0 \cap \Fad,
\]
by repeating the arguments of Propositions~\ref{prop:addmr} and \ref{prop:dAd}.
Thus, we establish the following Question which is the refinement of Question~\ref{qu:Jarossay}
\begin{qu}[ = Question~\ref{qu:AdDMR-grp}]\label{qu:AdDMR-grp-intro}
Is the map
\[
\Ad^R(x_1)\colon \DMR_0(R) \longrightarrow \AdDMR_0(R)\cap \FAd(R),\qquad \phi \longmapsto \phi^{-1} x_1 \phi
\]
a bijection for each unital commutative associative $\Q$-algebra $R$?
\end{qu}
Another main theorem provides further evidence for Questions \ref{qu:Jarossay-Lie-intro} and \ref{qu:AdDMR-grp-intro}.
If Question \ref{qu:Jarossay-Lie-intro} has a positive answer, then we expect that $\addmr_0 \cap \Fad$ carries a Lie algebra structure induced by $\ad(x_1)$.
In this theorem, we impose an additional condition on $\addmr_0 \cap \Fad$:
\[
V_{\strprty} := \{\Psi \in \Q\langle\langle x_0,x_1 \rangle\rangle \mid \Psi^{11} + \Psi^{10} + \Psi^{01} = 0\},
\]
where, for $\Psi \in \ide{h}^{\vee}$, we write
\[
\Psi = x_0 \Psi^{00} x_0 + x_1 \Psi^{10} x_0 + x_0 \Psi^{01} x_1 + x_1 \Psi^{11} x_1.
\]
Then, we obtain the following:
\begin{mt}[ = Main Theorem~\ref{mt:addmr-Lie}]
The pair $(\addmr_0 \cap \Fad \cap V_{\strprty}, \{\mathchar`-.\mathchar`-\}_1)$ forms a Lie algebra. See Remark \ref{rmk:addmr-bracker} for the definition of $\{\mathchar`-.\mathchar`-\}_1$.
\end{mt}

\section{Extended double shuffle relations}\label{sect:EDSR}

In this section, we first introduce the extended double shuffle relations.
The extended double shuffle relations were established independently in \cite{Ihara-Kaneko-Zagier} and \cite{Racinet}.
After setting up the algebraic framework, we review the works of \cite{Ihara-Kaneko-Zagier} and \cite{Racinet}, respectively.

\subsection{Shuffle algebra and harmonic product algebra}

In \cite{Hoffman}, Hoffman introduced certain commutative $\mathbb{Q}$-algebras for studying MZVs.
In this subsection, based on Hoffman's $\Q$-algebras, we represent two types of $\Q$-algebras.

Let $X := \{x_0,x_1\}$ be a set of letters, $X^*$ the free monoid generated by $X$ and $\ide{h}:= \mathbb{Q}\langle X\rangle$ the free associative $\mathbb{Q}$-algebra generated by $X$.
We define the weight $\wt w$ of $w \in X^*$ by the number of letters in $w$.

First, we set up the shuffle product $\shuffle$ on $\ide{h}$.
It is defined bilinearly and recursively by $1\shuffle w := w \shuffle 1 := w$ and 
\begin{align*}
l_1w_1\shuffle l_2w_2 =& l_1(w_1\shuffle l_2w_2) + l_2(l_1w_1 \shuffle w_2 )
\end{align*}
for all $l_1$, $l_2\in X$ and $w$, $w_1$, $w_2\in X^*$.
The pair $(\ide{h}, \shuffle)$ forms a unital and commutative $\mathbb{Q}$-algebra and we denote it by $\ide{h}^{\shuffle}$.

We define two $\mathbb{Q}$-subalgebras $\ide{h}^0$ and $\ide{h}^1$ of $\ide{h}$ by
\[
\ide{h}^1:= \mathbb{Q} + \bigoplus_{\substack{w\in x_{1}X^*}} \mathbb{Q} w
\supset
\ide{h}^0 := \mathbb{Q} +  \bigoplus_{\substack{w\in x_{1}X^*\\  w\notin X^*x_1 }} \mathbb{Q} w.
\]
Then, $\ide{h}^1$ and $\ide{h}^0$ are also closed under $\shuffle$ and become $\mathbb{Q}$-subalgebras of $\ide{h}^{\shuffle}$.
We denote them by $\ide{h}^{1,\shuffle}$ and $\ide{h}^{0,\shuffle}$, respectively.
\begin{rmk}{}{}
Equipped with the weight grading on $\ide{h}$, $\ide{h}^{\shuffle}$ is a graded $\mathbb{Q}$-algebra under the shuffle product.
\end{rmk}

Let $Y := \{y_{n}\}_{\substack{n\ge 1}}$ be the set of letters and $Y^*$ the free monoid generated by $Y$.
By abuse of notation, we also use $1$ to denote the empty word on $Y$. Let $\Q\langle Y \rangle$ be the free associative $\Q$-algebra generated by $Y$.
We define the $\Q$-subalgebra of $\Q\langle Y \rangle$ as follows: 
\[
\Q\langle Y \rangle^0 := \displaystyle \mathbb{Q} + \bigoplus_{\substack{k\in\mathbb{Z}_{>1}\\  w\in Y^*}} \Q wy_{k}.
\]

We define the $\mathbb{Q}$-bilinear binary operation $*$, called the harmonic product, on $\Q\langle Y \rangle$ inductively by $1 * w := w * 1 := w$ and
\begin{align*}
y_{k_1}w_1 * y_{k_2}w_2 :=& y_{k_1}(w_1 * y_{k_2} w_2) + y_{k_2}(y_{k_1}w_1 * w_2) + y_{k_1 + k_2}(w_1*w_2)
\end{align*}
for each letters $y_{k_1}$ and $y_{k_2} \in Y$ and each words $w$, $w_1$ and $w_2\in Y^*$.

The pair $(\Q\langle Y \rangle, *)$ constitutes a commutative associative $\mathbb{Q}$-algebra, which we denote by $\Q\langle Y \rangle^{*}$.
In the same way, $\Q\langle Y \rangle^0$ is also closed under $*$. Therefore, $\Q\langle Y \rangle^0$ constitutes a $\mathbb{Q}$-algebra with respect to $*$ and we denote it by $\Q\langle Y \rangle^{0,*}$.
\begin{rmk}{}{}
We define the weight of $\Q\langle Y \rangle$ by $\wt y_k := k$ for $k\in\mathbb{Z}_{>0}$,
which is preserved under the following natural $\Q$-linear map:
\[
\mathbf{p} : \Q\langle Y \rangle \rightarrow \ide{h}\quad ; \quad y_{k_1}\cdots y_{k_r} \mapsto x_1x_0^{k_1-1}\cdots x_1x_0^{k_r-1}.
\]
Thus, $\Q\langle Y \rangle$ constitutes a graded $\mathbb{Q}$-algebra with respect to the concatenation product and the harmonic product.
\end{rmk}

Before we finish this subsection, we define the $\mathbb{Q}$-linear map by
\begin{align*}
\mathbf{q} : \ide{h} \rightarrow \Q\langle Y \rangle \quad ; \quad \mathbf{q}(x_0^{k_0-1}x_1\cdots x_0^{k_1-1}x_1x_0^{k_{r}-1}) = 
\begin{cases}
y_{k_1}\cdots y_{k_r} & k_{0} = 1,\\
0& \text{otherwise.}
\end{cases}
\end{align*}
Then, $\mathbf{q}$ is the left inverse of $\mathbf{p}$.
Note that the notations $\mathbf{p}$ and $\mathbf{q}$ are based on \cite{Racinet}.

\subsection{Duals of $\ide{h}^{\shuffle}$ and $\Q\langle Y \rangle^*$}\label{subsect:dual-of-h}
This subsection discusses the dual of $\ide{h}^{\shuffle}$ and $\Q\langle Y \rangle^*$.
Let us denote by
\[
\ide{h}^{\vee} := \Hom(\ide{h},\Q)
\]
the $\Q$-linear dual of $\ide{h}$.
Then we have
\[
\ide{h}^{\vee} = \mathbb{Q}\langle \langle x_0,x_1 \rangle \rangle:= \varprojlim_{n}\ide{h}/ (x_0,x_1)^n,
\]
where $(x_0,x_1)$ is the two-sided ideal of $\ide{h}$.
 For $\Phi \in \ide{h}^{\vee}$, we denote $\Phi$ by
\[
\Phi = \sum_{w\in X^*} \langle \Phi \mid w \rangle w = \langle \Phi \mid 1 \rangle + \langle \Phi \mid x_0 \rangle x_0 + \langle \Phi \mid x_1 \rangle x_1 + \cdots \hspace{0.5cm}( \langle \Phi \mid w \rangle \in \mathbb{Q}).
\]
This notation induces a $\mathbb{Q}$-bilinear map $\langle \mathchar`- | \mathchar`-\rangle : \ide{h}^{\vee} \cotimes \ide{h} \rightarrow \mathbb{Q} \hspace{0.2cm} ; \hspace{0.2cm}  \Phi \otimes w \mapsto \langle \Phi \mid w \rangle$.
We define the shuffle coproduct $\Delta_{\shuffle} : \ide{h}^{\vee} \rightarrow  \ide{h}^{\vee\cotimes 2}$ by
\[
\Delta_{\shuffle}(\Phi) := \sum_{u,v\in X} \langle \Phi \mid u\shuffle v \rangle u\otimes v. 
\]
Then, $\Delta_{\shuffle}$ is a continuous algebra homomorphism and satisfies $\Delta_{\shuffle}(x_i) = x_i \otimes 1 + 1\otimes x_i$.
Let $S_X^{\vee} : \ide{h}^{\vee}\rightarrow \ide{h}^{\vee}$ be the anti-automorphism of $\ide{h}^{\vee}$ defined by $x_i \mapsto -x_i$ ($i = 0$, $1$).
Then, the tuple $(\ide{h}^{\vee}, \cdot, \Delta_{\shuffle},S^{\vee}_X)$ is a completed Hopf algebra that is topologically dual to $\ide{h}^{\shuffle}$.

Recall that $R$ is an unital, commutative, and associative $\Q$-algebra.
We often consider the coefficient expansion of $\ide{h}$ to $R$.
We endow $R$ with the discrete topology and define $R\langle \langle x_0,x_1 \rangle \rangle$ by
\[
R\langle \langle x_0,x_1 \rangle \rangle := \varprojlim_{n} R\otimes \ide{h} / (x_0,x_1)^n = \varprojlim_{n} R\langle x_0,x_1\rangle / (x_0,x_1)^n.
\]
We extend the operations defined on $\ide{h}^{\vee}$ to $R \langle \langle x_0,x_1 \rangle \rangle$ via the scalar expansions.

In the same way, we shall construct the dual of $\Q\langle Y \rangle^*$.
We define $\Q\langle\langle Y \rangle\rangle$ by
\begin{align*}
\Q\langle\langle Y \rangle\rangle := \varprojlim_{n} \Q\langle Y \rangle / Y_{\ge n},
\end{align*}
where, for $n \in \mathbb{Z}_{>0}$, we define  $Y_{\ge m}\subset \Q\langle Y\rangle$ as the $\Q$-subspace spanned by words of weight $\ge m$.
By abuse of notation, for $\Phi \in \Q\langle\langle Y \rangle\rangle$, we denote $\Phi$ by
\[
\Phi = \sum_{w\in Y^*} \langle \Phi \mid w \rangle w = \langle \Phi \mid 1 \rangle + \langle \Phi \mid y_1 \rangle y_1 + \langle \Phi \mid y_2 \rangle y_2 + \cdots \hspace{0.5cm}( \langle \Phi \mid w \rangle \in \mathbb{Q}).
\]
This notation induces a $\mathbb{Q}$-bilinear map $\langle \mathchar`- | \mathchar`- \rangle : \Q\langle\langle Y \rangle\rangle \cotimes \Q\langle Y \rangle \rightarrow \mathbb{Q} \hspace{0.2cm} ; \hspace{0.2cm}  \Phi \otimes w \mapsto \langle \Phi \mid w \rangle$.
We define the harmonic coproduct $\Delta_{*} : \Q\langle\langle Y \rangle\rangle \rightarrow  \Q\langle\langle Y \rangle\rangle^{\cotimes 2}$ by
\[
\Delta_{*}(\Phi) := \sum_{u,v\in Y^*} \langle \Phi \mid u* v \rangle u\otimes v. 
\]
Then, $\Delta_{*}$ is a continuous algebra homomorphism which satisfies $\Delta_{*}(y_k) = y_k \otimes 1 + 1\otimes y_k + \sum_{\substack{i + j = k \\ i,j>0}} y_i \otimes y_j$ for $k\in\mathbb{Z}_{>0}$.
A tuple $(\Q\langle\langle Y \rangle\rangle, \cdot, \Delta_{*})$ constitutes a completed Hopf algebra which is topologically dual to $\Q\langle  Y \rangle^*$.

Similarly to $\ide{h}^{\vee}$, we often consider the coefficient expansion $\Q\langle \langle Y \rangle \rangle$ to $R\langle \langle Y \rangle \rangle$.
Recall that we endow $R$ with the discrete topology, and we define
\begin{align*}
R\langle \langle Y \rangle \rangle := \varprojlim_{n} R\otimes \Q\langle Y \rangle / Y_{\ge n} = \varprojlim_{n} R \langle Y \rangle / Y_{\ge n}
\end{align*}
We extend the operations defined on $\Q \langle \langle Y \rangle \rangle$ to $R \langle \langle Y \rangle \rangle$.

Before we finish this subsection, we define a continuous $\Q$-linear map $\mathbf{q}^{\vee} : \ide{h}^{\vee} \rightarrow\Q\langle\langle Y \rangle\rangle$ by
\[
\mathbf{q}^{\vee}(x_0^{k_1-1}x_1\cdots x_0^{k_r-1}x_1x_0^{k_{r+1}-1}) 
=
\begin{cases}
y_{k_1}\cdots y_{k_r}& \text{if } k_{r+1}=1, \\
0& \text{otherwise.}
\end{cases}
\]

\subsection{Completed free Lie algebra $\ide{F}_2$}

We define the completed free Lie algebra $\ide{F}_2:= \{\Psi \in \ide{h}^{\vee} \mid \Delta_{\shuffle}(\Psi) = \Psi \otimes 1 + 1 \otimes \Psi\}$ generated by $X$, 
namely the subspace of primitive elements with respect to the shuffle coproduct. Equivalently, a series $\Psi\in\ide{h}^{\vee}$ lies in $\ide{F}_2$ if and only if $\langle \Psi\mid u\shuffle v\rangle=0$ for all nonempty words $u,v\in X^*$.
Then the universal enveloping ring of $\ide{F}_2$ is isomorphic to $\ide{h}^{\vee}$ as Hopf algebras.
Moreover, for $\Psi \in \ide{F}_2$,
\begin{align}
S_X^{\vee}(\Psi) = - \Psi \label{eq:antipode-X}
\end{align}
holds.
To finish this subsection, for each positive integer $k$, we define
\begin{align*}
\ide{F}_2^{k} :=& \{\Psi \in \ide{F}_2 \mid\langle \Psi  \mid w \rangle = 0\text{ for $w\in X^*$ with $\wt w\neq k$}\}, \\
\ide{F}_2^{\ge k} :=& \{ \Psi \in \ide{F}_2 \mid \langle \Psi \mid w \rangle = 0\text{ for } \wt w \le k-1   \}.
\end{align*}

\subsection{Extended double shuffle relations}

MZVs satisfy the shuffle relations and the stuffle (harmonic) relations, coming from their integral and series expressions.
Combining these two product-to-sum relations, one can derive the finite double shuffle relations. 
However, 
the finite double shuffle relations do not suffice all $\Q$-linear relations among MZVs (for instance, see \cite[pp.7]{Ihara-Kaneko-Zagier}).
The authors of \cite{Ihara-Kaneko-Zagier} introduced the extended double shuffle relations as a refinement of the finite double shuffle relations, and it has been conjectured that the extended double shuffle relations derive all $\Q$-linear relations among MZVs \cite[Conjecture 1]{Ihara-Kaneko-Zagier}. In this subsection, we briefly discuss this framework in a more general setting.

For a $\mathbb{Q}$-algebra homomorphism $Z_R : \ide{h}^{0,\shuffle} \rightarrow R$, we say $Z_R$ satisfies the \textit{finite double shuffle conditions} if $Z_R \circ \mathbf{p} : \Q\langle Y \rangle^{0,*} \rightarrow R$ is a $\Q$-algebra homomorphism.

Let $T$ be an indeterminate and $Z_R : \ide{h}^{0,\shuffle} \rightarrow R$ be a $\Q$-algebra homomorphism satisfying the finite double shuffle conditions.
Since $\ide{h}^{1,\shuffle} = \ide{h}^{0,\shuffle}[x_1]$ and $\Q\langle Y\rangle^{*} = \Q\langle Y\rangle^{0,*}[y_1]$, it follows that there exist two unique $\Q$-algebra homomorphisms associated with $Z_R$.
The first map is $\widetilde{Z}^{\shuffle}_R  : \ide{h}^1 \rightarrow R[T]$, defined by $\widetilde{Z}_R |_{\ide{h}^0} =  Z_R$ and $\widetilde{Z}_R(x_1) := T$. The second map is $ \widetilde{Z}^*_R : \Q\langle Y\rangle^{*} \rightarrow R$, which extends the $\mathbb{Q}$-algebra homomorphism $Z_R \circ \mathbf{p}$ with $\widetilde{Z}^*_R \mid_{\Q\langle Y\rangle^{0,*}} = Z_R \circ \mathbf{p}$ and $\widetilde{Z}^*_R(y_1) := T$.
The main theorem of \cite{Ihara-Kaneko-Zagier} is given as follows.

\begin{prop}[{\cite[Theorem 2]{Ihara-Kaneko-Zagier}}]\label{prop:IKZ-EDS}
Let $(R, Z_R)$ denotes a pair of a $\Q$-algebra $R$ and an element $Z_R$ of $\Hom_{\Qalg}(\ide{h}^{0,\shuffle},R)$ with finite double shuffle conditions. Then the following conditions are equivalent:
\begin{enumerate}
\item 
The following equality holds in $\Hom_{\text{$\Q$-lin}} (\ide{h}^{1}, R)$:
\[
\widetilde{Z}_R = \rho_R \circ \widetilde{Z}^*_R.
\]
Here, the $R$-module map $\rho_R : R[T] \rightarrow R[T]$ is defined by
\[
\rho_R (e^{Tu}) = \exp \left(\sum_{n\ge 2} \frac{(-1)^n}{n} Z_R(x_0^{n-1}x_1)u^n \right).
\]
\item For $w_1\in \ide{h}^{1,\shuffle}$ and $w_0 \in \ide{h}^{0,\shuffle}$, $Z_R^{\shuffle}(w_1\shuffle w_0 - \mathbf{p}(w_1 * w_0)) = 0$ holds.
\end{enumerate}
\end{prop}

Let $\EDS(R)$ be a subset of $\Hom_{\Qalg}(\ide{h}^{\vee}, R)$ satisfying the finite double shuffle condition and Proposition \ref{prop:IKZ-EDS}.

\begin{df}\label{df:EDS}
Let $Z_{R} : \ide{h}^{0,\shuffle} \rightarrow R$ be a $\mathbb{Q}$-algebra homomorphism satisfying the finite double shuffle conditions and Proposition \ref{prop:IKZ-EDS}. The \textit{extended double shuffle relations} are the $\Q$-linear relations derived from
\[
(\widetilde{Z}_{R} - \rho_{R} \circ \widetilde{Z}^*_{R})(w) |_{T = 0} = 0
\]
for all $w\in X^*$.
\end{df}

\begin{ex}\label{ex:reg-MZV}
The most fundamental examples of Proposition \ref{prop:IKZ-EDS} are MZVs.
Let 
\[
Z_{\mathbb{R}} : \ide{h}^{0,\shuffle} \rightarrow \mathbb{R} \quad ; \quad x_1x_0^{k_1-1}\cdots x_1x_0^{k_r-1} \mapsto \zeta(k_1,\ldots,k_r) \text{\hspace{1cm}($k_1,\ldots,k_r\in\mathbb{Z}_{>0}$ with $k_r>1$}).
\]
Since MZVs have the iterated integral expression, $Z_{\mathbb{R}}$ is a $\mathbb{Q}$-algebra homomorphism.
Moreover, $Z_{\mathbb{R}}$ satisfies the finite double shuffle property due to the multiple series representation. 
Therefore, there exist two $\mathbb{Q}$-algebra homomorphisms $\widetilde{Z}_{\mathbb{R}} : \ide{h}^{1,\shuffle} \rightarrow \mathbb{R}$ and $\widetilde{Z}^*_{\mathbb{R}} : \Q\langle Y \rangle^{*} \rightarrow \mathbb{R}$. 
Then, $\widetilde{Z}_{\mathbb{R}}$ and $\widetilde{Z}^*_{\mathbb{R}}$ satisfy the conditions of Proposition \ref{prop:IKZ-EDS}.
In connection with these two maps, we define certain regularized forms of MZVs. For $(k_1,\ldots,k_r)\in\mathbb{Z}_{>0}^r$, we respectively define shuffle-regularized MZVs $\zeta^{\shuffle}$ and harmonic-regularized MZVs $\zeta^*$ by
\begin{align*}
\zeta^{\shuffle}(k_1,\ldots,k_r) :=& \widetilde{Z}_{\mathbb{R}}(x_1x_0^{k_1-1}\cdots x_1x_0^{k_r-1})(0)\\
\zeta^{*}(k_1,\ldots,k_r) :=& \widetilde{Z}^*_{\mathbb{R}}(y_{k_1}\cdots y_{k_r})(0).
\end{align*}
\end{ex}

Let $I_{\EDS}$ be the ideal of $\ide{h}^{\shuffle}$ generated by $x_0$, $x_1$, and $w_1\shuffle w_0 - \mathbf{p}(\mathbf{q}(w_1) * \mathbf{q}(w_0))$ for $w_1\in \ide{h}^{1,\shuffle}$ and $w_0 \in \ide{h}^{0,\shuffle}$, and we define
\[
\mathcal{Z}^f := \ide{h}^{\shuffle}/I_{\EDS}.
\]
For $(k_1,\ldots,k_r)\in \mathbb{Z}_{>0}^r$ with $r\ge1$, we write $\zeta^f(k_1,\ldots,k_r)$ for the image of the word $x_1x_0^{k_1-1}\cdots x_1x_0^{k_r-1}$ in $\mathcal{Z}^f$, and we put $\zeta^f(\emptyset):=1$.
In particular, the classes $\zeta^f(k_1,\ldots,k_r)$ generate $\mathcal{Z}^f$ as a $\Q$-vector space.
In this framework, the authors of \cite{Ihara-Kaneko-Zagier} proposed the following conjecture.

\begin{conj}[{\cite{Ihara-Kaneko-Zagier}}]\label{conj:EDS}
The following $\Q$-algebra homomorphism
\[
\mathcal{Z}^f \rightarrow \mathcal{Z} : \zeta^f(k_1,\ldots,k_r)\mapsto  \zeta(k_1,\ldots,k_r)
\]
is a $\Q$-algebra isomorphism.
Namely, all $\Q$-linear relations among MZVs are obtained from the extended double shuffle relations.
\end{conj}

\subsection{The double shuffle group $\DMR_0(R)$}\label{subsect:DMR}

As mentioned in the previous subsection, Ihara, Kaneko, and Zagier \cite{Ihara-Kaneko-Zagier} introduced the $\shuffle$-regularized MZVs and $*$-regularized MZVs to show the connection between the shuffle product relations and the harmonic product relations among MZVs.
Racinet~\cite{Racinet} derived the regularization theorem (Proposition~\ref{prop:IKZ-EDS}) independently of
\cite{Ihara-Kaneko-Zagier} by using the two coproducts $\Delta_{\shuffle}$ and $\Delta_*$ respectively.
This is regarded as a dual formulation of the approach of \cite{Ihara-Kaneko-Zagier}.
Moreover, he proposed a specific subset $\DMR_0(R)$ and proved that $\DMR_0(R)$ forms a group under an operation
$\circledast$ \cite[Th\'eor\`eme~1]{Racinet}, which is referred to as the \textit{double shuffle group}.



Let
\begin{align*}
\Phi_{\shuffle} :=& \sum_{w\in X} \widetilde{Z}_{\mathbb{R}}(w) \overset{\leftarrow}{w} \in \mathcal{Z} \langle \langle x_0,x_1 \rangle\rangle, \text{ and }\\
\Phi_{*} :=& \sum_{w\in Y} \widetilde{Z}_{\mathbb{R}}^{*}(w) \overset{\leftarrow}{w} \in \mathcal{Z} \langle \langle Y \rangle\rangle,
\end{align*}
where, $\overset{\leftarrow}{w}$ is the reversal word of $w$.
Since $\widetilde{Z}_{\mathbb{R}}$ (respectively, $\widetilde{Z}_{\mathbb{R}}^*$) is a $\Q$-algebra homomorphism with respect to $\shuffle$ (respectively, $*$), and substituting $T = 0$ is also a $\Q$-algebra homomorphism, we have
\begin{align*}
\Delta_{\shuffle}(\Phi_{\shuffle}) =& \Phi_{\shuffle} \otimes \Phi_{\shuffle},\text{ and}\\
 \Delta_*(\Phi_*) =& \Phi_* \otimes \Phi_*.
\end{align*}
Related to these two generating functions, Racinet showed the following:
\begin{thm}[{\cite[Corollary 2.24]{Racinet}}]\label{thm:regularized-Racinet}
Let
\[
\Gamma_{\Phi_\shuffle} :=\exp \left( \sum_{n\ge 2} \frac{(-1)^{n-1}}{n} \langle \Phi_{\shuffle} \mid x_0^{n-1}x_1 \rangle y_1^n \right).
\]
Then, we have
\[
\Phi_* = \Gamma_{\Phi_\shuffle} \mathbf{q}^{\vee}(\Phi_{\shuffle}).
\]
\end{thm}

From Theorem \ref{thm:regularized-Racinet}, Racinet considered a certain subset of $R\langle \langle x_0,x_1 \rangle \rangle$. 
\begin{df}[{\cite[D\'{e}finition 3.2.1]{Racinet}}]\label{def:DMR}
Let $R$ be a $\mathbb{Q}$-algebra.
The double shuffle set $\DMR(R)$\footnote{The notation $\DMR$ is short for the \emph{double mélanges et r\'{e}gularisation}.} is the set of $\Phi \in R  \langle  \langle x_0,x_1 \rangle \rangle $ satisfying the following conditions
\begin{enumerate}
\item\label{def:DMR-1} $\langle \Phi \mid 1\rangle = 1$,
\item\label{def:DMR-2} $\langle \Phi \mid x_0\rangle = \langle \Phi \mid x_1\rangle = 0$,
\item\label{def:DMR-3} $\Delta_{\shuffle}(\Phi) = \Phi\otimes \Phi$, and
\item\label{def:DMR-4} $\Delta_{*}(\Phi_{\star}) = \Phi_{\star} \otimes \Phi_{\star}$.
\end{enumerate}
Here, $\Phi_{\star} := \Gamma_{\Phi}\mathbf{q}^{\vee}(\Phi)$ and $\Gamma_{\Phi}$ is defined by
\[
\Gamma_{\Phi} := \exp \left(\sum_{n\ge 2}\frac{(-1)^{n}}{n}\langle \Phi \mid x_0^{n-1}x_1 \rangle y_1^n\right).
\]
Let $\nu\in R$. By $\DMR_{\nu}(R)$, we denote the subset of $\DMR(R)$ satisfying the additional condition:
\begin{enumerate}
\setcounter{enumi}{4}
\item\label{def:DMR-5} $\langle \Phi \mid x_0x_1\rangle =- \frac{ \nu^2 }{24}$.
\end{enumerate} 
\end{df}

One can obtain an element of $\DMR_0(\C)$ from the KZ-associator (see \cite[Section 3.9.5]{from-numbers-to-motives}).
Let $F_0(z)$ and $F_1(z)$ be the fundamental solutions  of the differential equation
\[
\frac{\del F}{\del z}(z) = \left(\frac{x_0}{z} + \frac{x_1}{z-1} \right) F(z)
\]
around $z = 0$ and  $z = 1$.
Then, the KZ-associator $\Phi_{KZ}$ is the unique constant series such that
\[
F_0(z)=F_1(z)\,\Phi_{KZ}.
\]
More precisely, $\Phi_{KZ}$ is explicitly given by
\[
\Phi_{KZ} = \sum_{w \in X^*} (-1)^{\dep w} \zeta^{\shuffle}(w) \overset{\leftarrow}{w},
\]
where $\dep w$ is the number of $x_1$ in $w$.
Thus, for $\phi \in \C\langle \langle x_0,x_1\rangle\rangle$, we write $\phi(x_0,x_1)$, and we have
\[
\Phi_{KZ} = \Phi_{\shuffle}(x_0,-x_1).
\]
In particular, $\Phi_{KZ}(x_0,-x_1) = \Phi_{\shuffle} \in \DMR(\C)$ (see \cite[Section 3.2.3]{Racinet}).

\begin{rmk}
Let $\Phi \in R \langle \langle x_0,x_1 \rangle \rangle$. Then, $\Phi$ satisfies the condition (\ref{def:DMR-3}) in Definition \ref{def:DMR} if and only if $\langle \Phi \mid u \shuffle v \rangle = \langle \Phi \mid u\rangle \langle \Phi\mid v\rangle$ holds for $u$, $v\in X\*$. Similarly, $\Phi$ satisfies the condition (\ref{def:DMR-4}) in Definition \ref{def:DMR} if and only if $\langle \Phi_* \mid u * v \rangle = \langle \Phi_* \mid u \rangle \langle \Phi_* \mid  v \rangle$ for $u$, $v\in Y^*$. 
This observation implies the existence of a $\mathbb{Q}$-algebra map $Z_{R,\Phi} : \mathcal{Z}^f \rightarrow R ; \zeta^f(k_1,\ldots,k_r) \mapsto \langle \Phi \mid x_0^{k_r-1}x_1\cdots x_0^{k_1-1}x_1 \rangle$ for $\Phi \in \DMR(R)$. It can be verified that $Z_{R,\Phi} \in \EDS(R)$ induces a bijection $\DMR(R) \rightarrow \EDS(R)$. Further details can be found in \cite{Bachmann-Yaddaden}.
\end{rmk}

For $\phi, \psi\in \ide{h}^{\vee}$ with $\langle \phi \mid 1\rangle = \langle \psi \mid 1\rangle = 1$, we define 
\[
\phi \circledast \psi := \phi \cdot \kappa_{\phi^{-1}x_1\phi}(\psi),
\]
where $\kappa_{\phi}$ is an endomorphism of $\ide{h}^{\vee}$ as a $\mathbb{Q}$-algebra, associated with a fixed $f\in \ide{h}^{\vee}$ with $\langle f \mid 1\rangle = 0$, defined by
\[
x_0\mapsto x_0\hspace{1.0ex},\hspace{1.0ex} x_1\mapsto f.
\]

\begin{thm}[{\cite[Th\'{e}or\`{e}m I]{Racinet}}]\label{thm:Racinet-DMR}
For a $\mathbb{Q}$-algebra $R$, $\DMR_0(R)$ is a group with respect to $\circledast$.
\end{thm}
 
The group $\DMR_0(R)$ is often called \textit{double shuffle group}.

\section{Adjoint double shuffle set $\AdDMR_0(R)$} \label{sect:AdDMR}

In this section, we review Jarossay's research \cite{Jarossay}.

\subsection{Adjoint multiple zeta values}

\begin{df}[{\cite[Definition 2.3.1]{Jarossay}}]
For $(k_1,\ldots,k_r) \in \mathbb{Z}_{\ge 0}$ and $l \in\mathbb{Z}_{\ge 0}$, we define $\shuffle$-adjoint MZV ($\shuffle$-AdMZV) by
\begin{align*}
\sAdmzv{l}{k_1,\ldots,k_r} := \sum_{i = 0}^r (-1)^{k_{i+1}+\cdots + k_r + l}\zeta^{\shuffle}(k_1,\ldots,k_i)\zeta_l^{\shuffle}(k_r,\ldots,k_{i+1}) \in \mathcal{Z}.
\end{align*}
Here, 
\[
\zeta_{l}^{\shuffle}(k_r,\ldots,k_{i+1}) := (-1)^l \sum_{\substack{l_{i+1}+\cdots+l_r = l \\ l_{i+1},\ldots,l_r\ge 0}}
\left(\prod_{j = i+1}^r \binom{k_j + l_j-1}{l_j}\right)
\zeta^{\shuffle}(k_r+l_r,\ldots,k_{i+1}+l_{i+1}) \in \mathcal{Z}.
\]
We set $\zeta_{\Ad}^{\shuffle}(\emptyset;0) := 1$ and $\zeta_{\Ad}^{\shuffle}(\emptyset;l) := 0$ for $l \in \mathbb{Z}_{\ge 0}$.
\end{df}

In a similar manner, we define $*$-adjoint MZVs.

\begin{df}[{\cite[Definition 2.3.1]{Jarossay}}]
For $l\in\mathbb{Z}_{\ge 0}$ and $(k_1,\ldots,k_r)\in \mathbb{Z}_{>0}^r$,
we define \emph{$*$-adjoint MZV} ($*$-AdMZV) given by
\begin{align*}
\hAdmzv{l}{k_1,\ldots,k_r} :=& \sum_{i = 0}^r (-1)^{k_{i+1}+\cdots+k_r + l} \zeta^{*}(k_1,\ldots,k_i)\zeta^{*}_l(k_r,\ldots,k_{i+1}) \in \mathcal{Z}.
\end{align*}
Here, 
\[
\zeta_{l}^{*}(k_r,\ldots,k_{i+1}) := (-1)^l \sum_{\substack{l_{i+1}+\cdots+l_r = l \\ l_{i+1},\ldots,l_r\ge 0}}
\left(\prod_{j = i+1}^r \binom{k_j + l_j-1}{l_j}\right)
\zeta^{*}(k_r+l_r,\ldots,k_{i+1}+l_{i+1}).
\]
\end{df}
It should be noted that these two types of adjoint MZVs are congruent modulo $\zeta(2)$ (\cite[Proposition 3.2.4]{Jarossay}), i.e.,
\begin{align}
\sAdmzv{l}{k_1,\ldots,k_r} \equiv \hAdmzv{l}{k_1,\ldots,k_r} \mod \zeta(2)\mathcal{Z}. \label{eq:gen-AdMZV}
\end{align}
\begin{df}[{\cite[Definition 2.3.1]{Jarossay}}]\label{df:AdMZV}
For $l\in\mathbb{Z}_{\ge 0}$ and $(k_1,\ldots,k_r) \in \mathbb{Z}_{>0}^r$, we define an adjoint multiple zeta value as an element of $\mathcal{Z}/\zeta(2)\mathcal{Z}$ given by
\begin{align*}
\Admzv{l}{k_1,\ldots,k_r} := \sum_{i = 0}^r (-1)^{k_{i+1}+\cdots + k_r + l}\zeta^{\shuffle}(k_1,\ldots,k_i)\zeta_l^{\shuffle}(k_r,\ldots,k_{i+1}) \in \mathcal{Z}/\zeta(2)\mathcal{Z},
\end{align*}
where
\[
\zeta_{l}^{\shuffle}(k_r,\ldots,k_{i+1}) := (-1)^l \sum_{\substack{l_{i+1}+\cdots+l_r = l \\ l_{i+1},\ldots,l_r\ge 0}}
\left(\prod_{j = i+1}^r \binom{k_j + l_j-1}{l_j}\right)
\zeta^{\shuffle}(k_r+l_r,\ldots,k_{i+1}+l_{i+1}).
\]
\end{df}

Let $\iota:\mathcal Z \twoheadrightarrow \mathcal Z/\zeta(2)\mathcal Z$ be the canonical projection.
It induces a $\Q$-algebra homomorphism (by coefficientwise reduction)
\[
\iota_*:\mathcal Z\langle\!\langle x_0,x_1\rangle\!\rangle
\longrightarrow
(\mathcal Z/\zeta(2)\mathcal Z)\langle\!\langle x_0,x_1\rangle\!\rangle .
\]
Define
\[
\overline{\Phi}_{\shuffle}:=\iota_*(\Phi_{\shuffle})\in (\mathcal Z/\zeta(2)\mathcal Z)\langle\!\langle x_0,x_1\rangle\!\rangle .
\]

The fundamental idea of Jarossay's work is that the generating function of AdMZVs is given by
\begin{align*}
\overline{\Phi}_{\shuffle}^{-1}x_1\overline{\Phi}_{\shuffle} \in (\mathcal{Z}/\zeta(2)\mathcal{Z})\langle \langle x_0,x_1 \rangle\rangle,
\end{align*}
Based on this perspective, he established some properties of $\Phi_{\shuffle}^{-1}x_1\Phi_{\shuffle}$ and derived the $\Q$-linear relations among AdMZVs. He named these relations \textit{adjoint double shuffle relations}.
Additionally, he introduced a set $\AdDMR_0(R)$ whose defining equations are adjoint double shuffle relations, and raised the question \textit{whether the sets, $\DMR_0(R)$ and $\AdDMR_0(R)$, are actually bijective or not}.
In this section, we briefly review how the generating series $\overline{\Phi}_{\shuffle}^{-1}x_1\overline{\Phi}_{\shuffle}$ leads to the adjoint double shuffle relations, and we define Jarossay's $\AdDMR_0(R)$.
 
\subsection{The generating function of AdMZVs}

A key idea in Jarossay's work is that $\shuffle$-AdMZVs appear as coefficients of $\Phi_{\shuffle}^{-1} x_1 \Phi_{\shuffle}$.
The following proposition formalizes this observation.
\begin{prop}[{\cite[Proposition 3.2.2]{Jarossay}}]
For $l\in \mathbb{Z}_{\ge 0}$ and $k_1,\ldots,k_r\in\mathbb{Z}_{>0}$, the following holds:
\[
\langle \Phi_{\shuffle}^{-1}x_1\Phi_{\shuffle} \mid x_0^lx_1x_0^{k_r-1}x_1\cdots x_0^{k_1-1}x_1\rangle = \sAdmzv{l}{k_1,\ldots,k_r} .
\]
\end{prop}
\begin{proof}
For $l\in\mathbb{Z}_{\ge 0}$ and $k_1,\ldots,k_r\in\mathbb{Z}_{>0}$,
the claim follows from
{\small
\begin{align*}
\langle \Phi_{\shuffle}^{-1}x_1\Phi_{\shuffle} \mid x_0^{k_r-1}x_1\cdots x_0^{k_1-1}x_1x_0^l \rangle
=& \sum_{i = 0}^r \langle\Phi_{\shuffle}^{-1} \mid x_0^lx_1x_0^{k_r-1}x_1\cdots x_0^{k_i-1}\rangle \langle \Phi_{\shuffle}\mid x_0^{k_{i+1}-1}x_1 \cdots x_0^{k_r-1}x_1\rangle \\
=& \sum_{i = 0}^r(-1)^{k_1 + \cdots + k_i + l}\langle \Phi_{\shuffle} \mid x_0^{k_i-1}x_1 \cdots x_0^{k_1-1}x_1x_0^l \rangle\\
&\qquad\qquad \qquad\qquad \qquad\qquad   \times \langle \Phi_{\shuffle}\mid x_0^{k_{i+1}-1}x_1 \cdots x_0^{k_r-1}x_1\rangle \displaybreak[3] \\
=&\sum_{i = 0}^r(-1)^{k_1 + \cdots + k_i + l} \zeta_l^{\shuffle}(k_1,\ldots,k_i)\zeta^{\shuffle}(k_r,\ldots,k_{i+1})\\
=&\sAdmzv{l}{k_r,\ldots,k_1}.
\end{align*}
}
Here, in the second equality, we use the fact that $\Phi_{\shuffle}$ is group-like, which implies $S_X^{\vee}(\Phi_{\shuffle}) = \Phi_{\shuffle}^{-1}$.
Also, in the third equality, we use
\begin{align*}
\langle \Phi_{\shuffle} \mid x_0^{k_r-1}x_1\cdots x_0^{k_1-1}x_1x_0^l \rangle =  \zeta_{l}^{\shuffle} (k_1,\ldots,k_r).
\end{align*}
\end{proof}

\begin{df}
\begin{enumerate}
\item We define $\inv$ as the anti-automorphism of $\Q\langle\!\langle Y \rangle\!\rangle$ defined by $y_k \mapsto (-1)^k y_k$ for $k\in\mathbb{Z}_{>0}$.
\item
For $l \in \mathbb{Z}_{\ge 0}$, we define $\sft_{*.l}$ as the $\Q$-linear endomorphism of $\Q\langle\langle Y \rangle\rangle$ given by
{\small
\begin{align}
&\sft_{*,l}(\Phi) \notag \\
:=&\delta_{l,0}y_1 + (-1)^l \sum_{\substack{y_{k_1}\cdots y_{k_r} \in Y^* \\ k_1,\ldots,k_r\ge 1, r\ge 1}}\sum_{l_1+\cdots +l_r = l}\left(\prod_{i = 1}^r \binom{k_i + l_i -1}{l_i} \right)\langle\Phi \mid y_{k_1+l_1}\cdots y_{k_r+l_r} \rangle y_{l+1}y_{k_1}\cdots y_{k_r}. \label{eq:sft-l}
\end{align}
}
for $\Phi \in \Q\langle \langle Y \rangle\rangle$.
Here, $\delta_{l,0}$ is the Kronecker delta.
\end{enumerate}
We use the same symbol $\inv$ and $\sft_{*,l}$ after extension of scalars.
\end{df}
\begin{prop}\label{prop:Ad-har-exp}
We have
\begin{align*}
 \left(\sum_{l\ge 0} \inv\circ \sft_{*,l} (\Phi_{\star}) \right)\cdot \Phi_{\star} = y_1 + \sum_{\substack{y_{l+1}y_{k_r}\cdots y_{k_1} \in Y^* \\ k_1,\ldots,k_r\ge 1,  r\ge 1}}\hAdmzv{k_1,\ldots,k_r}{l}y_{l+1}y_{k_r}\cdots y_{k_1},
\end{align*}
where we recall 
\[
\Phi_{\star} := \sum_{w\in Y} \widetilde{Z}_{\mathbb{R}}^{*}(w) \overset{\leftarrow}{w} \in \mathcal{Z} \langle \langle Y \rangle\rangle.
\] 
\end{prop}

\begin{proof}

For $l\in \mathbb{Z}_{\ge 0}$ and $k_1,\ldots,k_r \in \mathbb{Z}_{>0}$, we have
{\small
\begin{align*}
&\left\langle \left(\sum_{l\ge 0} \inv\circ \sft_{*,l} (\Phi_{\star}) \right)\cdot \Phi_{\star} \middle| y_{l+1}y_{k_r}\cdots y_{k_1}  \right\rangle\\
=&\sum_{i = 0}^r \left\langle \sum_{l\ge 0} \inv\circ \sft_{*,l} (\Phi_{\star}) \middle| y_{l+1}y_{k_r}\cdots y_{k_{i+1}}  \right\rangle \cdot \left\langle\Phi_{\star} \middle| y_{k_i}\cdots y_{k_{1}}  \right\rangle\\
\stackrel{\eqref{eq:sft-l}}{=}&\hspace{0cm}\sum_{i = 0}^r \left((-1)^l\sum_{l_r+\cdots +l_{i+1} = l} \left(\prod_{j = 1}^r \binom{k_j + l_j -1}{l_j} \right) \left\langle \inv(\Phi_{\star}) \mid y_{k_r+l_r}\cdots y_{k_{i+1} + l_{i+1}}  \right\rangle \right)\cdot \left\langle\Phi_{\star} \middle| y_{k_i}\cdots y_{k_{1}}  \right\rangle\\
=&\sum_{i = 0}^r \left((-1)^{k_{i+1}+ \cdots + k_r}\sum_{l_r+\cdots +l_{i+1} = l} \left(\prod_{j = 1}^r \binom{k_j + l_j -1}{l_j} \right) \left\langle \Phi_{\star} \mid y_{k_r+l_r}\cdots y_{k_{i+1} + l_{i+1}}  \right\rangle \right)\cdot \left\langle\Phi_{\star} \middle| y_{k_i}\cdots y_{k_{1}}  \right\rangle\\
=& \sum_{i = 0}^r (-1)^{k_{i+1}+ \cdots + k_r}\sum_{l_r+\cdots +l_{i+1} = l} \left(\prod_{j = 1}^r \binom{k_j + l_j -1}{l_j} \right) \zeta^*_l(k_r,\ldots,k_{i+1})\zeta^*(k_i.\ldots.k_1),
\end{align*}
}
as claimed.
\end{proof}

\subsection{Adjoint double shuffle relations and $\AdDMR_0(R)$}

One of the purposes of this section is to prove Proposition~\ref{prop:t-adic}.
We first prepare some tools.
Let $T$ be an indeterminate and
we put
\begin{align*}
\Q[T;Y\rangle := \Q[T]\langle Y \rangle.
\end{align*}
Let
\begin{align*}
\mathbf{q}_{\#}\colon \ide{h}\longrightarrow\Q[T;Y\rangle
\end{align*}
be the $\Q$-linear map defined by
\[
\mathbf{q}_{\#}\bigl(x_0^{k_1-1}x_1 x_0^{k_2-1}x_1 \cdots x_0^{k_r-1}x_1 x_0^{k_{r+1}-1}\bigr)
:=
\begin{cases}
T^{k_1-1} y_{k_2}\cdots y_{k_r} & \text{if } k_{r+1}=1,\\
0 & \text{otherwise.}
\end{cases}
\]
Recall that $\Q\langle Y\rangle$ is endowed with the weight filtration, and $Y_{\ge m}\subset \Q\langle Y\rangle$ is the $\Q$-subspace spanned by words of weight $\ge m$.
We endow $\Q[T;Y\rangle$ with the grading determined by $\deg T=1$ and $\deg y_m=m$ for $m\in \mathbb{Z}_{\ge 1}$.
For $d\ge0$, let $\Q[T;Y\rangle^{(d)}$ denote the homogeneous component of total degree $d$.
For $m\ge1$, set
\[
A_m :=\Q[T;Y\rangle \Big/ \bigoplus_{d\ge m} \Q[T;Y\rangle ^{(d)} .
\]
Then the natural truncation maps $A_{m'}\rightarrow A_m$ $(m\le m')$ form a projective system, and we define
\begin{equation*}\label{eq:lim-com}
\Q[[T;Y \rangle\rangle
:= \varprojlim_{m\ge1} A_m .
\end{equation*}
Note that  canonical inclusion $\Q[T] \hookrightarrow \Q[T; Y \rangle$ induces a canonical action $\Q[T]/(T)^m \rightarrow  A_m$ for each $m \in \mathbb{Z}_{\ge 1}$.
Passing to the projective limit, we naturally see that $\Q[[T; Y \rangle\rangle$ is a $\Q[[T]]$-algebra

In the same manner, we define $R$-algebra $R[[T;Y \rangle\rangle$ as
\[
R[[T;Y \rangle\rangle
\cong \varprojlim_{m} R\otimes A_m
\]
and we can see $R[[T;Y \rangle\rangle$ as $R[[T]]$-algebra.

By definition, for each $n\ge 1$ we have
\[
\mathbf{q}_{\#}(X^{n})
\subset
\Q[T;Y\rangle^{(n-1)}
\subset\ \Q[T; Y\rangle .
\]
Thus, $\mathbf{q}_{\#}$ is a continuous $\Q$-linear map on $\ide{h}$ endowed with the $X$-adic topology,
and hence it extends uniquely to a continuous $\Q$-linear map
\[
\mathbf{q}_{\#}^{\vee}\colon \ide{h}^{\vee}\longrightarrow \Q[[T;Y \rangle\rangle .
\]
We use the same symbol $\mathbf{q}_{\#}^{\vee}$ via the scalar expansions.

\begin{prop}[{\cite[Proposition 3.2.5]{Jarossay}}]\label{prop:t-adic}
Let $R$ be a unital commutative associative $\Q$-algebra. For $\Phi \in \DMR_0(R)$, we have
\begin{enumerate}
\item\label{thm:t-adic-1} $\Delta_{\shuffle}(\Phi^{-1} x_1 \Phi)=\Phi^{-1} x_1 \Phi \otimes 1 + 1\otimes \Phi^{-1} x_1 \Phi$,
\item\label{thm:t-adic-2} $\Delta_{*}\bigl(\mathbf{q}_{\#}^{\vee}(\Phi^{-1} x_1 \Phi)\bigr)
=\mathbf{q}_{\#}^{\vee}(\Phi^{-1} x_1 \Phi)\otimes \mathbf{q}_{\#}^{\vee}(\Phi^{-1} x_1 \Phi)$,
where we view $\Delta_*$ as an $R[[T]]$-algebra homomorphism via the scalar extension to $R[[T]]$.
\end{enumerate}
\end{prop}

If we take $\overline{\Phi}_{\shuffle} \in \DMR_0(\mathcal{Z}/\zeta(2)\mathcal{Z})$, we obtain a class of $\Q$-algebraic relations among AdMZVs, referred to as the \emph{adjoint double shuffle relations}.
We shall prove Proposition \ref{prop:t-adic} through Propositions \ref{prop:ad-proj}--\ref{prop:Ad-har}.

Let $\Phi \in \DMR_0(R)$.
First, we prove Proposition \ref{prop:t-adic}\,(\ref{thm:t-adic-1}).

\begin{prop}[{ = Proposition \ref{prop:t-adic} (\ref{thm:t-adic-1}})]
For  $\Phi \in \DMR_0(R)$, we have
\begin{align}
\Delta_{\shuffle}(\Phi^{-1}x_1\Phi ) = \Phi^{-1}x_1\Phi  \otimes 1 + 1\otimes \Phi^{-1}x_1\Phi. \label{eq:ad-shuffle}
\end{align}

\end{prop}

\begin{proof}

Since $\Delta_{\shuffle}(\Phi) = \Phi \otimes \Phi$, we have
\begin{align*}
\Delta_{\shuffle}(\Phi^{-1}x_1\Phi ) =& \Delta_{\shuffle}(\Phi^{-1})\Delta_{\shuffle}(x_1) \Delta_{\shuffle}(\Phi )  \\
=&(\Phi^{-1} \otimes \Phi^{-1} ) \times (x_1\otimes 1 + 1\otimes x_1) \times (\Phi \otimes \Phi ) \\
=& \Phi^{-1}x_1\Phi  \otimes 1 + 1\otimes \Phi^{-1}x_1\Phi.
\end{align*}
\end{proof}
The group-like condition $\Delta_{\shuffle}(\Phi)=\Phi\otimes\Phi$ is equivalent to the shuffle product relations for $\shuffle$-regularized MZVs.
Similarly, \eqref{eq:ad-shuffle} gives rise to certain $\Q$-algebraic relations among AdMZVs coming from the shuffle product.
We call the $\Q$-algebraic relations among AdMZVs arising from \eqref{eq:ad-shuffle} the \textit{adjoint shuffle relations}.

We next investigate the $\Q$-algebraic relations among AdMZVs arising from the harmonic product relations.

\begin{prop}\label{prop:ad-proj}
For  $\Phi \in \DMR_0(R)$, we have
\[
\mathbf{q}^{\vee}(\Phi^{-1}x_1 \Phi) = \left(\sum_{l\ge 0}  \sft_{*,l} \circ \inv (\Phi_{\star}) \right) \cdot \Phi_{\star}.
\]
\end{prop}

\begin{proof}
To begin with, we show that, for $l\ge 0$,
\begin{align}
\sft_{*,l} \circ \inv (\Gamma_{\Phi}\mathbf{q}^{\vee}(\Phi))
= \sft_{*,l} \circ \inv(\mathbf{q}^{\vee}(\Phi)) \cdot \inv( \Gamma_{\Phi}) . \label{eq:sft-Gamma}
\end{align}
We recall the definition of $\sft_{*,l}$:
{\small
\begin{align}
\sft_{*,l}(\Psi) :=\delta_{l,0}y_1 + (-1)^l \sum_{\substack{y_{k_1}\cdots y_{k_r} \in Y^* \\ k_1,\ldots,k_r\ge 1, r\ge 1}}\sum_{l_1+\cdots +l_r = l}\left(\prod_{i = 1}^r \binom{k_i + l_i -1}{l_i} \right)\langle\Psi \mid y_{k_1+l_1}\cdots y_{k_r+l_r} \rangle y_{l+1}y_{k_1}\cdots y_{k_r} \label{eq:sft-l-2}
\end{align}
}
for $\Psi \in R\langle \langle Y \rangle\rangle$. For $l\ge 0$ and $w\in Y^*$, if we specialize to $\Psi=wy_1$, then the summand in \eqref{eq:sft-l-2} becomes
\begin{align*}
&\sum_{\substack{ k_1,\ldots,k_r\ge 1, r\ge 1}}\sum_{\substack{l_1+\cdots+l_r = l \\ l_1,\ldots,l_r\ge 0}} \left( \prod_{i = 1}^r \binom{k_i + l_i-1}{l_i}  \right)
\langle wy_1 \mid y_{k_1+l_1}\cdots y_{k_r+l_r} \rangle\\
=&\sum_{\substack{k_1,\ldots,k_{r-1}\ge 1, r\ge 1}}\sum_{\substack{l_1+\cdots+l_{r-1} = l \\ l_1,\ldots,l_{r-1}\ge 0}} \left( \prod_{i = 1}^{r-1} \binom{k_i + l_i-1}{l_i}  \right)
\langle wy_1 \mid  y_{k_1+l_1}\cdots y_{k_{r-1}+l_{r-1}}y_1 \rangle .
\end{align*}
It follows that 
\begin{align}
\sft_{*,l}(\Psi y_1) = \sft_{*,l}(\Psi)y_1 \label{eq:sft-y1}
\end{align}
 for $\Psi \in R\langle\langle Y \rangle\rangle$. 
Since the only positions at which the coefficients of $\Gamma_{\Phi}$ may be non-zero are the powers of $y_1$, we have
\begin{align*}
&\sft_{*,l} \circ \inv (\Gamma_{\Phi}\mathbf{q}^{\vee}(\Phi)) \displaybreak[3] \\
=&\sft_{*,l}\left(  \inv \circ \mathbf{q}^{\vee}(\Phi) \cdot \exp\left(-\sum_{n\ge 2} \frac{(-1)^{n-1}}{n}\langle \Phi \mid x_0^{n-1}x_1 \rangle \inv(y_1^n) \right)  \right) \displaybreak[3] \\
=&\sft_{*,l}\left(  \inv \circ \mathbf{q}^{\vee}(\Phi)  \cdot \sum_{k\ge 0} \frac{1}{k!}\left(-\sum_{n\ge 2} \frac{(-1)^{n-1}}{n}\langle \Phi \mid x_0^{n-1}x_1 \rangle \inv(y_1^n) \right)^k  \right) \displaybreak[3] \\
=&\sum_{k\ge 0} \frac{1}{k!} \-\sum_{n\ge 2} \sft_{*,l}\left( \inv \circ \mathbf{q}^{\vee}(\Phi)  \cdot  \left( \frac{(-1)^{n}}{n}\langle \Phi \mid x_0^{n-1}x_1 \rangle \inv(y_1^k) \right)^k \right) \displaybreak[3] \\
\stackrel{\eqref{eq:sft-y1}}{=}& \sft_{*,l} \inv \circ( \mathbf{q}^{\vee}(\Phi) )  \cdot \sum_{k\ge 0} \frac{1}{k!} \-\sum_{n\ge 2} \left( \frac{(-1)^{n}}{n}\langle \Phi \mid x_0^{n-1}x_1 \rangle \inv(y_1^n) \right)^k \displaybreak[3] \\
=& \sft_{*,l} \circ \inv( \mathbf{q}^{\vee}(\Phi)) \cdot \inv(\Gamma_{\Phi}). \displaybreak[3]
\end{align*}
and complete the proof of \eqref{eq:sft-Gamma}.
By using \eqref{eq:sft-Gamma}, we have
\begin{align*}
& \left(\sum_{l\ge 0}  \sft_{*,l} \circ \inv (\Phi_{\star}) \right) \cdot \Phi_{\star} \displaybreak[3] \\
=& \left(\sum_{l\ge 0} \sft_{*,l} \circ \inv  \left(\Gamma_{\Phi} \cdot (\mathbf{q}^{\vee}(\Phi))\right) \right) \cdot \Gamma_{\Phi} \cdot\mathbf{q}^{\vee}(\Phi) \displaybreak[3] \\
\stackrel{\eqref{eq:sft-Gamma}}{=}&\left(\sum_{l\ge 0}  \sft_{*,l} \circ \inv (\mathbf{q}^{\vee}(\Phi)) \cdot \inv(\Gamma_{\Phi} ) \right) \cdot \Gamma_{\Phi} \cdot\mathbf{q}^{\vee}(\Phi) \displaybreak[3] \\
=&\left(\sum_{l\ge 0} \sft_{*,l} \circ  \inv (\mathbf{q}^{\vee}(\Phi)) \right) \cdot  \inv(\Gamma_{\Phi}) \cdot \Gamma_{\Phi} \cdot \mathbf{q}^{\vee}(\Phi) 
\end{align*}
Here, we compute
\begin{align*}
&\inv\left(\Gamma_{\Phi}\right) \cdot \Gamma_{\Phi} \\
=&\exp\left(- \sum_{n\ge 2} \frac{(-1)^{n-1}}{n} \langle \Phi\mid x_0^{n-1}x_1 \rangle \inv(y_1^n) \right) \cdot \exp\left(- \sum_{n\ge 2} \frac{(-1)^{n-1}}{n} \langle \Phi\mid x_0^{n-1}x_1 \rangle y_1^n \right) \\
=&\exp  \left(- \sum_{n\ge 2} \frac{-1 + (-1)^{n-1}}{n} \langle \Phi\mid x_0^{n-1}x_1 \rangle y_1^n \right)\\
=&\exp  \left( \sum_{n\ge 2} \frac{2}{n} \langle \Phi\mid x_0^{n-1}x_1 \rangle y_1^n \right).
\end{align*}
Since $\Phi\in \DMR_0(R)$, we have $\langle \Phi \mid x_0^{n-1}x_1 \rangle = 0$ for arbitary positive even integer $n$ \cite[Proposition 3.3.3]{Racinet}.
Thus we have
\begin{align*}
\inv\left(\Gamma_{\Phi}\right) \cdot \Gamma_{\Phi} = 1.
\end{align*}
Therefore, it follows that
\begin{align*}
&\left(\sum_{l\ge 0}  \sft_{*,l} \circ \inv  (\Phi_{\star}) \right) \cdot \Phi_{\star} \\
=&\left(\sum_{l\ge 0}  \sft_{*,l} \circ \inv (\mathbf{q}(\Phi)) \right) \cdot \mathbf{q}(\Phi).
\end{align*}
Next, for $y_{l+1}y_{k_r}\cdots y_{k_1}\in Y^*$, we have
{\small
\begin{align*}
\bigl\langle \mathbf{q}^{\vee}(\Phi^{-1}x_1\Phi)\,\big|\, y_{l+1}y_{k_r}\cdots y_{k_1}\bigr\rangle
&=
\bigl\langle \Phi^{-1}x_1\Phi \,\big|\, x_0^{l}x_1x_0^{k_r-1}x_1\cdots x_0^{k_1-1}x_1\bigr\rangle \\
&=
\sum_{i=0}^{r}
\bigl\langle \Phi^{-1}\,\big|\, x_0^{l}x_1x_0^{k_r-1}x_1\cdots x_0^{k_{i+1}-1}\bigr\rangle\,
\bigl\langle \Phi\,\big|\, x_0^{k_i-1}x_1\cdots x_0^{k_1-1}x_1\bigr\rangle \\
&=
\sum_{i=0}^{r}(-1)^{k_{i+1}+\cdots+k_r+l}\,
\bigl\langle \Phi\,\big|\, x_0^{k_{i+1}-1}x_1\cdots x_0^{k_r-1}x_1x_0^{l}\bigr\rangle\,
\bigl\langle \Phi\,\big|\, x_0^{k_i-1}x_1\cdots x_0^{k_1-1}x_1\bigr\rangle .
\end{align*}
}
In the second equality, we note that $\langle \Phi^{-1} \mid x_0^l\rangle = \frac{1}{l!} \langle \Phi^{-1} \mid x_0^{\shuffle l}\rangle = 0$.
In the third equality, we use the antipode $S_X^{\vee}$.
Since $\Delta_{\shuffle}(\Phi)=\Phi\otimes\Phi$, for any $u\in X^*\setminus\{1\}$, we have
\begin{align}
\langle \Phi \mid x_0\shuffle u\rangle
=
\langle\Phi \mid x_0\rangle\cdot \langle\Phi \mid u\rangle
=
0. \label{eq:x0-shuffle}
\end{align}
By the general properties of the shuffle product, $ux_1x_0^l \in X^*$ can be written as
\begin{align}
ux_1x_0^l = \sum_{i = 0}^l u_i \shuffle x_0^{i} \label{eq:u-expl}
\end{align}
for some $u_i \in \ide{h}$, and $u_0$ is explicitly given by 
\[
u_0 = (-1)^l (u\shuffle x_0^{l})\cdot x_1.
\]
By  combining \eqref{eq:x0-shuffle} and \eqref{eq:u-expl}, we have
\begin{align*}
\langle \Phi \mid  x_0^{k_r-1}x_1\cdots x_0^{k_{i+1}-1}x_1x_0^l \rangle =& (-1)^l \langle \Phi \mid   ( x_0^{k_r-1}x_1\cdots x_0^{k_{i+1}-1}x_1 \shuffle x_0^{l})x_1 \rangle.
\end{align*}
More explicitly, we have
{\small
\begin{align*}
\bigl\langle\Phi \,\big|\, x_0^{k_r-1}x_1\cdots x_0^{k_{i+1}-1}x_1x_0^l  \bigr\rangle
=
(-1)^{l}\!\!\sum_{\substack{l_{i+1}+\cdots+l_r=l\\ l_{i+1},\ldots,l_r\ge0}}
\left(\prod_{j=i+1}^{r}\binom{k_j+l_j-1}{l_j}\right)
\bigl\langle \Phi \,\big|\, x_0^{k_{r}+l_{r}-1}x_1\cdots x_0^{k_{i+1}+l_{i+1}-1}x_1\bigr\rangle .
\end{align*}
}
Hence, we obtain
{\small
\begin{align*}
&\bigl\langle \mathbf{q}^{\vee}(\Phi^{-1}x_1\Phi)\,\big|\, y_{l+1}y_{k_r}\cdots y_{k_1}\bigr\rangle\\
&=
\sum_{i=0}^{r}
(-1)^{k_{i+1}+\cdots+k_r}
\left(
\sum_{\substack{l_1+\cdots+l_r=l\\ l_1,\ldots,l_r\ge0}}
\left(\prod_{j=i+1}^{r}\binom{k_j+l_j-1}{l_j}\right)
\bigl\langle \Phi \,\big|\, x_0^{k_{i+1}+l_{i+1}-1}x_1\cdots x_0^{k_r+l_r-1}x_1\bigr\rangle
\right) \\
&\hspace{9cm}\cdot
\bigl\langle \Phi\,\big|\, x_0^{k_i-1}x_1\cdots x_0^{k_1-1}x_1\bigr\rangle \\
&=
\sum_{i=0}^{r}
(-1)^{k_{i+1}+\cdots+k_r}
\left(
\sum_{\substack{l_1+\cdots+l_r=l\\ l_1,\ldots,l_r\ge0}}
\left(\prod_{j=i+1}^{r}\binom{k_j+l_j-1}{l_j}\right)
\bigl\langle \mathbf{q}^{\vee}(\Phi) \,\big|\,   y_{k_{i+1} + l_{i+1}}\cdots y_{k_r + l_r} \bigr\rangle
\right) \\
&\hspace{9cm}\cdot
\bigl\langle \mathbf{q}^{\vee}(\Phi)\,\big|\,  y_{k_i}\cdots y_{k_1}\bigr\rangle \\
&=
\sum_{i=0}^{r}
(-1)^l
\left(
\sum_{\substack{l_1+\cdots+l_r=l\\ l_1,\ldots,l_r\ge0}}
\left(\prod_{j=i+1}^{r}\binom{k_j+l_j-1}{l_j}\right)
\bigl\langle \inv \circ \mathbf{q}^{\vee}(\Phi) \,\big|\,   y_{k_r + l_r} \cdots y_{k_{i+1} + l_{i+1}} \bigr\rangle
\right) \\
&\hspace{9cm}\cdot
\bigl\langle \mathbf{q}^{\vee}(\Phi)\,\big|\,  y_{k_i}\cdots y_{k_1}\bigr\rangle \\
&=
\bigl\langle  \sft_{*,l}\circ  \inv \circ \mathbf{q}^{\vee}(\Phi)\,\big|\, y_{l+1}y_{k_r}\cdots y_{k_{i+1}}\bigr\rangle\,
\bigl\langle \mathbf{q}^{\vee}(\Phi)\,\big|\, y_{k_i}\cdots y_{k_1}\bigr\rangle .
\end{align*}
}
This proves the desired identity.
\end{proof}

\begin{rmk}
If we apply Proposition~\ref{prop:ad-proj} to
\[
\Phi =\overline{\Phi}_{\shuffle} \in \DMR_0(\mathcal{Z}/\zeta(2)\mathcal{Z}),
\]
Proposition \ref{prop:Ad-har-exp} implies \eqref{eq:gen-AdMZV}.
\end{rmk}

\begin{df}
We define the $\Q[[T]]$-algebra endomorphism $\sft_*^{\vee}$ of $\Q[[T;Y \rangle\rangle$ by setting, for each $k\ge 1$,
\[
y_k \longmapsto \sum_{j=0}^{k-1} \binom{k-1}{j} (-T)^{j} y_{k-j}.
\]
we use the same symbol $\sft_{*}^{\vee}$ after extension of scalars.
\end{df}
\begin{prop}[{In the proof of \cite[Proposition 3.2.5(ii)]{Jarossay}}]\label{eq:sft-coalg}
The algebra homomorphism $\sft_*^{\vee}$ is the coalgebra homomorphism with respect to $\Delta_*$, that is, we have
\begin{align*}
\sft_*^{\vee,\otimes 2}\circ \Delta_* = \Delta_*\circ \sft_*^{\vee}.
\end{align*}
\end{prop}

\begin{proof}

Since $\Delta_*$ is the algebra homomorphism with respect to the concatenation product and $\sft_*^{\vee}$ is also so, it is enough to verify the identity on the generators $y_k$ $(k\ge 1)$, namely
\begin{align*}
\sft_*^{\vee,\otimes 2}(\Delta_*(y_k)) = \Delta_*(\sft_*^{\vee}(y_k)).
\end{align*}
We have
\begin{align*}
&\sft_*^{\vee ,\otimes 2} (\Delta_*(y_k)) \displaybreak[3] \\
=&\sum_{s = 0}^{k-1} \binom{k-1}{s} (-T)^s (y_{k-s} \otimes 1 + 1\otimes y_{k-s})  \\
&\quad + \sum_{\substack{i+j = k \\ i,j>0}}\left(\sum_{\substack{0 \le i_1 \le i-1 \\ 0\le j_1 \le j-1}} (-T)^{i_1 + j_1} \binom{i-1}{i_1}\binom{j-1}{j_1}\cdot (y_{i-i_1} \otimes y_{j-j_1})  \right). \displaybreak[3]
\end{align*}
For the second term, we have
\begin{align*}
& \sum_{\substack{i+j = k \\ i,j>0}}\left(\sum_{\substack{0 \le i_1 \le i-1 \\ 0\le j_1 \le j-1}} (-T)^{i_1 + j_1} \binom{i-1}{i_1}\binom{j-1}{j_1}\cdot (y_{i-i_1} \otimes y_{j-j_1}) \right)\displaybreak[3]  \\
=&\sum_{s = 2}^k \left(\sum_{\substack{s_1 + s_2 = s \\ s_1,s_2\ge 1}}(-T)^{k-s} \left(\sum_{\substack{t_1 + t_2 = k-s \\ t_1,t_2\ge 0}}\binom{s_1 + t_1-1}{t_1}\binom{s_2 + t_2-1}{t_2}    \right)\right) \cdot  (y_{s_1} \otimes y_{s_2}). \displaybreak[3]
\end{align*}
Since, by the Chu--Vandermonde identity and fact $\displaystyle \binom{s_i+t_i-1}{t_i}=(-1)^{t_i}\binom{-s_i}{t_i}$, we have
\begin{align*}
&\sum_{\substack{t_1 + t_2 = k-s \\ t_1,t_2\ge 0}}\binom{s_1 + t_1-1}{t_1}\binom{s_2 + t_2-1}{t_2} \displaybreak[3] \\
=&\sum_{\substack{t_1 + t_2 = k-s \\ t_1,t_2\ge 0}} (-1)^{t_1 + t_2} \binom{-s_1}{t_1}\binom{-s_2}{t_2} \displaybreak[3] \\
=&(-1)^{k-s} \binom{-s}{k-s} \displaybreak[3] \\
=&\binom{k-1}{k-s} \displaybreak[3] \\
=&\binom{k-1}{s-1}, \displaybreak[3]
\end{align*}
The second term becomes
\begin{align*}
& \sum_{\substack{i+j = k \\ i,j>0}}\left(\sum_{\substack{0 \le i_1 \le i-1 \\ 0\le j_1 \le j-1}} (-T)^{i_1 + j_1} \binom{i-1}{i_1}\binom{j-1}{j_1}\cdot (y_{i-i_1} \otimes y_{j-j_1}) \right) \displaybreak[3] \\
=&\sum_{s = 2}^k \sum_{\substack{s_1 + s_2 = s \\ s_1,s_2\ge 1}}(-T)^{k-s} \binom{k-1}{s-1} \cdot  (y_{s_1} \otimes y_{s_2}).
\end{align*}
Therefore,
\begin{align*}
&\sft_*^{\vee ,\otimes 2} (\Delta_*(y_k)) \displaybreak[3] \\
=&\sum_{s = 0}^{k-1} \binom{k-1}{s} (-T)^s (y_{k-s} \otimes 1 + 1\otimes y_{k-s}) +  \sum_{s = 2}^k \sum_{\substack{s_1 + s_2 = s \\ s_1,s_2\ge 1}}(-T)^{k-s} \binom{k-1}{s-1} \cdot  (y_{s_1} \otimes y_{s_2}) \displaybreak[3] \\
=&\sum_{s = 1}^k \sum_{\substack{s_1 + s_2 = s \\ s_1,s_2\ge 0}}(-T)^{k-s} \binom{k-1}{s-1} \cdot  (y_{s_1} \otimes y_{s_2}) \displaybreak[3] \\
=&\sum_{s = 1}^k (-T)^{k-s} \binom{k-1}{s-1}  \Delta_*(y_s) \displaybreak[3] \\
=&\Delta_*\left( \sum_{s = 1}^k (-T)^{k-s} \binom{k-1}{s-1}  y_s  \right) \displaybreak[3] \\
=&\Delta_*(\sft_*^{\vee}(y_k)).
\end{align*}
This proves that $\sft_*^{\vee}$ is a coalgebra homomorphism with respect to $\Delta_*$.

\end{proof}

\begin{prop}\label{lem:sft-explicit}
For $\Psi \in R[[T; Y \rangle\rangle$, the element $\sft_*^{\vee}(\Psi)$ is explicitly given by
\begin{align}
\sft_*^{\vee}(\Psi) =1 +  \sum_{\substack{w = y_{k_1}\cdots y_{k_r} \in Y^* \\ k_1,\ldots,k_r\ge 1, r\ge 1}} \sum_{l\ge 0} (-T)^l
\sum_{\substack{l_1 + \cdots +l_r = l \\ l_1,\ldots, l_r\ge 0}}
\prod_{j=1}^r \binom{k_j + l_j -1}{l_j}
\langle\Psi \mid y_{k_1+l_1}\cdots y_{k_r+l_r} \rangle\, w.
\label{eq:sft-exc}
\end{align}
\end{prop}

\begin{proof}
Let $r\ge 0$ and $k_1,\ldots,k_r \in \mathbb{Z}_{\ge 1}$. 
It is enough to show that, for $k_1,\ldots,k_r \in \mathbb{Z}_{\ge 0}$,
\begin{align*}
\langle \sft_*^{\vee}(\Psi) \mid y_{k_1}\cdots y_{k_r} \rangle
=
\sum_{l\ge 0}\sum_{\substack{l_1 + \cdots +l_r = l \\ l_1,\ldots, l_r\ge 0}}  (-T)^{l}
\left( \prod_{i=1}^r \binom{k_i + l_i-1}{l_i} \right)
\langle \Psi \mid y_{k_1 + l_1}\cdots y_{k_r + l_r} \rangle .
\end{align*}
Indeed, this identity is exactly the coefficient formula corresponding to \eqref{eq:sft-exc}.

Let $(l_1,\ldots,l_r) \in \mathbb{Z}_{\ge 0}^r$. Then we have
\begin{align*}
&\langle \sft_*^{\vee}(y_{k_1 + l_1}\cdots y_{k_r + l_r}) \mid y_{k_1}\cdots y_{k_r} \rangle \displaybreak[3] \\
=&\langle  \sft_*^{\vee}(y_{k_1 + l_1})\cdots \sft_*^{\vee}(y_{k_r + l_r}) \mid y_{k_1}\cdots y_{k_r} \rangle \displaybreak[3] \\
=& {\footnotesize \left\langle 
\sum_{\substack{0\le j_i \le k_i+l_i-1 \\ 1\le i \le r}}
(-T)^{j_1+\cdots + j_r}
\left( \prod_{i=1}^r \binom{k_i + l_i-1}{j_i} \, y_{k_i + l_i - j_i} \right)
\middle|\, y_{k_1}\cdots y_{k_r} \right\rangle } \displaybreak[3] \\
=& (-T)^{l_1+\cdots + l_r} \left( \prod_{i=1}^r \binom{k_i + l_i-1}{l_i}\right). \displaybreak[3]
\end{align*}

Since $\sft_*^{\vee}$ maps each $y_m$ to an $R[[T]]$-linear combination of $y_{m-j}$ ($0\le j\le m-1$) and it does not change the length of the word. Hence, only words of the form $y_{k_1+l_1}\cdots y_{k_r+l_r}$ with $l_i\ge 0$ can contribute to the coefficient of $y_{k_1}\cdots y_{k_r}$ in $ \sft_*^{\vee}(\Psi)$. Therefore, we have
\begin{align*}
&\langle \sft_*^{\vee}(\Psi) \mid y_{k_1}\cdots y_{k_r} \rangle \displaybreak[3] \\
=&\sum_{l_1,\ldots,l_r\ge 0}  
\langle \Psi \mid y_{k_1 + l_1}\cdots y_{k_r + l_r} \rangle \cdot
\langle \sft_{*}^{\vee}(y_{k_1+l_1} \cdots y_{k_r + l_r}) \mid y_{k_1}\cdots y_{k_r} \rangle \displaybreak[3] \\
=& \sum_{l_1,\ldots,l_r\ge 0}  
(-T)^{l_1+\cdots + l_r} \left( \prod_{i=1}^r \binom{k_i + l_i-1}{l_i}\right)
\langle \Psi \mid y_{k_1 + l_1}\cdots y_{k_r + l_r} \rangle \displaybreak[3] \\
=&\sum_{l\ge 0}\sum_{\substack{l_1 + \cdots +l_r = l \\ l_1,\ldots, l_r\ge 0}}  (-T)^{l}
\left( \prod_{i=1}^r \binom{k_i + l_i-1}{l_i} \right)
\langle \Psi \mid y_{k_1 + l_1}\cdots y_{k_r + l_r} \rangle . \displaybreak[3]
\end{align*}
This completes the proof.
\end{proof}

\begin{prop}
For $\Psi \in R[[T ;  Y \rangle\rangle$ and $k_1,\ldots,k_r \in \mathbb{Z}_{>0}$, we have
\begin{align}
\langle \sft_*^{\vee}(\Psi) \mid  y_{k_1}\cdots y_{k_r} \rangle = \sum_{l\ge 0} T^l \langle \sft_{*,l} (\Psi) \mid y_{l+1}y_{k_1}\cdots y_{k_r} \rangle. \label{eq:sft-vee-l}
\end{align}
\end{prop}

\begin{proof}
For $y_{k_1}\cdots y_{k_r} \in Y^*$, we have
\begin{align*}
&\langle \sft_{*}^{\vee}(\Psi) \mid y_{k_r}\cdots y_{k_1} \rangle\\
\stackrel{\eqref{eq:sft-exc}}{=}&  \sum_{l\ge 0} (-T)^l
\sum_{\substack{l_1 + \cdots +l_r = l \\ l_1,\ldots, l_r\ge 0}}
\prod_{j=1}^r \binom{k_j + l_j -1}{l_j}
\langle\Psi \mid y_{k_1+l_1}\cdots y_{k_r+l_r} \rangle\, \\
\stackrel{\eqref{eq:sft-l}}{=}& \sum_{l\ge 0} T^l \langle \sft_{*,l}(\Psi)\mid y_{l+1}y_{k_1}\cdots y_{k_r} \rangle.
\end{align*}
\end{proof}

\begin{prop}\label{prop:Ad-explicit}
For  $\Phi \in \DMR_0(R)$, we have
\begin{align}
\mathbf{q}_{\#}^{\vee}(\Phi^{-1}x_1\Phi) = &\inv \circ \sft_*^{\vee}(\Phi_{\star})\cdot \Phi_{\star}  \label{eq:Ad-explicit}
\end{align}
\end{prop}
\begin{proof}
By \eqref{eq:sft-vee-l}, we obtain
{\normalsize
\begin{align*}
\mathbf{q}_{\#}^{\vee}(\Phi^{-1}x_1\Phi) =&
\sum_{w  \in Y^*}\left( \sum_{l\ge 0} T^l \langle \mathbf{q}^{\vee}(\Phi^{-1}x_1\Phi) \mid y_{l+1}w \rangle \right) w  \displaybreak[3] \\
=&\sum_{w \in Y^*}\left( \sum_{l_1\ge 0} T^{l_1} \left\langle \left(\sum_{l_2\ge 0}  \sft_{*,l_2} \circ \inv (\Phi_{\star}) \right) \cdot \Phi_{\star}  \middle| y_{{l_1}+1}w \right\rangle \right) w  \displaybreak[3] \\
=& \sum_{w \in Y^*}\left( \sum_{l_1\ge 0} T^{l_1} \sum_{uv = w} \left\langle \left(\sum_{l_2\ge 0}  \sft_{*,l_2} \circ \inv(\Phi_{\star}) \right) \middle| y_{l_1+1}u \right\rangle \cdot \left\langle  \Phi_{\star}  \middle| v \right\rangle \right)w. \\
\intertext{Since every word in $\sft_{*,l_2}\circ\inv(\Phi_\star)$ begins with $y_{l_2+1}$, we must have $l_1=l_2$. Hence,}
=& \sum_{w \in Y^*}\left( \sum_{l\ge 0} T^{l} \sum_{uv = w} \left\langle \left(\sum_{l\ge 0}  \sft_{*,l} \circ \inv(\Phi_{\star}) \right) \middle| y_{l+1}u \right\rangle \cdot \left\langle  \Phi_{\star}  \middle| v \right\rangle \right)w \\
=& \sum_{w \in Y^*}\left( \sum_{uv = w} \left\langle \left(  \sum_{l\ge 0} T^{l}   \sft_{*,l} \circ \inv(\Phi_{\star}) \right) \middle| y_{l+1}u \right\rangle \cdot \left\langle  \Phi_{\star}  \middle| v \right\rangle \right)w  \displaybreak[3] \\
\stackrel{\eqref{eq:sft-vee-l}}{=}&\sum_{w \in Y^*}\left( \sum_{uv = w} \left\langle \sft_{*}^{\vee} \circ \inv(\Phi_{\star}) \middle| u \right\rangle \cdot \left\langle  \Phi_{\star}  \middle| v \right\rangle \right)w \\
=&\sft_{*}^{\vee} \circ \inv (\Phi_{\star}) \cdot \Phi_{\star}.
\end{align*}
}
\end{proof}

\begin{prop}[ = Proposition \ref{prop:t-adic} (\ref{thm:t-adic-2})]\label{prop:Ad-har}
For $\Phi \in \DMR_0(R)$,
we have
\begin{align}
\Delta_*(\sft_{*}^{\vee} \circ \inv(\Phi_{\star}) \cdot \Phi_{\star}) = (\sft_{*}^{\vee} \circ \inv(\Phi_{\star}) \cdot \Phi_{\star}) \otimes (\sft_{*}^{\vee} \circ \inv(\Phi_{\star}) \cdot \Phi_{\star}). \label{eq:ad-harmonic-pre}
\end{align}
Moreover, by \eqref{eq:Ad-explicit}, we have
\begin{align}
\Delta_*\bigl(\mathbf{q}_{\#}^{\vee}(\Phi^{-1} x_1 \Phi)\bigr)
= \mathbf{q}_{\#}^{\vee}(\Phi^{-1} x_1 \Phi)\otimes
\mathbf{q}_{\#}^{\vee}(\Phi^{-1} x_1 \Phi).
\label{eq:ad-harmonic}
\end{align}
Here, we view $\Delta_*$ as a $R[[T]]$-algebra homomorphism via the scalar extension to $R[[T]]$.
\end{prop}

\begin{proof}
Since $\sft_*^{\vee}$ and $\inv$ is a coalgebra homomorphism with respect to $\Delta_*$ and $\Delta_*(\Phi_{\star}) = \Phi_{\star} \otimes \Phi_{\star}$ holds, we have
\begin{align*}
\Delta_*(\sft_*^{\vee}\circ \inv(\Phi_{\star}) \cdot \Phi_{\star}) =& \Delta_*( \sft_*^{\vee} \circ \inv(\Phi_{\star})) \cdot \Delta_*(\Phi_{\star}) \\
=&(\sft_*^{\vee,\otimes 2} \circ \inv^{\otimes 2} ) \Delta_*(\Phi_{\star}) \cdot \Delta_*(\Phi_{\star}) \\
=&( \sft_*^{\vee,\otimes 2} \circ  \inv^{\otimes 2} ) (\Phi_{\star} \otimes \Phi_{\star}) \cdot (\Phi_{\star} \otimes \Phi_{\star})\\
=& (\sft_*^{\vee}\circ \inv(\Phi_{\star})) \otimes (\sft_*^{\vee}\circ \inv(\Phi_{\star})) \cdot  (\Phi_{\star} \otimes \Phi_{\star}) \\
=& (\sft_*^{\vee}\circ \inv(\Phi_{\star}) \cdot \Phi_{\star}) \otimes (\sft_*^{\vee}\circ \inv(\Phi_{\star}) \cdot \Phi_{\star}).
\end{align*}
\end{proof}

Since the condition $\Delta_{*}(\Phi_{\star})=\Phi_{\star}\otimes \Phi_{\star}$ corresponds to the harmonic product relations for $*$-regularized MZVs, \eqref{eq:ad-harmonic-pre} yields $\Q$-algebraic relations among $*$-AdMZVs.
We call these $\Q$-algebraic relations among $*$-AdMZVs the \emph{adjoint harmonic product relations}.

Thus, we complete the proof of Proposition \ref{prop:t-adic}.

\begin{rmk}
Let $t$ be an indeterminate.
Ono, Seki, and Yamamoto \cite{Ono-Seki-Yamamoto} introduced the $t$-adic SMZVs.
More precisely, the $t$-adic SMZV is defined by
\[
\zeta_{\widehat{\S}}(k_1,\ldots,k_r)
:= \sum_{l\ge 0} \Admzv{l}{k_1,\ldots,k_r} t^l
\in (\mathcal{Z}/\zeta(2)\mathcal{Z})[[t]]
\]
for $(k_1,\ldots,k_r) \in \mathbb{Z}_{>0}^r$.
Independently, Jarossay \cite[Definition~2.3.1(ii)]{Jarossay} introduced the $\Lambda$-adjoint multiple zeta values. Here, $\Lambda$ is an indeterminate.
They coincide with the $t$-adic SMZVs after replacing $t$ by $\Lambda$.
Ono, Seki, and Yamamoto also discussed certain analogs of the double shuffle relations for $t$-adic SMZVs.
Comparing coefficients of powers of $t$, we see that any $\Q$-algebraic relations among $t$-adic SMZVs reduce to $\Q$-algebraic relations among AdMZVs.

Based on the discussion so far, we obtain two families of $\Q$-algebraic relations among AdMZVs, those due to Jarossay and those due to Ono, Seki, and Yamamoto.
They differ only in the shuffle relations.
More precisely, for some $\Phi \in \ide{h}^{\vee}$, Proposition~\ref{prop:t-adic}. (\ref{thm:t-adic-1}) gives Jarossay's linear shuffle relations, whereas the linear shuffle relations of Ono, Seki, and Yamamoto are given by
\begin{align*}
\langle \Phi \mid x_0^l x_1(u\shuffle v x_i)\rangle
= -\sum_{l_1+l_2 = l} \langle \Phi \mid x_0^{l_1}x_i(x_0^{l_2}\shuffle S_X(v))x_1u\rangle
\end{align*}
for $u$ and $v\in X^*$, and $i\in\{0,1\}$.
Jarossay proved that these two types of linear shuffle relations are equivalent, see \cite[Proposition 3.4.1]{Jarossay}.
\end{rmk}

Motivated by Proposition~\ref{prop:t-adic}, Jarossay introduced the \emph{adjoint double shuffle set}, defined by the adjoint double shuffle relations.
\begin{df}[{\cite[Definition 3.2.6]{Jarossay}}]\label{def:AdDMR}
Let $R$ be a unital commutative associative $\Q$-algebra.
We define $\AdDMR_0(R)$ to be the set of $\Phi \in R\langle\!\langle x_0,x_1 \rangle\!\rangle$ satisfying the following conditions:
\begin{enumerate}
\item\label{def:AdDMR-2} $\Phi - x_1 \in R\cotimes \ide{F}_2^{\ge 4}$,
\item\label{def:AdDMR-3} $\Delta_{*}(\Phi_{\#}) = \Phi_{\#} \otimes \Phi_{\#}$, where $\Phi_{\#} := \mathbf{q}_{\#}^{\vee}(\Phi)$.
\end{enumerate}
\end{df}

Due to Proposition \ref{prop:t-adic}, for each unital commutative associative $\Q$-algebra $R$, we can define
\begin{align}
\Ad^R(x_1)\colon \DMR_0(R) \longrightarrow \AdDMR_0(R), \qquad \phi \longmapsto \phi^{-1} x_1 \phi. \label{eq:Ad2}
\end{align}
Jarossay asked a question about whether this map is bijective.
\begin{qu}[{\cite[Question 3.2.8]{Jarossay}}]\label{qu:Jarossay}
Is the map \eqref{eq:Ad2} a bijection for any unital commutative associative $\Q$-algebra $R$?
\end{qu}

In the next section, we show that the nonsurjectivity of $\Ad^R(x_1)$ \eqref{eq:Ad2} for some $\Q$-algebra.

\section{Nonsurjectivity and Refined Questions}

In this section, we show that the map $\Ad^R(x_1)$ appearing in Question~\ref{qu:Jarossay} is not surjective.
This negative result naturally leads to refined versions of Question~\ref{qu:Jarossay}.
We also observe that, after imposing natural additional conditions, called adjoint conditions, one obtains a more precise formulation of the question.

Question~\ref{qu:Jarossay}

\subsection{Nonsurjectivity of $\Ad^R(x_1)$ for some $\Q$-algebra $R$}

In this subsection, we prove that the answer to Question~\ref{qu:Jarossay} is negative.
To study Question~\ref{qu:Jarossay}, we focus on the tangent spaces of $\DMR_0(\Q)$ and $\AdDMR_0(\Q)$.
Set
\[
\dmr_0 := \{\Psi \in \ide{h}^{\vee} \mid 1 + \varepsilon \Psi \in \DMR_0(\Q[\varepsilon]/(\varepsilon^2))\}.
\]
By Racinet \cite[Section 3.3.1, 3.3.8]{Racinet}, $\dmr_0$ coincides with
\[
\{\psi \in \ide{F}_2^{\ge 3} \mid \Delta_*(\psi_{\star}) = \psi_{\star} \otimes 1 + 1\otimes \psi_{\star} \},
\]
where $\psi_{\star} := \mathbf{q}(\psi) + \sum_{n\ge 2} \frac{(-1)^n}{n} \langle \psi \mid x_0^{n-1} x_1\rangle y_1^n$.

\begin{rmk}
For $\psi_1,\psi_2 \in \dmr_0$, define
\[
\{\psi_1,\psi_2\}
:= d_{[x_1,\psi_1]}(\psi_2)-d_{[x_1,\psi_2]}(\psi_1)+[\psi_1,\psi_2],
\]
where, for $f\in \ide{h}^{\vee}$, $d_f$ denotes the derivation of $\ide{h}^{\vee}$ determined by
\[
d_f(x_0) :=0,\qquad d_f(x_1) :=f.
\]
Then $\dmr_0$ is a Lie algebra with bracket $\{\mathchar`-,\mathchar`-\}$ which is shown in \cite[Proposition 4.A.1]{Racinet}.
\end{rmk}

On the adjoint side, we introduce:
\begin{df}\label{def:addmr}
We define
\[
\addmr_0 := \{\Psi \in  \ide{h}^{\vee} \mid x_1 + \varepsilon \Psi \in \AdDMR_0(\Q[\varepsilon]/(\varepsilon^2))\}.
\]
\end{df}

\begin{prop}\label{prop:addmr}
The tangent space $\addmr_0$ is the $\Q$-vector space of all $\Psi\in \ide{h}^{\vee}$ satisfying the following conditions:
\begin{enumerate}
\item\label{def:addmr--2} $\Psi \in \ide{F}_2^{\ge 4}$,
\item\label{def:addmr--3} $\Delta_{*}(\Psi_{\#}) = \Psi_{\#} \otimes 1 + 1 \otimes \Psi_{\#}$.
\end{enumerate}
\end{prop}

\begin{proof}
Assume that $x_1 + \varepsilon \Psi \in \AdDMR_0(\Q[\varepsilon]/(\varepsilon^{2}))$.
First, we show that $\Psi$ satisfies Conditions~(\ref{def:addmr--2}) and (\ref{def:addmr--3}).
By definition of $\AdDMR_0(\Q[\varepsilon]/(\varepsilon^{2}))$, we have
\begin{align}
(x_1 + \varepsilon \Psi) - x_1  &\in \Q[\varepsilon]/(\varepsilon^2) \cotimes \ide{F}_2^{\ge 4}, \label{eq:addmr-1}\\
\Delta_{*} ((x_1 + \varepsilon \Psi)_{\#})
&= (x_1 + \varepsilon \Psi)_{\#} \otimes (x_1 + \varepsilon \Psi)_{\#}. \label{eq:addmr-2}
\end{align}
By \eqref{eq:addmr-1}, we have $\varepsilon\Psi\in \Q[\varepsilon]/(\varepsilon^2)\cotimes \ide{F}_2^{\ge 4}$.
Taking the coefficient of $\varepsilon$, we obtain $\Psi\in \ide{F}_2^{\ge 4}$.

Next, we prove (\ref{def:addmr--3}).
Using $(x_1)_{\#}=1$ and $(x_1+\varepsilon\Psi)_{\#}=1+\varepsilon\Psi_{\#}$, we compute
\begin{align*}
\Delta_{*}((x_1 + \varepsilon \Psi)_{\#})
&\hspace{-0.15cm}\stackrel{\eqref{eq:addmr-2}}{=} (x_1 + \varepsilon \Psi)_{\#}\otimes (x_1 + \varepsilon \Psi)_{\#} \\
&=(1 + \varepsilon \Psi_{\#})\otimes (1 + \varepsilon \Psi_{\#})\\
&= 1\otimes 1 + \varepsilon(\Psi_{\#} \otimes 1 + 1\otimes \Psi_{\#}).
\end{align*}
Comparing the coefficients of $\varepsilon$ yields $\Delta_{*}(\Psi_{\#}) = \Psi_{\#} \otimes 1 + 1 \otimes \Psi_{\#}$.

Conversely, suppose that $\Psi \in\ide{h}^{\vee}$ satisfies
\begin{align}
\Psi &\in \ide{F}_2^{\ge4}, \notag\\
\Delta_*(\Psi_\#) &= \Psi_\#\otimes1+1\otimes\Psi_\#. \label{eq:addmr-3}
\end{align}
Then $(x_1+\varepsilon\Psi)-x_1=\varepsilon\Psi\in \Q[\varepsilon]/(\varepsilon^2)\cotimes \ide{F}_2^{\ge4}$.
Moreover, since $(x_1+\varepsilon\Psi)_\#=1+\varepsilon\Psi_\#$, we have
\begin{align*}
\Delta_{*}((x_1 + \varepsilon \Psi)_{\#})
&= \Delta_*(1 + \varepsilon \Psi_{\#}) \\
&= 1\otimes 1 + \varepsilon \Delta_* (\Psi_{\#})\\
&\hspace{-0.15cm}\stackrel{\eqref{eq:addmr-3}}{=} 1\otimes 1  + \varepsilon (\Psi_{\#} \otimes 1 + 1\otimes \Psi_{\#})\\
&= (1 + \varepsilon \Psi_{\#}) \otimes (1 + \varepsilon \Psi_{\#})\\
&= (x_1+\varepsilon\Psi)_\# \otimes (x_1+\varepsilon\Psi)_\#.
\end{align*}
Thus, $x_1+\varepsilon\Psi\in \AdDMR_0(\Q[\varepsilon]/(\varepsilon^2))$, which means that $\Psi\in\addmr_0$.
\end{proof}

\begin{prop}\label{prop:dAd}
The differential of the map
\[
\Ad^{\Q}(x_1):\DMR_0(\Q)\longrightarrow \AdDMR_0(\Q),\qquad \phi\longmapsto \phi^{-1} x_1 \phi
\]
at $1$ is the $\Q$-linear map
\[
d\Ad^{\Q}(x_1)_1:\dmr_0\longrightarrow \addmr_0,\qquad \Psi\longmapsto [x_1,\Psi].
\]
\end{prop}

\begin{proof}
By the definition of the tangent spaces, an element $\Psi\in\dmr_0$ corresponds to
\[
1+\varepsilon\Psi\in \DMR_0(\Q[\varepsilon]/(\varepsilon^2)).
\]
Since $(1+\varepsilon\Psi)^{-1}=1-\varepsilon\Psi$, we compute
\begin{align*}
\Ad^{\Q}(x_1)(1+\varepsilon\Psi)
&=(1+\varepsilon\Psi)^{-1}x_1(1+\varepsilon\Psi)\\
&=(1-\varepsilon\Psi)x_1(1+\varepsilon\Psi)\\
&=x_1+\varepsilon(x_1\Psi-\Psi x_1)\\
&=x_1+\varepsilon[x_1,\Psi].
\end{align*}
Thus the coefficient of $\varepsilon$ is $[x_1,\Psi]$. This shows that
$d\Ad^{\Q}(x_1)_1(\Psi)=[x_1,\Psi]$.
\end{proof}
Define
\[
\ad(x_1) := d \Ad^{\Q}(x_1)_1.
\]
At the level of tangent space, Question~\ref{qu:Jarossay} becomes the following statement.
\begin{qu}\label{qu:Jarossay-Lie}
Is the $\Q$-linear map
\[
\ad(x_1)\colon \dmr_0 \longrightarrow \addmr_0,\qquad \Psi \longmapsto [x_1,\Psi]
\]
the isomorphism?
\end{qu}

\begin{rmk}
Note that $\ad(x_1)$ is a homogeneous map which shifts weight by $1$.
\end{rmk}

To study Question~\ref{qu:Jarossay-Lie}, we decompose $\dmr_0$ and $\addmr_0$ into their weight-homogeneous components:
\[
\dmr_0 = \prod_{k\ge0}\dmr_0^{(k)},\qquad
\addmr_0 = \prod_{k\ge0}\addmr_0^{(k)}.
\]
Here,
\begin{align*}
\dmr_0^{(k)} := \dmr_0 \cap \ide{F}_2^{k},\qquad
\addmr_0^{(k)} := \addmr_0 \cap \ide{F}_2^{k}.
\end{align*}

\begin{prop}\label{prop:ex-5}
Set
\[
\omega
:=
-\bigl[[x_0,x_1],[x_0,[x_0,x_1]]\bigr]
+\bigl[[x_0,x_1],[x_1,[x_0,x_1]]\bigr]
-\bigl[x_1,[x_0,[x_0,[x_0,x_1]]]\bigr]
+\bigl[x_1,[x_1,[x_0,[x_0,x_1]]]\bigr].
\]
Then $\omega \in \addmr_0^{(5)}$.
\end{prop}

\begin{proof}
It is immediate that $\omega$ satisfies Condition~(\ref{def:addmr--2}) in Proposition~\ref{prop:addmr}.
Therefore, it suffices to verify Condition~(\ref{def:addmr--3}).
By direct calculation, we obtain
\begin{align*}
\omega_{\#} := \mathbf{q}_{\#}^{\vee}(\omega)
&= T^3y_1 + T^2y_1^2 - 2T^2y_2 + 3T\,y_1y_2 - 3T\,y_2y_1\\
&= T^3 y_1 + T^2 (y_1^2 - 2y_2) + 3T [y_1,y_2].
\end{align*}
Since the application $\Delta_*$ to $y_1$, $y_1^2 - 2y_2$ and $[y_1,y_2]$ gives
\begin{align*}
\Delta_*(y_1) &= y_1 \otimes 1 + 1\otimes y_1,\\
\Delta_*(y_1^2 - 2y_2) &= (y_1^2 - 2y_2) \otimes 1 + 1\otimes (y_1^2 - 2y_2),\\
\Delta_*([y_1,y_2]) &= [y_1,y_2] \otimes 1 + 1\otimes [y_1,y_2],
\end{align*}
we obtain
\begin{align*}
\Delta_*(\omega_{\#})
&= \Delta_*\bigl(T^3 y_1 + T^2 (y_1^2 - 2y_2) + 3T [y_1,y_2]\bigr)\\
&= T^3\Delta_*(y_1) + T^2\Delta_*(y_1^2 - 2y_2) + 3T\Delta_*\bigl([y_1,y_2]\bigr)\\
&= T^3(y_1 \otimes 1 + 1\otimes y_1)
 + T^2\bigl((y_1^2 - 2y_2)\otimes 1 + 1\otimes (y_1^2 - 2y_2)\bigr)\\
&\quad + 3T\bigl([y_1,y_2]\otimes 1 + 1\otimes [y_1,y_2]\bigr)\\
&= \omega_{\#} \otimes 1 + 1\otimes \omega_{\#}.
\end{align*}
\end{proof}

\begin{prop}
The $\Q$-linear map
\[
\ad(x_1)\colon \dmr_0\longrightarrow \addmr_0,\qquad
\phi \longmapsto [x_1,\phi]
\]
is not surjective.
\end{prop}

\begin{proof}
Since $\ad(x_1)$ increases the weight by $1$, Question~\ref{qu:Jarossay-Lie} would suggest that
$\dim_{\Q}\dmr_0^{(k)}=\dim_{\Q}\addmr_0^{(k+1)}$ for all $k\ge 0$.
However, the table in \cite[pp.5]{Schneps-ds} shows that $\dim_{\Q}\dmr_0^{(4)}=0$, whereas Proposition~\ref{prop:ex-5} gives $\dim_{\Q}\addmr_0^{(5)}\ge 1$.
This contradicts the dimension equality above, so $\ad(x_1)$ is not surjective.
\end{proof}

\begin{thm}[{Main Theorem~\ref{mt:Jarossay-negative}}]\label{thm:Jarossay-negative}
If we take $R=\Q[\varepsilon]/(\varepsilon^2)$, then the map
\begin{align}
\Ad^R(x_1)\colon \DMR_0(R)\longrightarrow \AdDMR_0(R),\qquad
\phi \longmapsto \phi^{-1}x_1\phi \label{eq:Ad3}
\end{align}
is not surjective.
\end{thm}

\begin{proof}
Assume that the map \eqref{eq:Ad3}
is bijective for $R=\Q[\varepsilon]/(\varepsilon^2)$.
By Proposition~\ref{prop:ex-5}, we have $\omega\in \addmr_0$.
Hence $x_1+\varepsilon \omega\in \AdDMR_0(R)$ by Definition~\ref{def:addmr}.
By the assumed bijectivity, there exists $\phi\in \DMR_0(R)$ such that
\[
\phi^{-1} x_1 \phi=x_1+\varepsilon \omega.
\]
Since $\phi\equiv 1\pmod{\varepsilon}$, we can write $\phi=1+\varepsilon\psi$ with $\psi\in \dmr_0$.
Using $(1+\varepsilon\psi)^{-1}=1-\varepsilon\psi$, we compute
\[
\phi^{-1} x_1 \phi=(1-\varepsilon\psi)x_1(1+\varepsilon\psi)
= x_1+\varepsilon[x_1,\psi].
\]
Comparing the coefficients of $\varepsilon$, we obtain $\omega=[x_1,\psi]$.

Now $\omega$ has weight $5$ and $\ad(x_1)$ increases the weight by $1$, so $\psi$ must have weight $4$, that is, $\psi\in \dmr_0^{(4)}$.
However, the table in \cite[pp.5]{Schneps-ds} shows that $\dim_{\Q}\dmr_0^{(4)}=0$, hence $\psi=0$.
This implies $\omega=0$, which contradicts Proposition~\ref{prop:ex-5}.
Therefore, $\Ad^{R}(x_1)$ cannot be bijective for all $\Q$-algebras $R$, and the answer to Question~\ref{qu:Jarossay} is negative.
\end{proof}

In the next section, we shall refine  Question \ref{qu:Jarossay} and \ref{qu:Jarossay-Lie} by imposing the additional conditions, respectively.

\subsection{Refinement of questions on $\AdDMR_0(R)$}

In this subsection, we formulate refined versions of the questions concerning $\AdDMR_0(R)$.
We now formulate a refined version of Questions~\ref{qu:Jarossay} and \ref{qu:Jarossay-Lie}.

First, we refine Question~\ref{qu:Jarossay-Lie} by introducing the following $\Q$-vector space:
\[
\Fad := \{\Psi \in \ide{F}_2^{\ge 3} \mid {}^\exists\, \psi \in \ide{F}_2^{\ge 2}\ \text{such that}\ \Psi = [x_1,\psi]\}.
\]

By comparing the dimensions of the weight-homogeneous components of $\dmr_0$ and $\addmr_0\cap \Fad$, we observe that
\[
\dim_{\Q}\dmr_0^{(k)} = \dim_{\Q}\bigl(\addmr_0^{(k+1)} \cap \Fad\bigr)
\]
for $k = 0,\ldots,10$ (see Table~\ref{tab:table-dimension}\footnote{The entries of Table \eqref{tab:table-dimension} were computed by a program written by the author, which is based on Hirose's library.}).

\begin{table}[h]
\centering
\begin{tabular}{c|ccccccccccccc}
$k$ & $0$ & $1$ & $2$ & $3$ & $4$ & $5$ & $6$ & $7$ & $8$ & $9$ & $10$ &  $\cdots$ \\
\hline
$\dim \dmr_0^{(k)}$ & $0$ & $0$ & $0$ & $1$ & $0$ & $1$ & $0$ & $1$ & $1$ & $1$ & $1$ &   $\cdots$ \\
$\dim \addmr_0^{(k+1)} \cap \Fad$
& $0$ & $0$ & $0$ & $1$ & $0$ & $1$ & $0$ & $1$ & $1$ & $1$ & $1$ & $\cdots$
 \\
\end{tabular}
\caption{Dimensions of $\dmr_0^{(k)}$ and $\addmr_0^{(k+1)} \cap \Fad$}
\label{tab:table-dimension}
\end{table}

By the above table, we establish the following, which is the refinement of Question~\ref{qu:Jarossay-Lie}.

\begin{qu}\label{qu:addmr-Lie}
Is the $\Q$-linear map
\[
\ad(x_1)\colon \dmr_0 \longrightarrow \addmr_0 \cap \Fad,\qquad \Psi \longmapsto [x_1,\Psi]
\]
an isomorphism?
\end{qu}

\begin{rmk}\label{rmk:addmr-bracker}
Assume that Question~\ref{qu:addmr-Lie} holds. Then, $\addmr_0 \cap \Fad$ is expected to be a Lie algebra with the Lie bracket given by
\[
\{\Psi_1,\Psi_2\}_1 := d_{\Psi_1}(\Psi_2) - d_{\Psi_2}(\Psi_1), \qquad (\Psi_1, \Psi_2 \in \addmr_0\cap \Fad).
\]
Note that $(\addmr_0 \cap \Fad, \{\mathchar`-,\mathchar`-\})$ is also expected to be a graded Lie algebra after shifting the weight by $-1$:
\[
\addmr_0 \cap \Fad = \bigoplus_{k\ge 0} (\addmr_0 \cap \Fad)^{[k]},\qquad (\addmr_0 \cap \Fad)^{[k]} := \addmr_0^{(k+1)} \cap \Fad.
\] 
\end{rmk}

Next, we refine Question~\ref{qu:Jarossay}.
For each unital commutative associative $\Q$-algebra $R$, we define
\[
\FAd(R)
:=
\left\{
\phi^{-1}x_1\phi
\ \middle|\
\phi\in 1+R\langle\langle x_0,x_1\rangle\rangle_{\ge 2},
\quad
\Delta_{\shuffle}(\phi)=\phi\otimes\phi
\right\}.
\]
and consider $\AdDMR_0(R) \cap \FAd(R)$. Then we have the following embedding:
\[
\Ad^R(x_1)\colon \DMR_0(R) \longrightarrow \AdDMR_0(R) \cap \FAd(R),\qquad \phi \longmapsto \phi^{-1} x_1 \phi.
\]

We note that the tangent space of $\AdDMR_0(\Q)\cap \FAd(\Q)$ at $x_1$ is
\[
T_{x_1}\bigl(\AdDMR_0(\Q)\cap \FAd(\Q)\bigr)=\addmr_0 \cap \Fad,
\]
by repeating the arguments of Propositions~\ref{prop:addmr} and \ref{prop:dAd}.
Thus, we establish the following Question which is the refinement of Question~\ref{qu:Jarossay}

\begin{qu}\label{qu:AdDMR-grp}
Is the map
\[
\Ad^R(x_1)\colon \DMR_0(R) \longrightarrow \AdDMR_0(R)\cap \FAd(R),\qquad \phi \longmapsto \phi^{-1} x_1 \phi
\]
a bijection for each unital commutative associative $\Q$-algebra $R$?
\end{qu}
\begin{rmk}
Assume that Question~\ref{qu:AdDMR-grp} holds. Then, for any $\Q$-algebra $R$, the set
$\AdDMR_0(R)\cap\FAd(R)$ is expected to be a group with the product given by
\[
\Phi_1 \circledast_1 \Phi_2 := \kappa_{\Phi_1}(\Phi_2)
\qquad(\Phi_1,\Phi_2\in \AdDMR_0(R)\cap\FAd(R)).
\]
\end{rmk}

\section{An additional result}

In this section, we further discuss another main theorem related to Question \ref{qu:addmr-Lie}.

\subsection{A Lie algebra $\tm_1$}

First, we introduce a new Lie algebra $\tm_1$ which includes $\addmr_0$ as a linear space.

We define
\[
\tm_1= \{\Psi \in  \ide{h}^{\vee} \mid \langle \Psi \mid x_1 \rangle = 0,  \langle \Psi \mid x_0^l \rangle = 0 \text{ ($l\in\mathbb{Z}_{\ge 0}$)}\}.
\]
For $\Psi_1$, $\Psi_2 \in \tm_1$, we define a $\Q$-bilinear binary operation $\{\mathchar`-,\mathchar`-\}_1$ \cite[Appendix A.1.2]{Jarossay2020} on $\tm_1$ by
\[
\{\Psi_1,\Psi_2\}_{1} := d_{\Psi_1}(\Psi_2) - d_{\Psi_2}(\Psi_1).
\]
Note that the binary operator $\{\mathchar`-,\mathchar`-\}_1$ coincides with the operator appearing in Remark \ref{rmk:addmr-bracker} after restriction to $\addmr \cap \Fad$.

\begin{prop}\label{prop:tm1}
A pair $(\tm_1,\{\mathchar`-,\mathchar`-\}_1)$ forms a Lie algebra.
\end{prop}

\begin{proof}
Let $\der(\ide{h}^{\vee})$ be a $\Q$-linear space consisting of all derivations on $\ide{h}^{\vee}$.
We define
\[
\widetilde{\tau}:\ide{h}^{\vee}\to \der(\ide{h}^{\vee}),\qquad
\Psi\mapsto d_\Psi.
\]
Since a derivation on \(\ide{h}^{\vee}\) is determined by its values on \(x_0,x_1\), $\widetilde{\tau}$ is an injection.

We first show that for any \(\Psi\in \tm_1\), one has
\[
d_\Psi(\tm_1)\subset \tm_1.
\]
Indeed, every word appearing in an element of \(\tm_1\) is neither of the form \(x_0^l\) (\(l\ge 0\)) nor equal to \(x_1\). Since \(d_\Psi(x_0)=0\) and \(d_\Psi(x_1)=\Psi\), applying \(d_\Psi\) to such a word replaces one occurrence of \(x_1\) by a word appearing in \(\Psi\), which is again neither of the form \(x_0^l\) nor equal to \(x_1\). Hence, the resulting word still belongs to \(\tm_1\), proving the claim.

Therefore, for \(\Psi_1,\Psi_2\in \tm_1\),
\[
\{\Psi_1,\Psi_2\}_1
=
d_{\Psi_1}(\Psi_2)-d_{\Psi_2}(\Psi_1)\in \tm_1,
\]
so the bracket is well-defined on \(\tm_1\).

Moreover,
\[
[d_{\Psi_1},d_{\Psi_2}](x_0)=0,\qquad
[d_{\Psi_1},d_{\Psi_2}](x_1)
=
d_{\Psi_1}(\Psi_2)-d_{\Psi_2}(\Psi_1)
=
\{\Psi_1,\Psi_2\}_1.
\]
Hence
\[
[d_{\Psi_1},d_{\Psi_2}]
=
d_{\{\Psi_1,\Psi_2\}_1}.
\]
Since $\der(\ide{h}^{\vee})$ is a Lie algebra with respect to $[\mathchar`-,\mathchar`-]$ and $\Im \widetilde{\tau}$ is closed under the commutator, $\mathrm{Im} \widetilde{\tau}$ is the Lie subalgebra of $\der(\ide{h}^{\vee})$.
Thus $\tm_1$ is a Lie algebra with respect to $\{\mathchar`-,\mathchar`-\}_1$ via the injective map $\widetilde{\tau}$.
\end{proof}

\subsection{Parity results}

Let us define $\mathcal{Z}_{>0} := \Span_{\mathbb{Q}} \{\zeta(\mathbf{k}) \mid \mathbf{k} \neq \emptyset\}$.
In this subsection, we put  $R = \mathcal{Z} / (\mathcal{Z}_{>0}^2 +\zeta(2)\mathcal{Z})$.
In our main theorem, we utilize the parity results, which state that $\zeta(k_1,\ldots,k_r)$ with $k_1+\dots + k_r + r \equiv 1 \mod 2$ lies in $\Q[\pi^2]$-span of a MZVs of depth less than $r$ (The depth of MZV means the number of entries in an argument of a MZV).
These parity results are first proved analytically in \cite{Tsumura} and later algebraically in \cite{Ihara-Kaneko-Zagier}.
In \cite{Hirose-parity}, Hirose provided an explicit formula of the parity results from the point of view of the multitangent functions.
The significant point of his proof is the following formula, which follows from the functional equation of multitangent functions
\begin{thm}[\cite{Hirose-parity}]
For $(k_1,\ldots,k_r)\in\mathbb{Z}_{>0}^r$, we have
\begin{align}
\begin{aligned}
&\sum_{j = 0}^r (-1)^{k_{j+1}+\cdots + k_r} \zeta^*(k_{j},\ldots,k_1)\zeta^*(k_{j+1},\ldots,k_r) \\
=&\delta^{k_1,\ldots,k_r} + \sum_{j = 1}^r \sum_{\substack{a+2m+b = k_j \\ a,b,m\ge 0}} (-1)^{b+k_{j+1}+\cdots +k_r + m + 1}\\
&\times \frac{(2\pi )^{2m}}{(2m)!}B_{2m}\zeta_a^*(k_{j-1},\ldots,k_1)\zeta_b^*(k_{j+1},\ldots,k_r).
\end{aligned}
\label{eq:parity}
\end{align}
\end{thm}
We omit the definition of $\delta^{k_1,\ldots,k_r}$ but only note that $\delta^{k_1,\ldots,k_r} \in \Q[\pi^2]$.
The left-hand side of \eqref{eq:parity} is no longer AdMZVs, but we must note that AdMZVs appear after taking modulo some sums of a product of MZVs and $\pi^2$.
From this perspective, we determine the $\Q$-linear relations of AdMZVs.

By definition of AdMZVs, we have
\begin{align*}
\Admzv{0}{k_1,\ldots,k_r} &\equiv \zeta(k_1,\ldots,k_r) + (-1)^{k_1+\cdots + k_r}\zeta(k_r,\ldots,k_1) \\
\Admzv{l}{k_1,\ldots,k_r} &\equiv (-1)^{k_1+\cdots + k_r + l} \zeta_l(k_r,\ldots,k_1)
\end{align*}
modulo $\mathcal{Z}_{>0}^2 +\zeta(2)\mathcal{Z}$ for arbitrary $k_1,\ldots,k_r\in\mathbb{Z}_{>0}$ and $l\in\mathbb{Z}_{ \ge 1 }$.

Considering \eqref{eq:parity} as the image of $R$, we obtain
\begin{align*}
&(-1)^{k_1+\cdots + k_r}(\zeta^*(k_1,\ldots,k_r) + (-1)^{k_1+\cdots + k_r} \zeta^*(k_r,\ldots,k_1))\\
\equiv &- \zeta_{k_r}^*(k_{r-1},\ldots,k_1) - (-1)^{k_1+\cdots+k_r} \zeta_{k_1}(k_2,\ldots,k_r)
\end{align*}
as an equality of $R$.
Hence, the parity result of AdMZVs can be written as
\begin{align}
\zeta_{\mathrm{Ad}}(k_1,\ldots,k_r;0)
\equiv&- \zeta_{\mathrm{Ad}}(k_1,\ldots,k_{r-1};k_r) - (-1)^{k_1+\cdots + k_r} \zeta_{\mathrm{Ad}}(k_r,\ldots,k_{2};k_1)  \label{eq:parity-s}
\end{align}
modulo $\mathcal{Z}_{>0}^2 +\zeta(2)\mathcal{Z}$.

Next, we determine the relation for the generating function derived from \eqref{eq:parity-s}.
Put
\[
\Psi_{\shuffle} := \sum_{w \in X^*} (\zeta^{\shuffle}(w) \mod \mathcal{Z}_{>0}^2 + \zeta(2)\mathcal{Z}) \overset{\leftarrow}{w} \in R\langle\langle x_0,x_1 \rangle \rangle.
\]
From \eqref{eq:parity-s}, we get the following equality in $R\langle\langle x_0,x_1 \rangle \rangle$:

\begin{prop}
Let $\Psi_{\ad} := [x_1,\Psi_{\shuffle}]$. Then we have
\begin{align*}
x_1\Psi_{\ad}^{11} x_1 = -x_1\Psi_{\ad}^{01}x_1 + S_X^{\vee}(x_1\Psi_{\ad}^{01}x_1). 
\end{align*}
Equivalently, 
\[
\Psi_{\ad}^{11} + \Psi_{\ad}^{01} - S_X^{\vee}(\Psi_{\ad}^{01}) = 0
\]
holds.
\end{prop}

\begin{proof}
Since
\[
\langle \Psi_{\ad} \mid x_0^lx_1x_0^{k_r}x_1\cdots x_0^{k_1}x_1 \rangle
= \zeta_{\mathrm{Ad}}(k_1+1,\ldots,k_r+1;l)
\]
holds, \eqref{eq:parity-s} turns out to be
\begin{align*}
& \langle \Psi_{\ad} \mid x_1x_0^{k_r}x_1\cdots x_0^{k_1}x_1 \rangle\\
=&-  \langle \Psi_{\ad} \mid x_0^{k_r+1}x_1x_0^{k_{r-1}}\cdots x_0^{k_1}x_1 \rangle - (-1)^{k_1+\cdots +k_r + r} \langle \Psi_{\ad} \mid x_0^{k_1+1}x_1x_0^{k_{2}}\cdots x_0^{k_r}x_1 \rangle.
\end{align*}

Therefore, we have
\begin{align*}
x_1\Psi_{\ad} x_1 =& \sum_{k_1,\ldots,k_r\ge 0}  \langle \Psi_{\ad} \mid x_1x_0^{k_r}x_1\cdots x_0^{k_1}x_1 \rangle x_1x_0^{k_r}x_1\cdots x_0^{k_1}x_1 \\
&
\begin{aligned}
\hspace{-0.3cm}=&\sum_{k_1,\ldots,k_r\ge 0} (-  \langle \Psi_{\ad} \mid x_0^{k_r+1}x_1x_0^{k_{r-1}}\cdots x_0^{k_1}x_1 \rangle\\
&\qquad - (-1)^{k_1+\cdots +k_r + r} \langle \Psi_{\ad} \mid x_0^{k_1+1}x_1x_0^{k_{2}}\cdots x_0^{k_r}x_1 \rangle ) x_1x_0^{k_r}x_1\cdots x_0^{k_1}x_1
\end{aligned}\\
&
\begin{aligned}
\hspace{-0.3cm}=&-\sum_{k_1,\ldots,k_r\ge 0}   \langle x_0\Psi_{\ad}^{01} x_1 \mid x_0^{k_r+1}x_1x_0^{k_{r-1}}\cdots x_0^{k_1}x_1 \rangle _1x_0^{k_r}x_1\cdots x_0^{k_1}x_1 \\
&\qquad + \sum_{k_1,\ldots,k_r\ge 0} \langle \Psi_{\ad} \mid x_0^{k_1+1}x_1x_0^{k_{2}}\cdots x_0^{k_r}x_1 \rangle S_X^{\vee} (x_1x_0^{k_1}x_1\cdots x_0^{k_r}x_1)\\
\end{aligned}\\
=& - x_1\Psi_{\ad}^{01}x_1 + S_X^{\vee}(x_1\Psi_{\ad}^{01}x_1).
\end{align*}

\end{proof}

From this perspective, we say that $\Psi \in \ide{h}^{\vee}$, with $\wt \Psi\ge 2$, satisfies  the \textbf{strong parity result} if $\Psi$ satisfies
\begin{align}
\Psi^{11} + \Psi^{01} - S_{X}^{\vee}(\Psi^{01}) = 0. \label{eq:str-parity-pri}
\end{align}

Assume that $\Psi \in \ide{F}_2$ satisfies Equality \eqref{eq:str-parity-pri}.
Then by Equality \eqref{eq:antipode-X}, we have $S_{X}^{\vee}(\Psi^{10}) = -\Psi^{01}$.
Therefore, the strong parity result turns out to be
\begin{align}
\Psi^{11} + \Psi^{10} + \Psi^{01} = 0. \label{eq:str-prty-3cycle}
\end{align}
 
To conclude, we define two $\Q$-linear subspaces of $\ide{h}^{\vee}$ as follows:
\begin{align*}
V_{\strprty}(R) =& \{\Psi \in  \ide{h}^{\vee} \mid \Psi^{11} + \Psi^{10} + \Psi^{01} = 0\}.
\end{align*}

\begin{prop}\label{prop:deriv-strprty}
We have a Lie subalgebra $V_{\strprty} \cap \Fad \subset \tm_1$.
\end{prop}

This proposition follows from the following lemma:

\begin{lem}
For $\Psi_1$, $\Psi_2 \in  \Fad \cap V_{\strprty}$, we have
\begin{align*}
&d_{\Psi_1}(\Psi_1)^{11} + d_{\Psi_1}(\Psi_2)^{01} + d_{\Psi_1}(\Psi_2)^{10}\\
&\begin{aligned}
=&(\Psi_{1}^{11}x_1\Psi_2^{11} + \Psi_{2}^{11}x_1\Psi_1^{11}) - ( \Psi_1^{01}x_1\Psi_2^{11} + \Psi_2^{01}x_1\Psi_1^{11}) \\
&\quad-(\Psi_1^{11}x_1\Psi_2^{10} + \Psi_2^{11}x_1\Psi_1^{10}) + (\Psi_{1}^{10}x_0\Psi_2^{01}  +  \Psi_{2}^{10}x_0\Psi_1^{01} ).
\end{aligned}
\end{align*}
\end{lem}

\begin{proof}

Since 
\begin{align*}
\begin{aligned}
&x_1d_{\Psi_1}(x_1\Psi_2^{11}x_1)^{11}x_1\\
\hspace{-0.1cm}=& x_1\Psi_1^{11}x_1 \Psi_2^{11}x_1 + x_1\Psi_1^{10}x_0 \Psi_2^{11}x_1 + x_1\Psi_2^{11}x_1\Psi_1^{11}x_1 + x_1\Psi_2^{11}x_0\Psi_1^{01}x_1 + x_1d_{\Psi_1}(\Psi_2^{11})x_1
\end{aligned}
\end{align*}
holds, we have
\begin{align*}
d_{\Psi_1}(\Psi_2)^{11} = \Psi_1^{11}x_1\Psi_2^{11} + \Psi_2^{11}x_1\Psi_1^{11} + \Psi_1^{10}x_0\Psi_2^{11} + \Psi_2^{11}x_0\Psi_1^{01} + d_{\Psi_1}(\Psi_2^{11}).
\end{align*}

In a similar manner, we have
\begin{align*}
&d_{\Psi_1}(x_0\Psi_2^{01}x_1 + x_1\Psi_2^{11}x_1)\\
\hspace{-0.1cm}=&x_0\Psi_2^{01}x_1\Psi_1^{11}x_1 + x_0\Psi_{2}^{01}x_0\Psi_1^{01}x_1 + x_0\Psi_{1}^{01}x_1\Psi_2^{11}x_1 + x_0\Psi_{1}^{00}x_0\Psi_2^{11}x_1 + x_0d_{\Psi_1}(\Psi_2^{01})x_1
\end{align*}
and
\begin{align*}
&d_{\Psi_1}(x_1\Psi_2^{10}x_0 + x_1\Psi_2^{11}x_1)\\
\hspace{-0.1cm}=&x_1\Psi_1^{11}x_1\Psi_2^{10}x_0 + x_1\Psi_{1}^{10}x_0\Psi_2^{10}x_0 + x_1\Psi_{2}^{11}x_1\Psi_2^{10}x_0 + x_1\Psi_{2}^{11}x_0\Psi_1^{00}x_0 + x_1d_{\Psi_1}(\Psi_2^{10})x_0.
\end{align*}
Thus we have

\begin{align*}
d_{\Psi_1}(\Psi_2)^{01} =& \Psi_2^{01}x_1\Psi_1^{11} + \Psi_{2}^{01}x_0\Psi_1^{01} + \Psi_{1}^{01}x_1\Psi_2^{11} + \Psi_{1}^{00}x_0\Psi_2^{11} + d_{\Psi_1}(\Psi_2^{01})
\end{align*}

and

\begin{align*}
d_{\Psi_1}(\Psi_2)^{10} =&\Psi_1^{11}x_1\Psi_2^{10} + \Psi_{1}^{10}x_0\Psi_2^{10} + \Psi_{2}^{11}x_1\Psi_2^{10} + \Psi_{2}^{11}x_0\Psi_1^{00} + d_{\Psi_1}(\Psi_2^{10}).
\end{align*}

Therefore, we get

\begin{align*}
&d_{\Psi_1}(\Psi_1)^{11} + d_{\Psi_1}(\Psi_2)^{01} + d_{\Psi_1}(\Psi_2)^{10}\\
&\begin{aligned}
=&(\Psi_{1}^{11}x_1\Psi_2^{11} + \Psi_{2}^{11}x_1\Psi_1^{11}) + ( \Psi_1^{01}x_1\Psi_2^{11} + \Psi_2^{01}x_1\Psi_1^{11}) \\
&\quad+(\Psi_1^{11}x_1\Psi_2^{10} + \Psi_2^{11}x_1\Psi_1^{10}) + (\Psi_1^{00}x_0\Psi_2^{11} + \Psi_2^{11}x_0\Psi_1^{00})\\
&\qquad + (\Psi_{1}^{10}x_0\Psi_2^{11} +  \Psi_{2}^{11}x_0\Psi_1^{01} + \Psi_1^{10}x_0\Psi_2^{10} + \Psi_2^{01}x_0 \Psi_1^{01})\\
&\qquad +(d_{\Psi_1}(\Psi_2)^{11} + d_{\Psi_1}(\Psi_2^{01}) + d_{\Psi_1}(\Psi_2^{10})).
\end{aligned}
\end{align*}

Since $\Psi_2 \in V_{\strprty}$, we have

\begin{align*}
 &\Psi_{1}^{10}x_0\Psi_2^{11} +  \Psi_{2}^{11}x_0\Psi_1^{01} + \Psi_1^{10}x_0\Psi_2^{10} + \Psi_2^{01}x_0 \Psi_1^{01}\\
\stackrel{\eqref{eq:str-prty-3cycle}}{=}&-\Psi_{1}^{10}x_0\Psi_2^{10} - \Psi_{1}^{10}x_0\Psi_2^{01}  -  \Psi_{2}^{10}x_0\Psi_1^{01} -  \Psi_{2}^{01}x_0\Psi_1^{01} + \Psi_1^{10}x_0\Psi_2^{10} + \Psi_2^{01}x_0 \Psi_1^{01}\\
=&-\Psi_{1}^{10}x_0\Psi_2^{01}  -  \Psi_{2}^{10}x_0\Psi_1^{01} 
\end{align*}

and 
\begin{align*}
&d_{\Psi_1}(\Psi_2)^{11} + d_{\Psi_1}(\Psi_2^{01})  d_{\Psi_1}(\Psi_2^{10}) \\
=&d_{\Psi_1}(\Psi_2^{11} + \Psi_2^{10} + \Psi_2^{01})\\
\stackrel{\eqref{eq:str-prty-3cycle}}{=}& 0.
\end{align*}

Since $\Psi \in \Fad$ satisfies $\Psi^{00} = 0$, we have
\[
\Psi_1^{00}x_0\Psi_2^{11} + \Psi_2^{11}x_0\Psi_1^{00} = 0.
\]

Therefore, it follows that
\begin{align*}
&d_{\Psi_1}(\Psi_1)^{11} - d_{\Psi_1}(\Psi_2)^{01} - d_{\Psi_1}(\Psi_2)^{10}\\
&\begin{aligned}
=&(\Psi_{1}^{11}x_1\Psi_2^{11} + \Psi_{2}^{11}x_1\Psi_1^{11}) + ( \Psi_1^{01}x_1\Psi_2^{11} + \Psi_2^{01}x_1\Psi_1^{11}) \\
&\quad +(\Psi_1^{11}x_1\Psi_2^{10} + \Psi_2^{11}x_1\Psi_1^{10}) - (\Psi_{1}^{10}x_0\Psi_2^{01}  +  \Psi_{2}^{10}x_0\Psi_1^{01} ).
\end{aligned}
\end{align*}

\end{proof}

\subsection{A Lie algebra associated with adjoint multiple zeta values}

Our second main theorem is the following.
\begin{mt}\label{mt:addmr-Lie}
The $\Q$-linear space $\addmr_0 \cap \Fad \cap V_{\strprty}$ is a Lie subalgebra of $\tm_1$.
\end{mt}

\begin{rmk}
We expect that 
\[
\addmr_0\cap \Fad \subset V_{\strprty}
\]
holds.
Therefore, Main Theorem~\ref{mt:addmr-Lie} means that, by imposing the assumption that $\addmr_0\cap \Fad \subset V_{\strprty}$, we were able to verify the statement of $\addmr_0 \cap \Fad$ described in Remark~\ref{rmk:addmr-bracker}.
\end{rmk}

\section{A proof of Main Theorem \ref{mt:addmr-Lie}}

\subsection{Reduction to Lemma~\ref{lem:essential}}

Let us prove the above Main Theorem \ref{mt:addmr-Lie} in the remaining pages.
The following lemma plays an essential role in the proof of Main Theorem \ref{mt:addmr-Lie}, and its proof is given later,

\begin{lem}\label{lem:essential}
For $\Psi_1$, $\Psi_2 \in \addmr_0 \cap \Fad \cap V_{\strprty}$, we have
\begin{align}
\begin{aligned}
\Delta_*(d_{\Psi_1}(\Psi_2)_{\#}) =&
d_{\Psi_1}(\Psi_2)_{\#} \otimes 1 + 1\otimes d_{\Psi_1}(\Psi_2)_{\#} + \Psi_{1,\#} \otimes  \Psi_{2,\#}  + \Psi_{2,\#} \otimes  \Psi_{1,\#}.  
\end{aligned}  \label{eq:essential}
\end{align}
\end{lem}

\begin{proof}[Proof of Main Theorem \ref{mt:addmr-Lie}]

Since, by the definition of $\addmr_0$ and Proposition \ref{prop:deriv-strprty}, $\addmr_0 \subset \tm_1$ as a $\Q$-linear space and $ (\Fad \cap  V_{\strprty}) \subset \tm_1$ as a Lie algebra, it suffices to show that $\Delta_*(\{\Phi,\Psi\}_{1,\#}) = \{\Phi,\Psi\}_{1,\#} \otimes 1 + 1\otimes \{\Phi,\Psi\}_{1,\#}$.
By \eqref{eq:essential}, we have
\[
\Delta_*(d_{\Psi_1}(\Psi_2)_{\#})  =d_{\Psi_1}(\Psi_2)_{\#} \otimes 1 + 1\otimes d_{\Psi_1}(\Psi_2)_{\#} +  \Psi_{1,\#} \otimes \Psi_{2,\#} + \Psi_{2,\#} \otimes \Psi_{1,\#}.
\]
Therefore, we have
\begin{align*}
\Delta_*(\{\Psi_1,\Psi_2\}_{1,\#}) &= \Delta_*(d_{\Psi_1}(\Psi_2)_{\#} - d_{\Psi_2}(\Psi_1)_{\#} ) \\
&\begin{aligned}
&=d_{\Psi_1}(\Psi_2)_{\#} \otimes 1 + 1\otimes d_{\Psi_1}(\Psi_2)_{\#} +  \Psi_{1,\#} \otimes \Psi_{2,\#} + \Psi_{2,\#} \otimes \Psi_{1,\#}\\
&\qquad-d_{\Psi_2}(\Psi_1)_{\#} \otimes 1 - 1\otimes d_{\Psi_2}(\Psi_1)_{\#}  -  \Psi_{2,\#} \otimes \Psi_{1,\#} + \Psi_{1,\#} \otimes \Psi_{2,\#}\\
\end{aligned}\\
&= d_{\Psi_1}(\Psi_2)_{\#} \otimes 1 + 1\otimes d_{\Psi_1}(\Psi_2)_{\#} - d_{\Psi_1}(\Psi_2)_{\#} \otimes 1 - 1\otimes d_{\Psi_1}(\Psi_2)_{\#}\\
&=\{\Psi_1,\Psi_2\}_1 \otimes1 + 1\otimes \{\Psi_1,\Psi_2\}_1,
\end{align*}
as claimed.
\end{proof}


\subsection{A proof of Lemma \ref{lem:essential}}\label{section:proof of lemma}

In this subsection, we prove Lemma \ref{lem:essential} to complete our main theorem.
We use the commutative power series $\vimo_{\Psi}$ which corresponds to $\Psi \in \ide{h}^{\vee}$.
Notations and operations are due to Ecalle's mould theory (cf. \cite{Ecalle}), but we will not discuss further here.

To begin with, let $U: = \{\u_i \mid i\ge 0\}$ be a set of indeterminates. Put 
\[
\mathcal{G} := \Q + \bigcup_{k\ge 0}\Q[[\u_0,\ldots,\u_k]].
\]

\nsubsection{Commutative power series corresponding to $\phi \in \ide{h}^{\vee}$}

For $\phi \in \ide{h}^{\vee}$, where we suppose that the number of occurrences of $x_1$ is equal to $r\in\mathbb{Z}_{\ge 0}$, we define $\vimo_{\phi}(\u_0,\ldots,\u_r)\in \mathcal{G}$ by
\[
\vimo_{\phi}(\u_0,\u_1,\ldots,\u_r) := \sum_{k_0,\ldots,k_r \ge 0} \langle \phi\mid x_0^{k_r}x_1\cdots x_0^{k_1}x_1x_0^{k_0} \rangle \u_0^{k_0}\cdots \u_r^{k_r}.
\]
Furthermore, we define $\mi_{\phi}$, $\ma_{\phi}$ by
\begin{align*}
\ma_{\phi}^r(\u_1,\ldots,\u_r) :=& \vimo_{\phi}^{r+1}(0,\u_1,\u_1+\u_2,\ldots,\u_1+\cdots + \u_r)\\
\mi_{\phi}^r(\u_1,\ldots,\u_r) :=& \vimo_{\phi}^{r+1}(0,\u_1,\ldots,\u_r).
\end{align*}

Let $\mathcal{U}_{\Z}$ be the $\Z$-module generated by $U$, and $\mathcal{U}_{\Z}^{\bullet}$ be the free monoid generated by $\mathcal{U}_{\Z}$ with the empty word $\emptyset$.
 For any word $u \in \mathcal{U}_{\Z}^{\bullet}$, we define $\deg u$ as the length of $u$. Let $A_{\mathcal{U}} := \Q \langle \mathcal{U}_{\Z}\rangle$.

Let $f \in \mathcal{G}$.
By the following $\Q$-linear map, we consider an element of $A_{\mathcal{U}}$ as the components of some commutative power series:
\[
A_{\mathcal{U}}\ni (\v_1)\cdots (\v_r) \mapsto f(\v_1,\ldots,\v_r) \in A_{\mathcal{U}}
\]
for $(\v_1)\cdots (\v_r) \in \mathcal{U}_{\Z}^{\bullet}$.
Additionally, put $\mathcal{K}:= \Q(\u_i \mid i\ge 0)$ and $A_{\mathcal{U}}^{\text{rat}} := \mathcal{K}\langle \mathcal{U}_{\Z}\rangle$.
If necessary, we consider the coefficient expansion of $\mathcal{K}$, we have a $\mathcal{K}$-linear map:
\[
A_{\mathcal{U}}^{\text{rat}}\ni h\cdot (\v_1)\cdots (\v_r) \mapsto h\cdot f(\v_1,\ldots,\v_r) \in A_{\mathcal{U}}^{\text{rat}}
\]
for $h \in \mathcal{K}$, $(\v_1)\cdots (\v_r) \in \mathcal{U}_{\Z}^{\bullet}$.

\nsubsection{Operations around commutative power series}

\vspace{0.5cm}
\hspace{-0.5cm}\textbf{Shuffle product}
\vspace{0.5cm}

We define a product $\shuffle$ on $A_{\mathcal{U}}$ as a $\Q$-bilinear binary operation given by $\emptyset \shuffle \w := \w\shuffle \emptyset = \w$ and
\[
a\w\shuffle b\v := a(\w\shuffle b\v) + b(a\w\shuffle \v)
\]
for $a$, $b\in \mathcal{U}_{\Z}$ and $\w$, $\v\in \mathcal{U}_{\Z}^{\bullet}$.
Note that a pair $(A_{\mathcal{U}},\shuffle)$ forms a unital, commutative, and associative $\mathbb{Q}$-algebra.

The following properties related to $\ma_{\phi}$ and $\mi_{\phi}$ are well-known.
\begin{prop}[{\cite[(\ref{sh-1-1}) = Lemma 45 (4), (\ref{sh-1-2}) = Proposition 27]{Furusho-Hirose-Komiyama}}]\label{prop:sh-1}
Let $\phi \in \ide{h}^{\vee}$.
\begin{enumerate}
\item \label{sh-1-1}  If $\langle \phi \mid x_1 \shuffle u \rangle = 0$ for $u\in X^*$ except for $1$, then for $r\in\mathbb{Z}_{>0}$, we have
\begin{align}
\vimo_{\phi}^{r+1}(\x_0,\x_1,\ldots,\x_r) = \vimo_{\phi}^r(0,\x_1-\x_0,\ldots,\x_r-\x_0). \label{eq:shift}
\end{align}
\item\label{sh-1-2} If $\langle \phi \mid u\shuffle v\rangle = 0$ for $u$ and $v\in X^{*}\setminus\{1\}$, then 
\[
\ma_{\phi} (\mathbb{X} \shuffle \mathbb{Y}) = 0
\]
holds for $\mathbb{X}$ and $\mathbb{Y} \in \mathcal{U}_{\Z}^{\bullet} \setminus\{\emptyset\}$.
\end{enumerate}
\end{prop}

\vspace{0.5cm}
\hspace{-0.5cm}\textbf{Harmonic product}
\vspace{0.5cm}

We consider another product $*$ on $A_{\mathcal{U}}^{\mathrm{rat}}$.
The $\Q$-bilinear binary operation $*$ on $A_{\mathcal{U}}^{\mathrm{rat}}$ is inductively defined by
\[
w * \emptyset = \emptyset * w = w
\]
and
\[
a\w * b\v := 
\begin{cases}
0 & a= b\\
a(\w * b\v) + b(a\w * b\v) + \frac{1}{a-b}(a(\w*\v) - b(\w * \v)) & a\neq b
\end{cases}
\]
for $a$, $b\in\mathcal{U}_{\Z}$ and $\w$, $\v\in \mathcal{U}_{\Z}^{\bullet}$.

\begin{rmk}
The harmonic product $*$ on $\Q\langle Y\rangle$ corresponds to the harmonic product $*$ on $A_{\mathcal{U}}^{\mathrm{rat}}$.
Namely, for $\Psi \in \ide{h}^{\vee}$, it follows that
\begin{align*}
&\vimo_{\Psi}(0,(\x_1,\ldots,\x_r) * (\x_{r+1},\ldots,\x_{r+l}))\\
=&\sum_{k_1,\ldots,k_{r+l}\ge 0} \langle \Pi_Y(\Psi) \mid (y_{k_{r+l}}\cdots y_{r+1}) * (y_{r} \cdots y_1)\rangle \x_1^{k_1}\cdots \x_{r+l}^{k_{r+l}}.
\end{align*} 
\end{rmk}

Put $z_k := x_0^{k-1}x_1$ for $k\in\mathbb{Z}_{\ge 1}$.
In terms of commutative power series, the interpretation of the harmonic product relations modulo products is as follows:

\begin{prop}[{\cite[Remark 31 (i) and Lemma 127]{Furusho-Hirose-Komiyama}}]\label{har-com}
Let $\Psi \in \ide{h}^{\vee}$. Then the condition $\langle \Pi_Y(\Psi) \mid u * v \rangle = 0$ for $u$, $v \in Y^*\setminus\{1\}$ is equivalent to the condition $\mi_{\Psi} (\mathbb{X} * \mathbb{Y}) = 0$ for $\mathbb{X}$, $\mathbb{Y} \in \mathcal{U}_{\Z}^{\bullet}\setminus\{\emptyset\}$.
\end{prop}

\nsubsection{ $\addmr_0$ via commutative power series}

In this section, we consider $\addmr_0$ via commutative power series. 
Based on Proposition~\ref{prop:sh-1} and Proposition~\ref{har-com}, $\addmr_0$ can be reformulated as follows:

\begin{prop}\label{prop:addmr-mould}
Let $\Psi \in  \ide{h}^{\vee}$. Then $\Psi \in \addmr_0$ if and only if $\Psi$ satisfies the following properties:
\begin{enumerate}
\item\label{enumi:ad--1} $\mi_{\Psi}(\v,0) = 0$ for $\v \in \mathcal{U}_{\Z}$,
\item\label{enumi:ad--2} $\ma_{\Psi}(\mathbb{X} \shuffle \mathbb{Y}) = 0$ for $\mathbb{X}$, $\mathbb{Y} \in \mathcal{U}_{\Z}^{\bullet}\setminus\{\emptyset\}$, and
\item\label{enumi:ad--3} It follows that 
\begin{align}
\vimo_{\Psi}(\z_1, \mathbb{X} * \mathbb{Y} , \z_2) = 0 \label{eq:adjoint-harmonic}
\end{align}
for $\mathbb{X}$, $\mathbb{Y} \in \mathcal{U}_{\Z}^{\bullet}\setminus\{\emptyset\}$ and $\z_1$ and $\z_2 \in \mathcal{U}_{\Z}\setminus\{\emptyset\}$.
\end{enumerate}
\end{prop}

\begin{proof}
Let $\Psi \in \addmr_0$.
The equivalence between $\langle \Psi \mid x_1x_0^{k}x_1\rangle = 0 $ for $k\in\mathbb{Z}_{\ge 0}$ and Condition (\ref{enumi:ad--1}) is clear.

The condition $\Delta_{\shuffle}(\Psi) = \Psi\otimes 1 + 1\otimes \Psi$ corresponds to Condition (\ref{enumi:ad--2}), as shown in  Proposition \ref{prop:sh-1} (\ref{sh-1-2}).

Let us denote $\mathbb{X} = (\x_1,\ldots,\x_r)$, $\mathbb{Y} = (\y_1,\ldots,\y_l)$.
The Condition (\ref{enumi:ad--3}) is equivalent to Proposition \ref{def:addmr} (\ref{def:addmr--3}) and Proposition \ref{prop:sh-1} (\ref{sh-1-1}). In fact, Proposition~\ref{def:addmr} (\ref{def:addmr--3}) is equivalent to
\begin{align*}
& \langle \mathbf{q}_{\#}^{\vee} (\Psi) \mid y_{k_1}\cdots y_{k_r} * y_{s_1}\cdots y_{s_t} \rangle = 0\\
\Longleftrightarrow&\sum_{l\ge 1} T^l \langle \mathbf{q}^{\vee}(\Psi) \mid y_l(y_{k_1}\cdots y_{k_r} * y_{s_1}\cdots y_{s_t})\rangle = 0\\
\Longleftrightarrow&  \langle \mathbf{q}^{\vee}(\Psi) \mid y_l(y_{k_1}\cdots y_{k_r} * y_{s_1}\cdots y_{s_t})\rangle = 0 \hspace{1ex} (l\in\mathbb{Z}_{\ge 1})
\end{align*}
for $k_1,\ldots,k_r, s_1,\ldots,s_t \in\mathbb{Z}_{>0}$, which implies that
\begin{align*}
&\vimo_{\Psi}(0,\mathbb{X} * \mathbb{Y},\z)\\
=&\sum_{\substack{l\ge 1 \\ k_1,\ldots,k_r\ge 1 \\ s_1,\ldots,s_t\ge 1}}\langle \mathbf{q}^{\vee}(\Psi) \mid y_l(y_{k_1}\cdots y_{k_r} * y_{s_1}\cdots y_{s_t})\rangle \y_{t}^{s_t-1}\cdots \y_{1}^{s_1-1}\x_r^{k_r-1}\cdots \x_1^{k_1-1}\z^{l-1} \\
=&0.
\end{align*}
Because of Proposition \ref{prop:sh-1} (\ref{sh-1-1}), 
\[
\vimo_{\Psi}(0,(\y_{s_t},\cdots,\y_{s_1}) * (\x_r,\cdots,\x_1),\z) = 0
\]
 turns out to be 
\[
\vimo_{\Psi}(\z_1,(\y_{s_t},\cdots,\y_{s_1}) * (\x_r,\cdots,\x_1),\z_2) = 0.
\]
This completes the proof.

\end{proof}

\begin{prop}\label{prop:adjoint-mould}
Let $\Psi \in \ide{h}^{\vee}$ with $\langle \Psi \mid w \rangle = 0$ for $w\in X^*$ with $\wt w \le 1$.
Then, $\Psi^{00} = 0$  if and only if 
\begin{align}
\begin{aligned}
&\vimo_{\Psi}(\x_0,\ldots,\x_r)\\ 
\hspace{1cm}&= \vimo_{\Psi}(\x_0,\ldots,\x_{r-1},0) + \vimo_{\Psi}(0,\x_1,\ldots,\x_r) - \vimo_{\Psi}(0,\x_1,\ldots,\x_{r-1},0)
\end{aligned}
 \label{eq:adjoint-mould}
\end{align}
 for $(\x_0,\ldots,\x_r) \in \mathcal{U}_{\Z}^{\bullet}$.
\end{prop}

\begin{proof}
Since $\Psi^{00} = 0$ is equivalent to $\langle \Psi \mid x_0^{k_0}x_1\cdots x_1x_0^{k_r} \rangle = 0$ for $(k_0,\ldots,k_r) \in \mathbb{Z}_{\ge 0}^{r+1}$ with $k_0,k_r > 0$. Therefore,
\begin{align*}
&\vimo_{\Psi}(\x_0,\ldots,\x_r)\\
&\begin{aligned}
=& \sum_{k_0,\ldots,k_r\ge 0} \langle \Psi \mid x_0^{k_r}x_1x_0^{r-1}\cdots x_0^{k_1}x_1x_0^{k_0}\rangle\\
\end{aligned}\\
&\begin{aligned}
=&\bigg(\sum_{k_1,\ldots,k_r\ge 0} + \sum_{k_0,\ldots,k_{r-1}\ge 0} - \sum_{k_1,\ldots,k_{r-1}\ge 0} \bigg) \langle \Psi \mid x_0^{k_r}x_1x_0^{k_{r-1}}\cdots x_0^{k_1}x_1x_0^{k_0}\rangle\\
\end{aligned}\\
&\begin{aligned}
=&\vimo_{\Psi}(\x_0,\ldots,\x_{r-1},0) + \vimo_{\Psi}(0,\x_1,\ldots,\x_r) - \vimo_{\Psi}(0,\x_1,\ldots,\x_{r-1},0)
\end{aligned}
\end{align*}
as claimed.
\end{proof}

\begin{cor}\label{prop:adjoint-mould}
If $\Psi \in \Fad$, then \eqref{eq:adjoint-mould} holds.
\end{cor}
\begin{proof}
Immediately, this condition follows from the definition of $\Fad$.
\end{proof}


\begin{prop}\label{prop:strong-parity}
Let $\Psi \in \ide{h}^{\vee}$. Then $\Psi \in V_{\strprty}$ is equivalent to
\begin{align}
\begin{aligned}
&\vimo_{\Psi}(0,\x_1,\ldots,\x_r,0)\\
=& -\frac{\vimo_{\Psi }(0,\x_1,\ldots,\x_r) - \vimo_{\Psi }(0,\x_1,\ldots,\x_{r-1},0)}{\x_r} \\
&- \frac{\vimo_{\Psi }(\x_1,\ldots,\x_r,0) - \vimo_{\Psi }(0,\x_2,\ldots,\x_{r},0)}{\x_1}
\end{aligned}
 \label{eq:strong-parity}
\end{align}
for $r\in\mathbb{Z}_{>0}$.
\end{prop}

\begin{proof}

Since $\Psi \in V_{\strprty}$ satisfies
\begin{align*}
& \langle \Psi \mid x_1x_0^{k_r}x_1\cdots x_0^{k_1}x_1 \rangle\\
=&-  \langle \Psi \mid x_0^{k_r+1}x_1x_0^{k_{r-1}}\cdots x_0^{k_1}x_1 \rangle - \langle \Psi \mid x_1x_0^{k_r}\cdots x_0^{k_2}x_1x_0^{k_1+1}\rangle,
\end{align*}
for $k_1,\ldots,k_r \in \mathbb{Z}_{\ge 0}$, we have

\begin{align*}
&\vimo_{\Psi}(0,\x_1,\ldots,\x_r,0)\\
=&\sum_{k_1,\ldots,k_r\ge 0} \x_1^{k_1}\cdots \x_r^{k_r} \langle \Psi \mid  x_1x_0^{k_r}x_1\cdots x_0^{k_1}x_1 \rangle \\
=&- \sum_{k_1,\ldots,k_r\ge 0}  \x_1^{k_1}\cdots \x_r^{k_r} \bigg( \langle \Psi \mid x_0^{k_r+1}x_1x_0^{k_{r-1}}\cdots x_0^{k_1}x_1 \rangle + \langle \Psi \mid x_1x_0^{k_r}\cdots x_0^{k_2}x_1x_0^{k_1+1}\rangle  \bigg)\\
=& -\frac{\vimo_{\Psi }(0,\x_1,\ldots,\x_r) - \vimo_{\Psi}(0,\x_1,\ldots,\x_{r-1},0)}{\x_r} \\
&- \frac{\vimo_{\Psi }(\x_1,\ldots,\x_r,0) - \vimo_{\Psi }(0,\x_2,\ldots,\x_{r},0)}{\x_1}
\end{align*}
as claimed.

\end{proof}

\begin{cor}
If $\Psi \in \ide{F}_2 \cap V_{\strprty}$, we have
\begin{align}
\begin{aligned}
&\vimo_{\Psi}(\y,\x_1,\ldots,\x_r,\y)\\
&\begin{aligned}
=& -\frac{\vimo_{\Psi }(0,\x_1,\ldots,\x_r) - \vimo_{\Psi }(0,\x_1,\ldots,\x_{r-1},\y)}{\x_r - \y}\\
&\qquad - \frac{\vimo_{\Psi }(\x_1,\ldots,\x_r,0) - \vimo_{\Psi }(\y,\x_2,\ldots,\x_{r},0)}{\x_1 - \y}.
\end{aligned}
\end{aligned} \label{eq:strong-parity-revised}
\end{align}
\end{cor}

\begin{proof}
It follows that
\begin{align*}
&\vimo_{\Psi}(\y,\x_1,\ldots,\x_r,\y) \\
\stackrel{\eqref{eq:shift}}{=} & \vimo_{\Psi}(0,\x_1-\y,\ldots,\x_r-\y,0)\\
\stackrel{\eqref{eq:strong-parity}}{=}& -\frac{\vimo_{\Psi }(0,\x_1-\y,\ldots,\x_r-\y) - \vimo_{\Psi }(0,\x_1-\y,\ldots,\x_{r-1}-\y,0)}{\x_r-\y} \\
&- \frac{\vimo_{\Psi }(\x_1-\y,\ldots,\x_r-\y,0) - \vimo_{\Psi }(0,\x_2-\y,\ldots,\x_{r}-\y,0)}{\x_1-\y}\\
\stackrel{\eqref{eq:shift}}{=} & -\frac{\vimo_{\Psi }(\y,\x_1,\ldots,\x_r) - \vimo_{\Psi }(\y,\x_1,\ldots,\x_{r-1},\y)}{\x_r-\y} \\
&- \frac{\vimo_{\Psi }(\x_1,\ldots,\x_r,\y) - \vimo_{\Psi }(\y,\x_2,\ldots,\x_{r},\y)}{\x_1-\y}\\
\stackrel{\eqref{eq:adjoint-mould}}{=} &  -\frac{\vimo_{\Psi }(0,\x_1,\ldots,\x_r) + \vimo_{\Psi }(\y,\x_1,\ldots,\x_{r-1},0) - \vimo_{\Psi }(0,\x_1,\ldots,\x_{r-1},0)}{\x_r-\y}\\
& + \frac{\vimo_{\Psi }(0,\x_1,\ldots,\x_{r-1},\y) + \vimo_{\Psi }(\y,\x_1,\ldots,\x_{r-1},0) - \vimo_{\Psi }(0,\x_1,\ldots,\x_{r-1},0)}{\x_r-\y}\\
&- \frac{\vimo_{\Psi }(0,\x_2,\ldots,\x_r,\y) + \vimo_{\Psi }(\x_1,\ldots,\x_r,0) - \vimo_{\Psi }(0,\x_2,\ldots,\x_r,0)}{\x_1-\y}\\
&+ \frac{\vimo_{\Psi }(0,\x_2,\ldots,\x_{r},\y) + \vimo_{\Psi }(\y,\x_2,\ldots,\x_{r},0) - \vimo_{\Psi }(0,\x_2,\ldots,\x_{r},0) }{\x_1-\y}\\
=& -\frac{\vimo_{\Psi }(0,\x_1,\ldots,\x_r) - \vimo_{\Psi }(0,\x_1,\ldots,\x_{r-1},\y)}{\x_r - \y} \\
&- \frac{\vimo_{\Psi }(\x_1,\ldots,\x_r,0) - \vimo_{\Psi }(\y,\x_2,\ldots,\x_{r},0)}{\x_1 - \y}
\end{align*}
as claimed.
\end{proof}


\nsubsection{Preliminaries}

To prove Lemma \ref{lem:essential}, we use special conventions for some variables.
Hereafter, unless otherwise stated, capital bold letters denote elements of $\mathcal{U}_{\Z}^{\bullet}$, lowercase script letters denote elements of $\mathcal{U}_{\Z} \setminus\{\emptyset\}$, and lowercase Fraktur letters denote elements of $\mathcal{U}_{\Z}$.

First, we consider $d_{\Psi_1}(\Psi_2)$ via commutative power series for $\Psi_1$, $\Psi_2 \in \tm_1$.
Let $\mathbb{X}$, $\mathbb{Y} \in \mathcal{U}_{\Z}^{\bullet}$.

\begin{prop}\label{prop:deriv-explicit}
For $\Psi_1$, $\Psi_2\in  \ide{h}^{\vee}$, it follows that
\[
\vimo_{d_{\Psi_1}(\Psi_2)}(\x_0,\ldots,\x_r) = \sum_{0\le a < b \le r} \vimo_{\Psi_1}(\x_a,\ldots,\x_b) \vimo_{\Psi_2}(\x_0,\ldots,\x_a,\x_b,\ldots,\x_r).
\]
In other words, we have
\begin{align}
\begin{aligned}
&\vimo_{d_{\Psi_1}(\Psi_2)}(\z,\mathbb{X},\z')\\
& =  \sum_{\mathbb{V}_L(\v)\mathbb{V}_M(\v')\mathbb{V}_R = (\z)\cdot \mathbb{X} \cdot (\z') } f_1(\v,\mathbb{V}_M,\v')f_2(\delta_{\z}(\v),(\mathbb{V}_L\cdot (\v)\cdot(\v') \cdot \mathbb{V}_R),\delta_{\z'}(\v')),
\end{aligned}  \label{eq:deriv-explicit}
\end{align}
where, for any letters $a$, $b \in \mathcal{U}_{\Z}$, define
\[
\delta_{a}(b):=
\begin{cases}
a & a \neq b\\
\emptyset & a = b.
\end{cases}
\]
\end{prop}

\begin{proof}
Let $k\in\mathbb{Z}_{>0}$, $r\in\mathbb{Z}_{\ge k}$, and $s_0,\ldots,s_k\in\mathbb{Z}_{\ge 0}$. Since the following
{\small
\begin{align*}
&\sum_{n_0,n_1,\ldots,n_r\ge 0} \langle d_{\Psi_1}(x_0^{s_k}x_1\cdots x_1x_0^{s_0}) \mid x_0^{n_r}x_1\cdots x_0^{n_1}x_1x_0^{n_0} \rangle \x_0^{n_0}\cdots \x_r^{n_r}\\
\displaybreak[3]
=&\sum_{n_0,n_1,\ldots,n_r\ge 0}\sum_{j = 1}^k \langle  x_0^{s_k}x_1\cdots x_0^{s_j}\Psi x_0^{s_{j-1}} \cdots x_0^{s_1}x_1x_0^{s_0} \mid x_0^{n_r}x_1\cdots x_0^{n_1}x_1x_0^{n_0} \rangle \x_0^{n_0}\cdots \x_r^{n_r} \\
\displaybreak[3]
=&\sum_{j = 1}^r \x_0^{s_0}\cdots \x_{j-2}^{s_{j-2}}\x_{r + j -k +1}^{s_{j+1}} \cdots \x_r^{s_k} \times \\
& \sum_{\substack{n_j,\ldots,n_{r+j-k-1}\ge 0 \\ n_{j-1}\ge s_{j-1}, n_{r+j-k}\ge s_j}} \langle \Psi \mid x_0^{n_{r + j -k }-s_{j}}x_1x_0^{n_{r+j-k-1}} \cdots x_0^{n_j}x_1x_0^{n_{j-1} - s_{j-1}}\rangle
\x_{j-1}^{n_{j-1} - s_{j-1}} \x_j^{n_j} \cdots \x_{r + j -k -1}^{n_{r+j-k-1}} \x_{r + j -k}^{n_{r + j -k }-s_{j}}  \\
\displaybreak[3]
=&\sum_{j = 1}^r \x_0^{s_0}\cdots \x_{j-2}^{s_{j-2}}\x_{j-1}^{s_{j-1}} \cdot \x_{r+j-k}^{s_j}\x_{r + j -k +1}^{s_{j+1}} \cdots \x_r^{s_k} \\
&\times \sum_{n_{j-1},\ldots,n_{r+j-k}\ge 0 } \langle \Psi \mid x_0^{n_{r + j -k }}x_1x_0^{n_{r+j-k-1}} \cdots x_0^{n_j}x_1x_0^{n_{j-1} }\rangle \x_{j-1}^{n_{j-1} } \x_j^{n_j} \cdots \x_{r + j -k -1}^{n_{r+j-k-1}} \x_{r + j -k}^{n_{r + j -k }}  \\
\displaybreak[3]
=& \sum_{j = 1}^r \x_0^{s_0}\cdots \x_{j-2}^{s_{j-2}}\x_{j-1}^{s_{j-1}} \cdot \x_{r+j-k}^{s_j}\x_{r + j -k +1}^{s_{j+1}} \cdots \x_r^{s_k} \times\vimo_{\Psi_1}(\x_{j-1},\x_j,\ldots,\x_{r+j-k-1},\x_{r+j-k})
\end{align*}
}
holds, we have
{\small
\begin{align*}
&\vimo_{d_{\Psi_1}(\Psi_2)}(\x_0,\ldots,\x_r) \displaybreak[3] \\
=&\sum_{n_0,\ldots,n_r\ge 0} \sum_{k = 1}^r \sum_{s_0,\ldots,s_j\ge 0} \langle \Psi_2\mid x_0^{s_k}x_1\cdots x_1x_0^{s_0} \rangle \langle d_{\Psi_1}(x_0^{s_k}x_1\cdots x_0^{s_1}x_1x_0^{s_0})\mid x_0^{n_r}x_1\cdots x_0^{n_1}x_1x_0^{n_0}\rangle \x_0^{n_0}\cdots \x_r^{n_r} \displaybreak[3] \\
=& \sum_{k = 1}^r \sum_{s_0,\ldots,s_j\ge 0} \langle \Psi_2\mid x_0^{s_k}x_1\cdots x_1x_0^{s_0} \rangle \\
&\times  \sum_{j = 1}^r \x_0^{s_0}\cdots \x_{j-2}^{s_{j-2}}\x_{j-1}^{s_{j-1}} \cdot \x_{r+j-k}^{s_j}\x_{r + j -k +1}^{s_{j+1}} \cdots \x_r^{s_k} \times \vimo_{\Psi_1}(\x_{j-1},\x_j,\ldots,\x_{r+j-k-1},\x_{r+j-k})\\
=&\sum_{j = 1}^r\sum_{k = 1}^r \vimo_{\Psi_1}(\x_{j-1},\x_j,\ldots,\x_{r+j-k-1},\x_{r+j-k}) \vimo_{\Psi_2}(\x_0,\ldots,\x_{j-1},\x_{r+j-k},\ldots,\x_r) \displaybreak[3] \\
=&\sum_{0\le a < b \le r} \vimo_{\Psi_1}(\x_{a},\x_{a+1},\ldots,\x_{b-1},\x_{b}) \vimo_{\Psi_2}(\x_0,\ldots,\x_{a},\x_{b},\ldots,\x_r).
\end{align*}
}
\end{proof}


For $\x \in \mathcal{U}_{\Z}\setminus\{\emptyset\}$ and $\ide{y} \in \mathcal{U}_{\mathbb{Z}}$, we define
\[
c(\x, \ide{y}) := 
\begin{cases}
1 & \ide{y} = \emptyset \\
0& \ide{y} =  x \\
\frac{1}{x - \ide{\eta}} & \text{otherwise.}
\end{cases}
\]

\begin{lem}\label{lem:har-inductive}
Let $\mathbb{X} = \mathbb{X}_L (\x) \mathbb{X}_R$, $\mathbb{Y} \in \mathcal{U}_{\mathbb{Z}}^{\bullet}$ with distinct letters.
Then we have
\begin{align}
\begin{aligned}
& \mathbb{X}_L (\x) \mathbb{X}_R * \mathbb{Y}\\
=& \sum_{\mathbb{Y}_L\cdot (\ide{y})\cdot \mathbb{Y}_R =  \mathbb{Y}} c(\x,\ide{y}) (\mathbb{X}_L * \mathbb{Y}_L) \cdot((\x) - (\ide{y}))\cdot (\mathbb{X}_R * \mathbb{Y}_R).
\end{aligned} \label{eq:harmonic-composition}
\end{align}
\end{lem}

\begin{proof}
Let $\mathbb{Y} := (\y_1)\cdots (\y_l)$.
We prove our claim by induction on $\deg\mathbb{X}_L$.
If $\deg \mathbb{X}_L = 0$, that is, $\mathbb{X}_L = \emptyset$, the definition of the harmonic product implies that
\begin{align*}
(\x) \mathbb{X}_R * \mathbb{Y} =& (\x) (\mathbb{X}_R * \mathbb{Y}) + \frac{1}{\x - \y_1} ((\x) - (\y_1)) (\mathbb{X}_R * ((\y_2)\cdots (\y_l)) ) \\
&+(\y_1) ((\x) (\mathbb{X}_R) * ((\y_2)\cdots (\y_l))) \displaybreak[3] \\
=&\\
\vdots&\\
=&\sum_{j =0}^l (\y_1)\cdots(\y_j)(\x) (\mathbb{X}_R *( (\y_{j+1})\cdots(\y_l))) \\
&+ \sum_{j = 1}^l  \frac{1}{\x - \y_j} (\y_1)\cdots(\y_{j-1})((\x) - (\y_j)) (\mathbb{X}_R *( (\y_{j+1})\cdots(\y_l))) \displaybreak[3] \\
=& \sum_{\mathbb{Y}_L\cdot (\ide{y})\cdot \mathbb{Y}_R =  \mathbb{Y}} c(\x,\ide{y}) \mathbb{Y}_L \cdot((\x) - (\ide{y}))\cdot (\mathbb{X}_R * \mathbb{Y}_R). \displaybreak[3]
\end{align*}
Assume that $\deg \mathbb{X}_L>0$ and our claim holds for all positive integers less than $\deg \mathbb{X}_L$.
Put $\mathbb{X}_L = (\x_1)\mathbb{X}_{M}$. By the induction hypothesis, we have
\begin{align*}
&((\x_1)\mathbb{X}_{M}(\x)\mathbb{X}_R) * \mathbb{Y}\\
\stackrel{\eqref{eq:harmonic-composition}}{=}&\sum_{\mathbb{Y}_L(\ide{y}') \mathbb{Y}_R = \mathbb{Y}} c(\x_1,\ide{y}') \mathbb{Y}_L((\x_1) - (\ide{y}')) ((\mathbb{X}_{M}(\x)\mathbb{X}_R) * \mathbb{Y}_R).
\end{align*}
Since $\deg \mathbb{X}_M < \deg \mathbb{X}_L$, applying the induction hypothesis implies that
\begin{align*}
&(\mathbb{X}_{M}(\x)\mathbb{X}_R) * \mathbb{Y}_R\\
=&\sum_{\mathbb{Y}_M (\ide{y})\mathbb{Y}_{RR} = \mathbb{Y}_R} c(\x,\ide{y})(\mathbb{Y}_M * \mathbb{X}_{M})((\x) - (\ide{y})) (\mathbb{X}_{R} * \mathbb{Y}_{RR}).
\end{align*}
Therefore, we have
\begin{align*}
&((\x_1)\mathbb{X}_{M}(\x)\mathbb{X}_R) * \mathbb{Y} \displaybreak[3] \\
&\begin{aligned}
=&\sum_{\substack {\mathbb{Y}_L(\ide{y}') \mathbb{Y}_R = \mathbb{Y} \\ \mathbb{Y}_M (\ide{y})\mathbb{Y}_{RR} = \mathbb{Y}_R}} c(\x_1,\ide{y}') c(\x,\ide{y}) \mathbb{Y}_L((\x_1) - (\ide{y}'))(\mathbb{Y}_M * \mathbb{X}_{M})((\x) - (\ide{y})) (\mathbb{X}_{R} * \mathbb{Y}_{RR}) 
\end{aligned}\\
\displaybreak[3] \\
&\begin{aligned}
=&\sum_{\mathbb{Y}_L(\ide{y})\mathbb{Y}_R = \mathbb{Y}}c(\x,\ide{y})\Bigg\{ \sum_{\mathbb{Y}_{LL}(\ide{y}')\mathbb{Y}_M = \mathbb{Y}_L}c(\x_1,\ide{y}') \mathbb{Y}_{LL} ((\x_1) - (\ide{y}')) (\mathbb{X}_M * \mathbb{Y}_M) \Bigg\} ((\x) - (\ide{y})) (\mathbb{X}_{R} * \mathbb{Y}_{R}) 
\end{aligned} \displaybreak[3]\\
&
\begin{aligned}\hspace{-0.1cm}\stackrel{\eqref{eq:harmonic-composition}}{=}&
\sum_{\mathbb{Y}_L(\ide{y})\mathbb{Y}_R = \mathbb{Y}} c(\x,\ide{y})  ((\x_1 \mathbb{X}_M) * \mathbb{Y}_L) ((\x) - (\ide{y}))  (\mathbb{X}_{R} * \mathbb{Y}_{R})
\end{aligned}\\
&\hspace{-0.1cm}=\sum_{\mathbb{Y}_L(\ide{y})\mathbb{Y}_R = \mathbb{Y}} c(\x,\ide{y}) ( \mathbb{X}_L * \mathbb{Y}_L) ((\x) - (\ide{y}))  (\mathbb{X}_{R} * \mathbb{Y}_{R})
\end{align*}
as claimed.
\end{proof}

\begin{lem}\label{lem:-har-expl}
For $\Psi_1$, $\Psi_2 \in  \ide{h}^{\vee}$, we have
\[
\vimo_{d_{\Psi_1}(\Psi_2)}(\z,\,\mathbb{X} * \mathbb{Y},\,\z')
= \sum_{(\V,\V')\in S} A_{\V\V'},
\]
where $S:=\{(\z,\z'),(\z,\X'),(\z,\Y'),(\X,\z'),(\Y,\z'),(\X,\X'),(\X,\Y'),(\Y,\X'),(\Y,\Y')\}$, and we define

\begin{align*}
&\begin{aligned}
A_{\z\z'} =&  f_1(\z, \mathbb{X} * \mathbb{Y},\z')f_2(\z,\z')
\end{aligned}\displaybreak[3]\\
&\begin{aligned}
A_{\z X'} := &\sum_{\substack{\mathbb{X}_L (\x) \mathbb{X}_R = \mathbb{X} \\  \mathbb{Y}_L(\ide{y}) \mathbb{Y}_R = \mathbb{Y}}}
c(\x,\ide{y}) f_1(\z,\mathbb{X}_L * \mathbb{Y}_L, \x)f_2(\z,\x,\mathbb{X}_R * \mathbb{Y}_R,\z') 
\end{aligned}\displaybreak[3]\\
&\begin{aligned}
A_{\z Y'} := &\sum_{\substack{\mathbb{X}_L (\ide{x}) \mathbb{X}_R = \mathbb{X} \\  \mathbb{Y}_L (\y) \mathbb{Y}_R = \mathbb{Y}}}
c(\y, \ide{x}) f_1(\z,\mathbb{X}_L * \mathbb{Y}_L, \y)f_2(\z,\y,\mathbb{X}_R * \mathbb{Y}_R,\z')
\end{aligned}\displaybreak[3]\\
&\begin{aligned}
A_{X\z'} :=&\sum_{\substack{\mathbb{X}_L(x)\mathbb{X}_R = \mathbb{X} \\ \mathbb{Y}_L(\ide{y})\mathbb{Y}_R = \mathbb{Y}}}c(\x,\ide{y}) f_1(\x,\mathbb{X}_R * \mathbb{Y}_R, \z')f_2(\z,\mathbb{X}_L * \mathbb{Y}_L, \x ,\z') 
\end{aligned}\displaybreak[3]\\
&\begin{aligned}
A_{Y\z'} :=&\sum_{\substack{\mathbb{X}_L(\ide{x})\mathbb{X}_R = \mathbb{X} \\ \mathbb{Y}_L(\y)\mathbb{Y}_R = \mathbb{Y}}}c(\y,\ide{x}) f_1(\y,\mathbb{X}_R * \mathbb{Y}_R, \z')f_2(\z,\mathbb{X}_L * \mathbb{Y}_L, \y,\z') 
\end{aligned}\displaybreak[3]\\
&\begin{aligned}
A_{XX'} :=& \sum_{\substack{\mathbb{X}_L (\x) \mathbb{X}_M(\x')\mathbb{X}_R = \mathbb{X} \\ \mathbb{Y}_L (\ide{y}) \mathbb{Y}_M(\ide{y}')\mathbb{Y}_R = \mathbb{Y} }} c(\x,\ide{y})c(\x',\ide{y}')f_1(\x,\mathbb{X}_M * \mathbb{Y}_M, \x')\\
&\hspace{5cm}\times f_2(\z, \mathbb{X}_L * \mathbb{Y}_L , \x, \x', \mathbb{X}_R * \mathbb{Y}_R , \z')
\end{aligned}\displaybreak[3]\\
&\begin{aligned}
A_{XY'} :=& \sum_{\substack{\mathbb{X}_L (\x) \mathbb{X}_M(\ide{x}')\mathbb{X}_R = \mathbb{X} \\ \mathbb{Y}_L (\ide{y}) \mathbb{Y}_M(\y')\mathbb{Y}_R = \mathbb{Y} }} c(\x,\ide{y})c(\y',\ide{x}')f_1(\x,\mathbb{X}_M * \mathbb{Y}_M, \y')\\
&\hspace{5cm} \times f_2(\z, \mathbb{X}_L * \mathbb{Y}_L , \x, \y', \mathbb{X}_R * \mathbb{Y}_R , \z')
\end{aligned}\displaybreak[3]\\
&\begin{aligned}
A_{YX'} :=& \sum_{\substack{\mathbb{X}_L (\ide{x}) \mathbb{X}_M(\x')\mathbb{X}_R = \mathbb{X} \\ \mathbb{Y}_L (\y) \mathbb{Y}_M(\ide{y}')\mathbb{Y}_R = \mathbb{Y} }} c(\y,\ide{x})c(\x',\ide{y}')f_1(\y,\mathbb{X}_M * \mathbb{Y}_M, \x')\\
&\hspace{5cm} \times f_2(\z, \mathbb{X}_L * \mathbb{Y}_L , \y, \x', \mathbb{X}_R * \mathbb{Y}_R , \z')
\end{aligned}\displaybreak[3]\\
&\begin{aligned}
A_{YY'} :=& \sum_{\substack{\mathbb{X}_L (\ide{x}) \mathbb{X}_M(\ide{x}')\mathbb{X}_R = \mathbb{X} \\ \mathbb{Y}_L (\y) \mathbb{Y}_M(\y')\mathbb{Y}_R = \mathbb{Y} }} c(\y,\ide{x})c(\y',\ide{x}')f_1(\y,\mathbb{X}_M * \mathbb{Y}_M, \y')\\
&\hspace{5cm}\times f_2(\z, \mathbb{X}_L * \mathbb{Y}_L , \y, \y', \mathbb{X}_R * \mathbb{Y}_R , \z').
\end{aligned}
\end{align*}

\end{lem}

\begin{proof}

By Proposition \ref{prop:deriv-explicit}, $\vimo_{d_{\Psi_1}(\Psi_2)}(\z,\mathbb{X} * \mathbb{Y},\z')$ is decomposed as
\begin{align*}
&\vimo_{d_{\Psi_1}(\Psi_2)}(\z,\mathbb{X} * \mathbb{Y},\z')\\
=& \sum_{\substack{\mathbb{V}_L (\v) \mathbb{V}_M (\v') \mathbb{V}_R \\ = (\z) (\mathbb{X} * \mathbb{Y})(\z')}}
f_1(\v,\mathbb{V}_M, \v')f_2(\delta_{\z}(\v), \mathbb{V}_L, (\v), (\v'), \mathbb{V}_R,\delta_{\z'}(\v')).
\end{align*}

We divide it into cases based on the values of $\v$ and $\v'$. 
First, we can easily see
\[
f_1(\z, \mathbb{X} * \mathbb{Y},\z')f_2(\z,\z') = A_{\z\z'}
\]
for the case where $\v = \z$, $\v' = \z'$.
For the remaining cases, we shall show the case where $\v = \z$ and $\v'$ appears in $\mathbb{X}$ is equal to $A_{\z\X'}$ and the case where $\v$ and $\v'$ appear in $\mathbb{X}$ is equal to $A_{\X\X'}$.
It suffices to demonstrate these two cases because of the discussion of the symmetrality. Namely,
\[
\xymatrix{
  A_{\z\X'} \ar[d]_{\mathbb{X}\leftrightarrow \mathbb{Y}}
             \ar[r]^{\substack{\text{counterwise}\\ \text{discussion}}}
  & A_{\X\z'} \ar[d]^{\mathbb{X}\leftrightarrow \mathbb{Y}} \ar[l] \\
  A_{\z\Y'} \ar[u]
            \ar[r]_{\substack{\text{counterwise}\\ \text{discussion}}}
  & A_{\Y\z'} \ar[u] \ar[l]
}
\hspace{1.0cm}
\xymatrix{
  A_{\X\X'} \ar[r]^{\mathbb{X}\leftrightarrow \mathbb{Y}}
            \ar[d]_{\substack{\mathbb{X}\leftrightarrow \mathbb{Y}\\ \text{for } \x',\,\y'}}
  & A_{\Y\Y'} \ar[l]
              \ar[d]^{\substack{\mathbb{X}\leftrightarrow \mathbb{Y}\\ \text{for } \x',\,\y'}} \\
  A_{\X\Y'} \ar[u]
            \ar[r]_{\substack{\text{counterwise}\\ \text{discussion}}}
  & A_{\Y\X'} \ar[u] \ar[l]
}
\]

\begin{itemize}
\item First, consider the case $\v=\z$ and $\v'$ occurs in $\mathbb{X}$.  
Write $\mathbb{X}=\mathbb{X}_L(\x')\mathbb{X}_R$.  
By Lemma \ref{lem:har-inductive} and Proposition \ref{prop:deriv-explicit},
\begin{align*}
&\vimo_{d_{\Psi_1}(\Psi_2)}(\z, \mathbb{X} * \mathbb{Y},\z') \displaybreak[3]  \\
& \stackrel{\eqref{eq:harmonic-composition}}{=} \sum_{\mathbb{Y}_L(\ide{y}) \mathbb{Y}_R = \mathbb{Y}}c(\x',\ide{y})\vimo_{d_{\Psi_1}(\Psi_2)}((\z)(\mathbb{X}_L * \mathbb{Y}_L)((\x') - (\ide{y})) (\mathbb{X}_R * \mathbb{Y}_R)(\z')) \displaybreak[3]  \\
& 
\begin{aligned}
 & \stackrel{\eqref{eq:deriv-explicit}}{=} \sum_{\mathbb{Y}_L(\ide{y}) \mathbb{Y}_R = \mathbb{Y}} c(\x',\ide{y}) \sum_{\substack{\mathbb{V}_L (\v) \mathbb{V}_M (\v') \mathbb{V}_R \\ = (\z)(\mathbb{X}_L * \mathbb{Y}_L)((\x') - (\ide{y})) (\mathbb{X}_R * \mathbb{Y}_R)(\z')  }} f_1(\v,\mathbb{V}_M, \v')f_2(\delta_{\z}(\v), \mathbb{V}_L, (\v), (\v'), \mathbb{V}_R,\delta_{\z'}(\v')) 
\end{aligned}
\displaybreak[3]  \\
&=\sum_{\mathbb{Y}_L(\ide{y}) \mathbb{Y}_R = \mathbb{Y}} c(\x',\ide{y}) f_1(\z,\mathbb{X}_L * \mathbb{Y}_L, \x')f_2(\z,\x',\mathbb{X}_R * \mathbb{Y}_R,\z')\\
&\qquad + (\text{the cases where $(\v,\v')\neq (\z,\x')$}).
\end{align*}
Summing over all choices of $\x'$ in the decomposition $\mathbb{X}=\mathbb{X}_L(\x')\mathbb{X}_R$ yields that
\begin{align*}
&\vimo_{d_{\Psi_1}(\Psi_2)}(\z,\mathbb{X}*\mathbb{Y},\z')\\
=&
\sum_{\substack{\mathbb{X}_L(\x')\,\mathbb{X}_R=\mathbb{X}\\[2pt]
\mathbb{Y}_L(\ide{y})\,\mathbb{Y}_R=\mathbb{Y}}}
c(\x',\ide{y})\,
f_1(\z,\mathbb{X}_L*\mathbb{Y}_L,\x')\,
f_2(\z,\x',\mathbb{X}_R*\mathbb{Y}_R,\z')\\
&\quad + \text{(all remaining terms with $(\v,\v')\neq(\z,\x')$ for every $\x'$ occurring in $\mathbb{X}$)}.
\end{align*}
The first sum on the right-hand side is precisely $A_{\z\X'}$.

\item Next, consider the case where both $\v$ and $\v'$ occur in $\mathbb{X}$.  
Write $\mathbb{X}=\mathbb{X}_L(\x)\mathbb{X}_M(\x')\mathbb{X}_R$.  
By Lemma \ref{lem:har-inductive} and Proposition \ref{prop:deriv-explicit},
\begin{align*}
&\vimo_{d_{\Psi_1}(\Psi_2)}(\z, \mathbb{X} * \mathbb{Y},\z') \\
&\begin{aligned}
\stackrel{\eqref{eq:harmonic-composition}}{=}&
\sum_{\mathbb{Y}_L(\ide{y}) \mathbb{Y}_R = \mathbb{Y}}c(\x,\ide{y})\vimo_{d_{\Psi_1}(\Psi_2)}((\z)(\mathbb{X}_L * \mathbb{Y}_L)((\x) - (\ide{y})) (\mathbb{X}_R * \mathbb{Y}_R)(\z')) 
\end{aligned}
\displaybreak[3]\\
&\begin{aligned}
\stackrel{\eqref{eq:harmonic-composition}}{=}&\sum_{\mathbb{Y}_L(\ide{y}) \mathbb{Y}_M(\ide{y}')\mathbb{Y}_R = \mathbb{Y}}c(\x,\ide{y})c(\x',\ide{y}')\\
&\hspace{0.1cm}\times\vimo_{d_{\Psi_1}(\Psi_2)}((\z)(\mathbb{X}_L * \mathbb{Y}_L)((\x) - (\ide{y})) (\mathbb{X}_M * \mathbb{Y}_M)((\x') - (\ide{y}')) (\mathbb{X}_R * \mathbb{Y}_R)(\z'))
\end{aligned}
 \displaybreak[3] \\
&\begin{aligned}
\stackrel{\eqref{eq:deriv-explicit}}{=}& \sum_{\mathbb{Y}_L(\ide{y}) \mathbb{Y}_M(\ide{y}')\mathbb{Y}_R = \mathbb{Y}}c(\x,\ide{y})c(\x',\ide{y}') \\
&\qquad \times \sum_{\substack{\mathbb{V}_L(\v) \mathbb{V}_M (\v') \mathbb{V}_R \\ = (\z)(\mathbb{X}_L * \mathbb{Y}_L)((\x) - (\ide{y})) (\mathbb{X}_M * \mathbb{Y}_M)((\x') - (\ide{y}')) (\mathbb{X}_R * \mathbb{Y}_R)(\z')}}\\
&\hspace{5cm} f_1(\v,\mathbb{V}_M, \v')f_2(\delta_{\z}(\v), \mathbb{V}_L, (\v), (\v'), \mathbb{V}_R,\delta_{\z'}(\v'))
\end{aligned}   \displaybreak[3]\\
&\begin{aligned}
=& \sum_{\mathbb{Y}_L(\ide{y}) \mathbb{Y}_M(\ide{y}')\mathbb{Y}_R = \mathbb{Y}}c(\x,\ide{y})c(\x',\ide{y}') f_1(\x,\mathbb{X}_M * \mathbb{Y}_M, \x')f_2(\z, \mathbb{X}_L * \mathbb{Y}_L , \x, \x', \mathbb{X}_R * \mathbb{Y}_R , \z')\\
&+  (\text{the cases where $(\v,\v') \neq (\x,\x')$}).
\end{aligned}
\end{align*}

Summing over all choices of $\x,\x'$ in the decomposition $\mathbb{X}=\mathbb{X}_L(\x)\mathbb{X}_M(\x')\mathbb{X}_R$ yields that
\begin{equation*}
\begin{aligned}
&\vimo_{d_{\Psi_1}(\Psi_2)}(\z,\mathbb{X}*\mathbb{Y},\z')\\
=&
\sum_{\substack{\mathbb{X}_L(\x)\,\mathbb{X}_M(\x')\,\mathbb{X}_R=\mathbb{X}\\
\mathbb{Y}_L(\ide{y})\,\mathbb{Y}_M(\ide{y}')\,\mathbb{Y}_R=\mathbb{Y}}}
c(\x,\ide{y})\,c(\x',\ide{y}')\,
f_1(\x,\mathbb{X}_M*\mathbb{Y}_M,\x')\,
f_2(\z,\mathbb{X}_L*\mathbb{Y}_L,\x,\x',\mathbb{X}_R*\mathbb{Y}_R,\z')\\
&\qquad + \text{(all remaining terms with $(\v,\v')\neq(\x,\x')$ for every $\x,\x'$ in $\mathbb{X}$)}.
\end{aligned}
\end{equation*}
\end{itemize}
\end{proof}

\begin{lem}\label{lem:vimo-B}
If $\Psi_1 \in \addmr_0\cap \Fad$ and $\Psi_2 \in \ide{h}^{\vee}$, we have
\begin{align*}
&\vimo_{d_{\Psi_1(\Psi_2)}}(\z, \mathbb{X} * \mathbb{Y}, \z')\\
=&B_{\z0} + B_{0\z'}  - B_{00}^z + B_{\x0} + B_{0\x} + B_{\y0} + B_{0\y} - B_{00},
\end{align*}
where, 
\begin{align*}
B_{\z0} =&\sum_{\mathbb{X}_L\mathbb{X}_R = \mathbb{X}} f_1(\z,\mathbb{X}_L,0)f_2(\z,\mathbb{X}_R * \mathbb{Y},\z')
+\sum_{\mathbb{Y}_L\mathbb{Y}_R = \mathbb{Y}} f_1(\z,\mathbb{Y}_L,0)f_2(\z,\mathbb{X} * \mathbb{Y}_R,\z') \displaybreak[3] \\
B_{0\z'} =&\sum_{\mathbb{X}_L\mathbb{X}_R = \mathbb{X}} f_1(0,\mathbb{X}_R,\z')f_2(\z,\mathbb{X}_L * \mathbb{Y},\z')
+\sum_{\mathbb{Y}_L\mathbb{Y}_R = \mathbb{Y}} f_1(0,\mathbb{Y}_R,\z')f_2(\z,\mathbb{X} * \mathbb{Y}_L,\z') \displaybreak[3] \\
B_{00}^{\z} =&\sum_{\mathbb{X}_L\mathbb{X}_R = \mathbb{X}} f_1(0,\mathbb{X}_L,0)f_2(\z,\mathbb{X}_R * \mathbb{Y},\z')
+\sum_{\mathbb{Y}_L\mathbb{Y}_R = \mathbb{Y}} f_1(0,\mathbb{Y}_L,0)f_2(\z,\mathbb{X} * \mathbb{Y}_R,\z') \displaybreak[3] \\
B_{\X0} =& \sum_{\substack{\mathbb{X}_L(\x)\mathbb{X}_M\mathbb{X}_R = \mathbb{X} \\ \mathbb{Y}_L(\ide{y}) \mathbb{Y}_R = \mathbb{Y}}} c(\x , \ide{y}) f_1(\x,\mathbb{X}_M, 0)f_2(\z,\mathbb{X}_L * \mathbb{Y}_L , \x, \mathbb{X}_R * \mathbb{Y}_R, \z') \\
&\hspace{0.2cm} + \sum_{\substack{\mathbb{X}_L(\x)\mathbb{X}_R = \mathbb{X} \\ \mathbb{Y}_L(\ide{y}) \mathbb{Y}_M \mathbb{Y}_R = \mathbb{Y}}} c(\x , \ide{y}) f_1(\x,\mathbb{Y}_M, 0)f_2(\z,\mathbb{X}_L * \mathbb{Y}_L , \x, \mathbb{X}_R * \mathbb{Y}_R, \z') \displaybreak[3] \\ 
B_{0\X'} =& \sum_{\substack{\mathbb{X}_L\mathbb{X}_M(\x)\mathbb{X}_R = \mathbb{X} \\ \mathbb{Y}_L(\ide{y}) \mathbb{Y}_R = \mathbb{Y}}} c(\x , \ide{y}) f_1(0,\mathbb{X}_M, \x)f_2(\z,\mathbb{X}_L * \mathbb{Y}_L , \x, \mathbb{X}_R * \mathbb{Y}_R, \z')  \\
&\hspace{0.2cm} + \sum_{\substack{\mathbb{X}_L(\x)\mathbb{X}_R = \mathbb{X} \\ \mathbb{Y}_L \mathbb{Y}_M (\ide{y})\mathbb{Y}_R = \mathbb{Y}}} c(\x , \ide{y}) f_1(\x,\mathbb{Y}_M, 0)f_2(\z,\mathbb{X}_L * \mathbb{Y}_L , \x, \mathbb{X}_R * \mathbb{Y}_R, \z') \displaybreak[3] \\
B_{\Y0} =&  \sum_{\substack{\mathbb{X}_L(\ide{x})\mathbb{X}_M\mathbb{X}_R = \mathbb{X} \\ \mathbb{Y}_L(\y) \mathbb{Y}_R = \mathbb{Y}}} c(\y , \ide{x}) f_1(\y,\mathbb{X}_M, 0)f_2(\z,\mathbb{X}_L * \mathbb{Y}_L , \y, \mathbb{X}_R * \mathbb{Y}_R, \z') \\
&\hspace{0.2cm} + \sum_{\substack{\mathbb{X}_L(\ide{x})\mathbb{X}_R = \mathbb{X} \\ \mathbb{Y}_L(\y) \mathbb{Y}_M \mathbb{Y}_R = \mathbb{Y}}} c(\y,\ide{x}) f_1(\y,\mathbb{Y}_M, 0)f_2(\z,\mathbb{X}_L * \mathbb{Y}_L , \y, \mathbb{X}_R * \mathbb{Y}_R, \z') \displaybreak[3] \\
B_{0\Y'} =& \sum_{\substack{\mathbb{X}_L\mathbb{X}_M(\ide{x})\mathbb{X}_R = \mathbb{X} \\ \mathbb{Y}_L(\y) \mathbb{Y}_R = \mathbb{Y}}} c(\y , \ide{x}) f_1(0,\mathbb{X}_M, \y) f_2(\z,\mathbb{X}_L * \mathbb{Y}_L , \y, \mathbb{X}_R * \mathbb{Y}_R, \z') \\
&\hspace{0.2cm} + \sum_{\substack{\mathbb{X}_L(\ide{x})\mathbb{X}_R = \mathbb{X} \\ \mathbb{Y}_L \mathbb{Y}_M (\y)\mathbb{Y}_R = \mathbb{Y}}} c(\y , \ide{x}) f_1(\x,\mathbb{Y}_M, 0)f_2(\z,\mathbb{X}_L * \mathbb{Y}_L , \y, \mathbb{X}_R * \mathbb{Y}_R, \z') \displaybreak[3] \\
B_{00} =&\sum_{\substack{\mathbb{X}_L(\x)\mathbb{X}_M\mathbb{X}_R = \mathbb{X} \\ \mathbb{Y}_L(\ide{y}) \mathbb{Y}_R = \mathbb{Y}}} c(\x , \ide{y}) f_1(0,\mathbb{X}_M, 0)f_2(\z,\mathbb{X}_L * \mathbb{Y}_L , \x, \mathbb{X}_R * \mathbb{Y}_R, \z')  \\
&\hspace{0.2cm} + \sum_{\substack{\mathbb{X}_L(\x)\mathbb{X}_R = \mathbb{X} \\ \mathbb{Y}_L(\ide{y}) \mathbb{Y}_M \mathbb{Y}_R = \mathbb{Y}}} c(\x , \ide{y}) f_1(0,\mathbb{Y}_M, 0)f_2(\z,\mathbb{X}_L * \mathbb{Y}_L , \x, \mathbb{X}_R * \mathbb{Y}_R, \z') \\ 
&+\sum_{\substack{\mathbb{X}_L(\ide{x})\mathbb{X}_M\mathbb{X}_R = \mathbb{X} \\ \mathbb{Y}_L(\y) \mathbb{Y}_R = \mathbb{Y}}} c(\y , \ide{x}) f_1(0,\mathbb{X}_M, 0)f_2(\z,\mathbb{X}_L * \mathbb{Y}_L , \y, \mathbb{X}_R * \mathbb{Y}_R, \z') \\
&\hspace{0.2cm} + \sum_{\substack{\mathbb{X}_L(\ide{x})\mathbb{X}_R = \mathbb{X} \\ \mathbb{Y}_L(\y) \mathbb{Y}_M \mathbb{Y}_R = \mathbb{Y}}} c(\y,\ide{x}) f_1(0,\mathbb{Y}_M, 0)f_2(\z,\mathbb{X}_L * \mathbb{Y}_L , \y, \mathbb{X}_R * \mathbb{Y}_R, \z').
\end{align*}

\end{lem}

\begin{proof}
Since $\Psi_1 \in \addmr_0$, $A_{\z\z'} =  f_1(\z, \mathbb{X} * \mathbb{Y},\z')f_2(\z,\z') \stackrel{\eqref{eq:adjoint-harmonic}}{=} 0$.
Let us consider $A_{\z\X'}$.
By Proposition \ref{prop:adjoint-mould}, we have

\begin{align*}
&A_{\z\X'} \\
&\begin{aligned}
=& \sum_{\substack{\mathbb{X}_L (\x) \mathbb{X}_R = \mathbb{X} \\  \mathbb{Y}_L\ide{y} \mathbb{Y}_R = \mathbb{Y}}}
c(\x,\ide{y}) f_1(\z,\mathbb{X}_L * \mathbb{Y}_L, \x)f_2(\z,\x,\mathbb{X}_R * \mathbb{Y}_R,\z')
\end{aligned}\\
&\begin{aligned}
\stackrel{\eqref{eq:adjoint-mould}}{=}& \sum_{\substack{\mathbb{X}_L (\x) \mathbb{X}_R = \mathbb{X} \\  \mathbb{Y}_L(\ide{y}) \mathbb{Y}_R = \mathbb{Y}}}
c(\x,\ide{y})f_2(\z,\x,\mathbb{X}_R * \mathbb{Y}_R,\z')\{ f_1(\z,\mathbb{X}_L * \mathbb{Y}_L, 0) +  f_1(0,\mathbb{X}_L * \mathbb{Y}_L, \x) -  f_1(0,\mathbb{X}_L * \mathbb{Y}_L, 0)  \}.
\end{aligned}
\end{align*}
By $A_{\z\X'\rightarrow \z0}$ (respectively, $A_{\z\X'\rightarrow 0\X'}$, $A_{\z\X'\rightarrow 00}$ ), we denote the first (respectively, second, $(-1) \times $ third) term of the right-hand side above. This means transforming either end component (or both ends) of $f_1$ to $0$.

In a similar manner,  the same convention applies to the cases $(\V,\V') \in S$, i.e., put
\[
A_{\V\V'} := A_{\V\V'\to \V 0}+A_{\V\V'\to 0\V'}-A_{\V\V'\to 00}.
\]

Let us prove the following equality:
\begin{align*}
B_{\V0} = \sum_{\V' \in \{\z',\X',\Y'\}}A_{\V\V'\rightarrow \V0}
\end{align*}
 for $\V \in \{\z, \X,\Y\}$,
\begin{align*}
B_{0\V'}=\sum_{\V \in \{\z,\X,\Y\}}A_{\V\V'\rightarrow 0\V'}
\end{align*}
 for $\V' \in \{\z', \X',\Y'\}$,
\begin{align*}
B_{00}^{\z} =  \sum_{\V' \in \{\z',\X',\Y'\}}A_{\z\V'\rightarrow 00},
\end{align*}
and
\begin{align*}
B_{00} = \sum_{\V' \in \{\z',\X',\Y'\}}A_{\X\V'\rightarrow 00} +  \sum_{\V' \in \{\z',\X',\Y'\}}A_{\Y\V'\rightarrow 00}.
\end{align*}
Essentially, it is enough to show
\begin{align}
B_{\z0}=&A_{\z\X' \rightarrow \z0} + A_{\z\Y' \rightarrow \z0}  \label{eq:Bz0} \\
B_{\X0}=&A_{\X\z' \rightarrow \X0} + A_{\X\X' \rightarrow \X0} +  A_{\X\Y' \rightarrow \X0} \label{eq:Bx0} 
\end{align}
because of the symmetry. Namely,
\[
\xymatrix{
B_{\X0} \ar[r]^{\mathbb{X}\leftrightarrow \mathbb{Y}} & B_{\Y0} \ar[l]
},
\xymatrix{
B_{\V0} \ar[r]^{\begin{matrix} \text{counterwise} \\ \text{discussion} \end{matrix}} & B_{0\V'} \ar[l]
}
,\text{ and}
\xymatrix{
B_{\z0} \ar[r]^{\begin{matrix} \text{counterwise} \\ \text{discussion} \end{matrix}} & B_{0\z'}  \ar[l]
}
\]
and
\[
b_0(B_{\x0} + B_{\y0}) =  B_{00}, b_0'(B_{\z0}) = B_{00}^z
\]
($b_0 (A_{? s\rightarrow ? 0}) = A_{? s \rightarrow 00}$, $b_0' (A_{s? \rightarrow 0?}) = A_{s ? \rightarrow 00}$).

First, let us prove \eqref{eq:Bz0}. 
By direct calculation, we have
\begin{align}
 &\sum_{\V' \in \{\z',\X',\Y'\}}A_{\z\V'\rightarrow \z0} \displaybreak[3] \notag\\
=&\sum_{\substack{\mathbb{X}_L (\x) \mathbb{X}_R = \mathbb{X} \\  \mathbb{Y}_L\ide{y} \mathbb{Y}_R = \mathbb{Y}}}
c(\x,\ide{y}) f_1(\z,\mathbb{X}_L * \mathbb{Y}_L, 0)f_2(\z,\x,\mathbb{X}_R * \mathbb{Y}_R,\z') \notag \\
&+\sum_{\substack{\mathbb{X}_L (\ide{x}) \mathbb{X}_R = \mathbb{X} \\  \mathbb{Y}_L (\y) \mathbb{Y}_R = \mathbb{Y}}}
c(\y, \ide{x}) f_1(\z,\mathbb{X}_L * \mathbb{Y}_L, 0)f_2(\z,\y,\mathbb{X}_R * \mathbb{Y}_R,\z') \notag \displaybreak[3] \\
=& \sum_{\substack{\mathbb{X}_L \mathbb{X}_R = \mathbb{X} \\ \mathbb{Y}_L \mathbb{Y}_R = \mathbb{Y}}}
f_1(\z,\mathbb{X}_L * \mathbb{Y}_L, 0) \notag \\
&\hspace{-0.2cm}\times\Bigg\{ \sum_{\substack{(\x)\mathbb{X}_{RR} = \mathbb{X}_R \\ (\ide{y})\mathbb{Y}_{RR} = \mathbb{Y}_R}} c(\x,\ide{y})f_2(\z,\x,\mathbb{X}_{RR} * \mathbb{Y}_{RR},\z')
+ \sum_{\substack{(\ide{x})\mathbb{X}_{RR} = \mathbb{X}_R \\ (\y)\mathbb{Y}_{RR} = \mathbb{Y}_R}} 
c(\y, \ide{x}) f_2(\z,\y,\mathbb{X}_{RR} * \mathbb{Y}_{RR},\z') \Bigg\} \notag \displaybreak[3] \\
=&\sum_{\substack{\mathbb{X}_L \mathbb{X}_R = \mathbb{X} \\ \mathbb{Y}_L \mathbb{Y}_R = \mathbb{Y}}}
f_1(\z,\mathbb{X}_L * \mathbb{Y}_L, 0) f_2(\z , \mathbb{X}_R * \mathbb{Y}_R , \z').  \displaybreak[3] \label{eq:Bz0-1}
\end{align}
In the last equality, we use the definition of the harmonic product.
Thus, due to $\Psi_1 \in \addmr_0$, we conclude

\begin{align*}
&\sum_{\V' \in \{\z',\X',\Y'\}}A_{\V\V'\rightarrow \z0} \displaybreak[3]\\
\stackrel{\eqref{eq:Bz0-1}}{=}& \sum_{\substack{\mathbb{X}_L \mathbb{X}_R = \mathbb{X} \\ \mathbb{Y}_L \mathbb{Y}_R = \mathbb{Y}}}f_1(\z,\mathbb{X}_L * \mathbb{Y}_L,0)f_2(\z,  \mathbb{X}_{R} * \mathbb{Y}_{R}, \z') \displaybreak[3]\\
\stackrel{\eqref{eq:adjoint-harmonic}}{=}&
\sum_{\substack{\mathbb{X}_L \mathbb{X}_R = \mathbb{X} }}f_1(\z,\mathbb{X}_L ,0)f_2(\z,  \mathbb{X}_{R} * \mathbb{Y}, \z')+\sum_{\substack{\mathbb{Y}_L \mathbb{Y}_R = \mathbb{Y} }}f_1(\z,\mathbb{Y}_L ,0)f_2(\z,  \mathbb{X} * \mathbb{Y}_R, \z') \displaybreak[3] \\
=&B_{\z0}.
\end{align*}

Next, let us prove \eqref{eq:Bx0}. Direct calculations show
\begin{align}
& A_{\X\z' \rightarrow \X0} + A_{\X\X' \rightarrow \X0} +  A_{\X\Y' \rightarrow \X0} \notag \\
=&\sum_{\substack{\mathbb{X}_L(x)\mathbb{X}_R = \mathbb{X} \\ \mathbb{Y}_L(\ide{y})\mathbb{Y}_R = \mathbb{Y}}}c(\x,\ide{y}) f_1(\x,\mathbb{X}_R * \mathbb{Y}_R, 0)f_2(\z,\mathbb{X}_L * \mathbb{Y}_L, \x ,\z') \notag \\
&+ \sum_{\substack{\mathbb{X}_L (\x) \mathbb{X}_M(\x')\mathbb{X}_R = \mathbb{X} \\ \mathbb{Y}_L (\ide{y}) \mathbb{Y}_M(\ide{y}')\mathbb{Y}_R = \mathbb{Y} }} c(\x,\ide{y})c(\x',\ide{y}')f_1(\x,\mathbb{X}_M * \mathbb{Y}_M, 0)f_2(\z, \mathbb{X}_L * \mathbb{Y}_L , \x, \x', \mathbb{X}_R * \mathbb{Y}_R , \z') \notag \\
&+\sum_{\substack{\mathbb{X}_L (\x) \mathbb{X}_M(\ide{x}')\mathbb{X}_R = \mathbb{X} \\ \mathbb{Y}_L (\ide{y}) \mathbb{Y}_M(\y')\mathbb{Y}_R = \mathbb{Y} }} c(\x,\ide{y})c(\y',\ide{x}')f_1(\x,\mathbb{X}_M * \mathbb{Y}_M, 0)f_2(\z, \mathbb{X}_L * \mathbb{Y}_L , \x, \y', \mathbb{X}_R * \mathbb{Y}_R , \z') \displaybreak[3]  \notag \\
=&\sum_{\substack{\mathbb{X}_L(x)\mathbb{X}_R = \mathbb{X} \\ \mathbb{Y}_L(\ide{y})\mathbb{Y}_R = \mathbb{Y}}}c(\x,\ide{y}) f_1(\x,\mathbb{X}_R * \mathbb{Y}_R, 0)f_2(\z,\mathbb{X}_L * \mathbb{Y}_L, \x ,\z') \notag \\
&+\sum_{\substack{ \mathbb{X}_L (\x) \mathbb{X}_M \mathbb{X}_R = \mathbb{X} \\\mathbb{Y}_L (\ide{y}) \mathbb{Y}_M \mathbb{Y}_R = \mathbb{Y}  }} c(\x,\ide{y}) f_1(\x, \mathbb{X}_M  * \mathbb{Y}_M ,0) \notag \\
& 
\begin{aligned}
\times &\Bigg\{  \sum_{\substack{(\x')\mathbb{X}_{RR} = \mathbb{X}_R \\ (\ide{y}')\mathbb{Y}_{RR} = \mathbb{Y}_R }} c(\x',\ide{y}')f_2(\z, \mathbb{X}_L * \mathbb{Y}_L , \x, \x', \mathbb{X}_{RR} * \mathbb{Y}_{RR} , \z') + \sum_{\substack{(\ide{x}')\mathbb{X}_{RR} = \mathbb{X}_R \\ (\y')\mathbb{Y}_{RR} = \mathbb{Y}_R }}  c(\y',\ide{x}')f_2(\z, \mathbb{X}_L * \mathbb{Y}_L , \x, \y', \mathbb{X}_{RR} * \mathbb{Y}_{RR} , \z')  \Bigg\} 
\end{aligned}\notag 
\displaybreak[3] \\
&\begin{aligned}
\hspace{-0.3cm}=&\sum_{\substack{\mathbb{X}_L(x)\mathbb{X}_R = \mathbb{X} \\ \mathbb{Y}_L(\ide{y})\mathbb{Y}_R = \mathbb{Y}}}c(\x,\ide{y}) f_1(\x,\mathbb{X}_R * \mathbb{Y}_R, 0)f_2(\z,\mathbb{X}_L * \mathbb{Y}_L, \x ,\z') +\sum_{\substack{ \mathbb{X}_L (\x) \mathbb{X}_M \mathbb{X}_R = \mathbb{X} \\\mathbb{Y}_L (\ide{y}) \mathbb{Y}_M \mathbb{Y}_R = \mathbb{Y} \\ (\mathbb{X}_R, \mathbb{Y}_R) \neq (\emptyset, \emptyset) }} c(\x,\ide{y}) f_1(\x, \mathbb{X}_M,0)f_2(\z, \mathbb{X}_L * \mathbb{Y}_L , \x, \mathbb{X}_{R} * \mathbb{Y}_{R} , \z')
\end{aligned}
 \notag \\
=&\sum_{\substack{ \mathbb{X}_L (\x) \mathbb{X}_M \mathbb{X}_R = \mathbb{X} \\\mathbb{Y}_L (\ide{y}) \mathbb{Y}_M \mathbb{Y}_R = \mathbb{Y}}} c(\x,\ide{y}) f_1(\x, \mathbb{X}_M * \mathbb{Y}_M,0)f_2(\z, \mathbb{X}_L * \mathbb{Y}_L , \x, \mathbb{X}_{R} * \mathbb{Y}_{R} , \z').  \label{eq:BV0-1}
\end{align}

In the above third equality, we use the definition of the harmonic product.
Therefore, due to $\Psi_1 \in \addmr_0$, we have
\begin{align*}
&A_{\X\z'\rightarrow \X0} + A_{\X\X'\rightarrow \X0} + A_{\X\Y'\rightarrow \X0} \\
\stackrel{\eqref{eq:BV0-1}}{=}& \sum_{\substack{ \mathbb{X}_L (\x) \mathbb{X}_M \mathbb{X}_R = \mathbb{X} \\\mathbb{Y}_L (\ide{y}) \mathbb{Y}_M \mathbb{Y}_R = \mathbb{Y}}} c(\x,\ide{y}) f_1(\x, \mathbb{X}_M  * \mathbb{Y}_M ,0)f_2(\z, \mathbb{X}_L * \mathbb{Y}_L , \x, \mathbb{X}_{R} * \mathbb{Y}_{R} , \z') \displaybreak[3] \\
\stackrel{\eqref{eq:adjoint-harmonic}}{=}&
 \sum_{\substack{ \mathbb{X}_L (\x) \mathbb{X}_M \mathbb{X}_R = \mathbb{X} \\\mathbb{Y}_L (\ide{y}) \mathbb{Y}_R = \mathbb{Y}}} c(\x,\ide{y}) f_1(\x, \mathbb{X}_M  ,0)f_2(\z, \mathbb{X}_L * \mathbb{Y}_L , \x, \mathbb{X}_{R} * \mathbb{Y}_{R} , \z') \\
&+ \sum_{\substack{ \mathbb{X}_L (\x) \mathbb{X}_R = \mathbb{X} \\\mathbb{Y}_L (\ide{y}) \mathbb{Y}_M \mathbb{Y}_R = \mathbb{Y}}} c(\x,\ide{y}) f_1(\x,  \mathbb{Y}_M ,0)f_2(\z, \mathbb{X}_L * \mathbb{Y}_L , \x, \mathbb{X}_{R} * \mathbb{Y}_{R} , \z') 
\displaybreak[3]\\
=&B_{\X0}
\end{align*}
as claimed.

\end{proof}

\begin{lem}\label{eq:Conv-harmonic-z}
For $\Psi_1$ $\in \addmr_0 \cap \Fad $ and $\Psi_2\in \addmr_0$,
\begin{align*}
B_{\z0} + B_{0\z'} - B_{00}^{\z}=&f_1(\z,\mathbb{X},\z')f_2(\z,\mathbb{Y},\z') +  f_1(\z,\mathbb{Y},\z')f_2(\z,\mathbb{X},\z').
\end{align*}
\end{lem}

\begin{proof}
Since $\Psi_2 \in \addmr_0$, we have
\begin{align*}
B_{\z0} =& f_1(0,\mathbb{X},\z')f_2(\z,\mathbb{Y},\z') +  f_1(0,\mathbb{Y},\z')f_2(\z,\mathbb{X},\z')\\ 
B_{0\z'} =& f_1(\z,\mathbb{X},0)f_2(\z,\mathbb{Y},\z') +  f_1(\z,\mathbb{Y},0)f_2(\z,\mathbb{X},\z')\\ 
B_{00}^{\z} =& f_1(0,\mathbb{X},0)f_2(\z,\mathbb{Y},\z') +  f_1(0,\mathbb{Y},0)f_2(\z,\mathbb{X},\z').
\end{align*}
Since $\Psi_1 \in \Fad$, we have
\begin{align*}
&B_{\z0} + B_{0\z'} - B_{00}^z \\
=&\{ f_1(0,\mathbb{X},\z') + f_1(\z,\mathbb{X},0) - f_1(0,\mathbb{X},0) \}f_2(\z,\mathbb{Y},\z') \\
&\qquad \qquad+ \{ f_1(0,\mathbb{Y},\z') + f_1(\z,\mathbb{Y},0) - f_1(0,\mathbb{Y},0) \}f_2(\z,\mathbb{X},\z')\\
\stackrel{\eqref{eq:adjoint-mould}}{=} &f_1(\z,\mathbb{X},\z')f_2(\z,\mathbb{Y},\z') +  f_1(\z,\mathbb{Y},\z')f_2(\z,\mathbb{X},\z')
\end{align*}
as claimed.
\end{proof}

For $(\V,\V') \in \{(\X,0),(\Y,0),(0,\X'),(0,\Y'),(0,0)\}$, decompose $B_{\V\V'}$ as
\[
B_{\V\V'} = B_{\V\V'}^{\X} + B_{\V\V'}^{\Y},
\]
where $B_{\V\V'}^{\X}$ (respectively, $B_{\V\V'}^{\Y}$) is the sum of those terms of $B_{\V\V'}$ for which the middle part of $f_{1}$ (i.e., $f_{1}$ with its first and last letters removed) is a word in $\mathbb{X}$ (respectively, over $\mathbb{Y}$).

\begin{lem}\label{lem:Conv-harmonic}
For $\Psi_1$ $\in \addmr_0 \cap \Fad$ and $\Psi_2\in \addmr_0$, we have\begin{align}
&B_{\X0}^{\X} =  \sum_{\substack{\mathbb{X}_L(\x)\mathbb{X}_M\mathbb{X}_R = \mathbb{X} \\ \mathbb{Y}_L(\y) \mathbb{Y}_R = \mathbb{Y}}} c(\x,\y ) f_1(\x,\mathbb{X}_M, 0)f_2(\z,\mathbb{X}_L * \mathbb{Y}_L , \y, \mathbb{X}_R * \mathbb{Y}_R, \z'), \label{eq:Conv-harmonic-1} \\
&B_{0\X'}^{\X} =  \sum_{\substack{\mathbb{X}_L\mathbb{X}_M(\x)\mathbb{X}_R = \mathbb{X} \\ \mathbb{Y}_L(\y) \mathbb{Y}_R = \mathbb{Y}}} c(\x,\y ) f_1(0,\mathbb{X}_M, \x)f_2(\z,\mathbb{X}_L * \mathbb{Y}_L , \y, \mathbb{X}_R * \mathbb{Y}_R, \z'), \text{ and } \label{eq:Conv-harmonic-2}\\
&\begin{aligned}
 &\sum_{\substack{\mathbb{X}_L(\x)\mathbb{X}_M\mathbb{X}_R = \mathbb{X} \\ \mathbb{Y}_L(\ide{y}) \mathbb{Y}_R = \mathbb{Y}}} c(\x , \ide{y}) f_1(0,\mathbb{X}_M, 0)f_2(\z,\mathbb{X}_L * \mathbb{Y}_L , \x, \mathbb{X}_R * \mathbb{Y}_R, \z')\\
&\hspace{3cm} =  \sum_{\substack{\mathbb{X}_L(\x)\mathbb{X}_M\mathbb{X}_R = \mathbb{X} \\ \mathbb{Y}_L(\y) \mathbb{Y}_R = \mathbb{Y}}} c(\x,\y ) f_1(0,\mathbb{X}_M, 0)f_2(\z,\mathbb{X}_L * \mathbb{Y}_L , \y, \mathbb{X}_R * \mathbb{Y}_R, \z').
\end{aligned}
 \label{eq:Conv-harmonic-3} 
\end{align}

\end{lem}

\begin{proof}

We prove only \eqref{eq:Conv-harmonic-1}. The others can be shown in a similar manner.

Consider the following decomposition: $\mathbb{X}_L (\x)\mathbb{X}_M\mathbb{X}_R = \mathbb{X}$.
Then, applying Lemma \ref{lem:har-inductive} to $(\mathbb{X}_L(\x)\mathbb{X}_R) * \mathbb{Y}$, we have
\begin{align}
(\mathbb{X}_L(\x)\mathbb{X}_R) * \mathbb{Y}
\stackrel{\eqref{eq:harmonic-composition}}{=}\sum_{\mathbb{Y}_L(\ide{y})\mathbb{Y}_R = \mathbb{Y}} c (\x,\ide{y})(\mathbb{X}_L * \mathbb{Y}_L)((\x) - (\ide{y}))\mathbb{X}_R * \mathbb{Y}_R. \label{eq:Conv-harmonic-10}
\end{align}
Since $\Psi_2 \in \addmr_0$, it follows that
\begin{align}
0=& f_2(\z , (\mathbb{X}_L\cdot (\x)\cdot \mathbb{X}_R) * \mathbb{Y},\z') \notag \\
\stackrel{\eqref{eq:Conv-harmonic-10}}{=}&\sum_{\mathbb{Y}_L(\ide{y})\mathbb{Y}_R = \mathbb{Y}} c(\x,\ide{y}) f_2(\z,( \mathbb{X}_L * \mathbb{Y}_L), ((\x) - (\ide{y})), (\mathbb{X}_R * \mathbb{Y}_R),\z') \notag \\
\Longleftrightarrow& \notag \\
&
\begin{aligned}
&\sum_{\mathbb{Y}_L(\ide{y})\mathbb{Y}_R = \mathbb{Y}} c(\x,\ide{y}) f_2(\z,( \mathbb{X}_L * \mathbb{Y}_L), (\x), (\mathbb{X}_R * \mathbb{Y}_R),\z')= \sum_{\mathbb{Y}_L(\y)\mathbb{Y}_R = \mathbb{Y}} c(\x,\y) f_2(\z,( \mathbb{X}_L * \mathbb{Y}_L), (\y), (\mathbb{X}_R * \mathbb{Y}_R),\z').
\end{aligned}
 \label{eq:Conv-harmonic-11}
\end{align}
Thus we have
\begin{align*}
&B_{\X0}^{\X}\\
=&\sum_{\substack{\mathbb{X}_L(\x)\mathbb{X}_M \mathbb{X}_R \\ \mathbb{Y}_L(\ide{y})\mathbb{Y}_R = \mathbb{Y}}} c(\x,\ide{y}) f_1(\x,\mathbb{X}_M,0)f_2(\z,(\mathbb{X}_L * \mathbb{Y}_L),\x,(\mathbb{X}_R * \mathbb{Y}_R),\z')\\
\stackrel{\eqref{eq:Conv-harmonic-11}}{=}&
\sum_{\substack{\mathbb{X}_L(\x)\mathbb{X}_M \mathbb{X}_R \\ \mathbb{Y}_L(\y)\mathbb{Y}_R = \mathbb{Y}}} c(\x,\y) f_1(\x,\mathbb{X}_M,0)f_2(\z,(\mathbb{X}_L * \mathbb{Y}_L),\y,(\mathbb{X}_R * \mathbb{Y}_R),\z')
\end{align*}
as claimed.
\end{proof}

\begin{lem}\label{lem:Conv-harmonic-00}
For $\Psi_1$ $\in \addmr_0 \cap \Fad$ and $\Psi_2\in \addmr_0$, we have
\begin{align}
B_{00}^{\X}
=\sum_{\substack{\mathbb{X}_L\mathbb{X}_M\mathbb{X}_R = \mathbb{X} \\ \mathbb{Y}_L(\y) \mathbb{Y}_R = \mathbb{Y}}}  f_1(0,\mathbb{X}_M, 0)f_2(\z,\mathbb{X}_L * \mathbb{Y}_L , \y, \mathbb{X}_R * \mathbb{Y}_R, \z'). \label{eq:Conv-harmonic-00}
\end{align}
\end{lem}

\begin{proof}

By Lemma \ref{lem:Conv-harmonic}, we have

\begin{align*}
&B_{00}^{\X}\\
&\begin{aligned}
=&\sum_{\substack{\mathbb{X}_L(\x)\mathbb{X}_M\mathbb{X}_R = \mathbb{X} \\ \mathbb{Y}_L(\ide{y}) \mathbb{Y}_R = \mathbb{Y}}} c(\x , \ide{y}) f_1(0,\mathbb{X}_M, 0)f_2(\z,\mathbb{X}_L * \mathbb{Y}_L , \x, \mathbb{X}_R * \mathbb{Y}_R, \z') \\
& \qquad + \sum_{\substack{\mathbb{X}_L(\ide{x})\mathbb{X}_M\mathbb{X}_R = \mathbb{X} \\ \mathbb{Y}_L(\y) \mathbb{Y}_R = \mathbb{Y}}} c(\y , \ide{x}) f_1(0,\mathbb{X}_M, 0)f_2(\z,\mathbb{X}_L * \mathbb{Y}_L , \y, \mathbb{X}_R * \mathbb{Y}_R, \z')
\end{aligned}
\\ 
&\begin{aligned}
\stackrel{\eqref{eq:Conv-harmonic-3}}{=}& \sum_{\substack{\mathbb{X}_L(\x)\mathbb{X}_M\mathbb{X}_R = \mathbb{X} \\ \mathbb{Y}_L(\y) \mathbb{Y}_R = \mathbb{Y}}} c(\x,\y ) f_1(0,\mathbb{X}_M, 0)f_2(\z,\mathbb{X}_L * \mathbb{Y}_L , \y, \mathbb{X}_R * \mathbb{Y}_R, \z')\\
& \hspace{2cm} +\sum_{\substack{\mathbb{X}_L(\ide{x})\mathbb{X}_M\mathbb{X}_R = \mathbb{X} \\ \mathbb{Y}_L(\y) \mathbb{Y}_R = \mathbb{Y}}} c(\y , \ide{x}) f_1(0,\mathbb{X}_M, 0)f_2(\z,\mathbb{X}_L * \mathbb{Y}_L , \y, \mathbb{X}_R * \mathbb{Y}_R, \z')
\end{aligned}
 \displaybreak[3] \\
&\begin{aligned}
=& \sum_{\substack{\mathbb{X}_L(\x)\mathbb{X}_M\mathbb{X}_R = \mathbb{X} \\ \mathbb{Y}_L(\y) \mathbb{Y}_R = \mathbb{Y}}} c(\x,\y ) f_1(0,\mathbb{X}_M, 0)f_2(\z,\mathbb{X}_L * \mathbb{Y}_L , \y, \mathbb{X}_R * \mathbb{Y}_R, \z') \\
&\hspace{2cm} +\sum_{\substack{\mathbb{X}_L(\x)\mathbb{X}_M\mathbb{X}_R = \mathbb{X} \\ \mathbb{Y}_L(\y) \mathbb{Y}_R = \mathbb{Y}}} c(\y , \x) f_1(0,\mathbb{X}_M, 0)f_2(\z,\mathbb{X}_L * \mathbb{Y}_L , \y, \mathbb{X}_R * \mathbb{Y}_R, \z') \\
&\hspace{4cm} +\sum_{\substack{\mathbb{X}_L\mathbb{X}_M\mathbb{X}_R = \mathbb{X} \\ \mathbb{Y}_L(\y) \mathbb{Y}_R = \mathbb{Y}}}  f_1(0,\mathbb{X}_M, 0)f_2(\z,\mathbb{X}_L * \mathbb{Y}_L , \y, \mathbb{X}_R * \mathbb{Y}_R, \z') 
\end{aligned}
\displaybreak[3] \\
&\begin{aligned}
=&\sum_{\substack{\mathbb{X}_L\mathbb{X}_M\mathbb{X}_R = \mathbb{X} \\ \mathbb{Y}_L(\y) \mathbb{Y}_R = \mathbb{Y}}}  f_1(0,\mathbb{X}_M, 0)f_2(\z,\mathbb{X}_L * \mathbb{Y}_L , \y, \mathbb{X}_R * \mathbb{Y}_R, \z')
\end{aligned}
\end{align*}
In the above third equality, we distinguish the case where $\ide{x}$ lies in $\mathbf{X}$ or not.
Thus, we obtain our desired equality.
\end{proof}

\begin{lem}\label{lem:Conv-harmonic-3}
For $\Psi_1$ $\in \addmr_0 \cap \Fad \cap \V_{\strprty}$ and $\Psi_2\in \addmr_0$, we have
\begin{align}
&B_{\Y0}^{\X} + B_{0\Y'}^{\X} - B_{00}^{\X} \notag  \\
=&\sum_{\substack{\mathbb{X}_L\mathbb{X}_M \mathbb{X}_R = \mathbb{X} \\ \mathbb{Y}_L(\y)\mathbb{Y}_R = \mathbb{Y}}} f_1(\y,\mathbb{X}_M,\y)f_2(\z,(\mathbb{X}_L * \mathbb{Y}_L),\y,(\mathbb{X}_R * \mathbb{Y}_R),\z') \notag \\
&\begin{aligned}
&\hspace{1cm}+\sum_{\substack{\mathbb{X}_L(\x)\mathbb{X}_M \mathbb{X}_R = \mathbb{X} \\ \mathbb{Y}_L(\y)\mathbb{Y}_R = \mathbb{Y}}} c(\y,\x) f_1(\y,\mathbb{X}_M,0)f_2(\z,(\mathbb{X}_L * \mathbb{Y}_L),\y,(\mathbb{X}_R * \mathbb{Y}_R),\z')\\
&\hspace{3cm}+\sum_{\substack{\mathbb{X}_L\mathbb{X}_M(\x) \mathbb{X}_R \\ \mathbb{Y}_L(\y)\mathbb{Y}_R = \mathbb{Y}}} c(\y,\x) f_1(0,\mathbb{X}_M,\y)f_2(\z,(\mathbb{X}_L * \mathbb{Y}_L),\y,(\mathbb{X}_R * \mathbb{Y}_R),\z'). \label{eq:Conv-harmonic-4}
\end{aligned}
\end{align}
\end{lem}
\begin{proof}
By direct calculation, we have
\begin{align*}
&B_{\Y0}^{\X} + B_{0\Y'}^{\X} - B_{00}^{\X} \displaybreak[3] \\
&\begin{aligned}
=&\sum_{\substack{\mathbb{X}_L(\ide{x})\mathbb{X}_M \mathbb{X}_R = \mathbb{X} \\ \mathbb{Y}_L(\y)\mathbb{Y}_R = \mathbb{Y}}} c(\y,\ide{x}) f_1(\y,\mathbb{X}_M,0)f_2(\z,(\mathbb{X}_L * \mathbb{Y}_L),\y,(\mathbb{X}_R * \mathbb{Y}_R),\z') \\
&\hspace{1cm}+\sum_{\substack{\mathbb{X}_L\mathbb{X}_M(\ide{x}) \mathbb{X}_R = \mathbb{X} \\ \mathbb{Y}_L(\y)\mathbb{Y}_R = \mathbb{Y}}} c(\y,\ide{x}) f_1(0,\mathbb{X}_M,\y)f_2(\z,(\mathbb{X}_L * \mathbb{Y}_L),\y,(\mathbb{X}_R * \mathbb{Y}_R),\z')-B_{00}^{\X} 
\end{aligned} \displaybreak[3] \\
&\begin{aligned}
\stackrel{\eqref{eq:Conv-harmonic-00}}{=}&\sum_{\substack{\mathbb{X}_L(\ide{x})\mathbb{X}_M \mathbb{X}_R = \mathbb{X} \\ \mathbb{Y}_L(\y)\mathbb{Y}_R = \mathbb{Y}}} c(\y,\ide{x}) f_1(\y,\mathbb{X}_M,0)f_2(\z,(\mathbb{X}_L * \mathbb{Y}_L),\y,(\mathbb{X}_R * \mathbb{Y}_R),\z') \\
&\hspace{1cm} +\sum_{\substack{\mathbb{X}_L\mathbb{X}_M(\ide{x}) \mathbb{X}_R = \mathbb{X} \\ \mathbb{Y}_L(\y)\mathbb{Y}_R = \mathbb{Y}}} c(\y,\ide{x}) f_1(0,\mathbb{X}_M,\y)f_2(\z,(\mathbb{X}_L * \mathbb{Y}_L),\y,(\mathbb{X}_R * \mathbb{Y}_R),\z')\\
&\hspace{3cm} -\sum_{\substack{\mathbb{X}_L\mathbb{X}_M\mathbb{X}_R = \mathbb{X} \\ \mathbb{Y}_L(\y) \mathbb{Y}_R = \mathbb{Y}}}  f_1(0,\mathbb{X}_M, 0)f_2(\z,\mathbb{X}_L * \mathbb{Y}_L , \y, \mathbb{X}_R * \mathbb{Y}_R, \z') 
\end{aligned} \displaybreak[3] \\
&\begin{aligned}
=&\sum_{\substack{\mathbb{X}_L\mathbb{X}_M \mathbb{X}_R = \mathbb{X} \\ \mathbb{Y}_L(\y)\mathbb{Y}_R = \mathbb{Y}}}  f_1(\y,\mathbb{X}_M,0)f_2(\z,(\mathbb{X}_L * \mathbb{Y}_L),\y,(\mathbb{X}_R * \mathbb{Y}_R),\z') \\
&\hspace{1cm}+\sum_{\substack{\mathbb{X}_L(\x)\mathbb{X}_M \mathbb{X}_R = \mathbb{X} \\ \mathbb{Y}_L(\y)\mathbb{Y}_R = \mathbb{Y}}} c(\y,\x) f_1(\y,\mathbb{X}_M,0)f_2(\z,(\mathbb{X}_L * \mathbb{Y}_L),\y,(\mathbb{X}_R * \mathbb{Y}_R),\z')\\
&+\sum_{\substack{\mathbb{X}_L\mathbb{X}_M\mathbb{X}_R = \mathbb{X} \\ \mathbb{Y}_L(\y)\mathbb{Y}_R = \mathbb{Y}}} f_1(0,\mathbb{X}_M,\y)f_2(\z,(\mathbb{X}_L * \mathbb{Y}_L),\y,(\mathbb{X}_R * \mathbb{Y}_R),\z')\\
&\hspace{1cm}+\sum_{\substack{\mathbb{X}_L\mathbb{X}_M(\x) \mathbb{X}_R = \mathbb{X} \\ \mathbb{Y}_L(\y)\mathbb{Y}_R = \mathbb{Y}}} c(\y,\x) f_1(0,\mathbb{X}_M,\y)f_2(\z,(\mathbb{X}_L * \mathbb{Y}_L),\y,(\mathbb{X}_R * \mathbb{Y}_R),\z') \\
&-\sum_{\substack{\mathbb{X}_L\mathbb{X}_M\mathbb{X}_R = \mathbb{X} \\ \mathbb{Y}_L(\y) \mathbb{Y}_R = \mathbb{Y}}}  f_1(0,\mathbb{X}_M, 0)f_2(\z,\mathbb{X}_L * \mathbb{Y}_L , \y, \mathbb{X}_R * \mathbb{Y}_R, \z') 
\end{aligned}
\displaybreak[3] \\
&\begin{aligned}
=&\sum_{\substack{\mathbb{X}_L\mathbb{X}_M \mathbb{X}_R = \mathbb{X} \\ \mathbb{Y}_L(\y)\mathbb{Y}_R = \mathbb{Y}}}\{  f_1(\y,\mathbb{X}_M,0) +  f_1(0,\mathbb{X}_M,\y)  -  f_1(0,\mathbb{X}_M,0)\}\\
&\hspace{6cm}\times f_2(\z,(\mathbb{X}_L * \mathbb{Y}_L),\y,(\mathbb{X}_R * \mathbb{Y}_R),\z')\\
&+\sum_{\substack{\mathbb{X}_L(\x)\mathbb{X}_M \mathbb{X}_R = \mathbb{X} \\ \mathbb{Y}_L(\y)\mathbb{Y}_R = \mathbb{Y}}} c(\y,\x) f_1(\y,\mathbb{X}_M,0)f_2(\z,(\mathbb{X}_L * \mathbb{Y}_L),\y,(\mathbb{X}_R * \mathbb{Y}_R),\z')\\
&\hspace{1cm} +\sum_{\substack{\mathbb{X}_L\mathbb{X}_M(\x) \mathbb{X}_R \\ \mathbb{Y}_L(\y)\mathbb{Y}_R = \mathbb{Y}}} c(\y,\x) f_1(0,\mathbb{X}_M,\y)f_2(\z,(\mathbb{X}_L * \mathbb{Y}_L),\y,(\mathbb{X}_R * \mathbb{Y}_R),\z')
\end{aligned}
\displaybreak[3]\\
&\begin{aligned}
\stackrel{\eqref{eq:adjoint-harmonic}}{=}&\sum_{\substack{\mathbb{X}_L\mathbb{X}_M \mathbb{X}_R = \mathbb{X} \\ \mathbb{Y}_L(\y)\mathbb{Y}_R = \mathbb{Y}}} f_1(\y,\mathbb{X}_M,\y)f_2(\z,(\mathbb{X}_L * \mathbb{Y}_L),\y,(\mathbb{X}_R * \mathbb{Y}_R),\z')
\\
&+\sum_{\substack{\mathbb{X}_L(\x)\mathbb{X}_M \mathbb{X}_R = \mathbb{X} \\ \mathbb{Y}_L(\y)\mathbb{Y}_R = \mathbb{Y}}} c(\y,\x) f_1(\y,\mathbb{X}_M,0)f_2(\z,(\mathbb{X}_L * \mathbb{Y}_L),\y,(\mathbb{X}_R * \mathbb{Y}_R),\z') \\
&\hspace{1cm}+\sum_{\substack{\mathbb{X}_L\mathbb{X}_M(\x) \mathbb{X}_R \\ \mathbb{Y}_L(\y)\mathbb{Y}_R = \mathbb{Y}}} c(\y,\x) f_1(0,\mathbb{X}_M,\y)f_2(\z,(\mathbb{X}_L * \mathbb{Y}_L),\y,(\mathbb{X}_R * \mathbb{Y}_R),\z')
\end{aligned}
 \displaybreak[3]
\end{align*}
as claimed.
\end{proof}

\begin{lem}\label{eq:Conv-harmonic-parity}
For $\Psi_1$ $\in \addmr_0 \cap \Fad \cap \V_{\strprty}$ and $\Psi_2\in \addmr_0$, we have
\begin{align*}
B_{\X0}^{\X} + B_{0\X'}^{\X} + B_{\Y0}^{\X} + B_{0\Y'}^{\X} - B_{00}^{\X} = 0.
\end{align*}
\end{lem}

\begin{proof}
By Lemma \ref{lem:Conv-harmonic} and Lemma \ref{lem:Conv-harmonic-3}, we have
\begin{align*}
&B_{\X0}^{\X} + B_{0\X'}^{\X} + B_{\Y0}^{\X} + B_{0\Y'}^{\X} - B_{00}^{\X}  \displaybreak[3] \\
&\begin{aligned}
\stackrel{\substack{\eqref{eq:Conv-harmonic-1} \\ \eqref{eq:Conv-harmonic-2} }}{=} & 
 \sum_{\substack{\mathbb{X}_L(\x)\mathbb{X}_M\mathbb{X}_R = \mathbb{X} \\ \mathbb{Y}_L(\y) \mathbb{Y}_R = \mathbb{Y}}} c(\x,\y ) f_1(\x,\mathbb{X}_M, 0)f_2(\z,\mathbb{X}_L * \mathbb{Y}_L , \y, \mathbb{X}_R * \mathbb{Y}_R, \z') \\
&\hspace{2cm} + \sum_{\substack{\mathbb{X}_L\mathbb{X}_M(\x)\mathbb{X}_R = \mathbb{X} \\ \mathbb{Y}_L(\y) \mathbb{Y}_R = \mathbb{Y}}} c(\x,\y ) f_1(0,\mathbb{X}_M, \x)f_2(\z,\mathbb{X}_L * \mathbb{Y}_L , \y, \mathbb{X}_R * \mathbb{Y}_R, \z') \\
&\hspace{5cm}+B_{\Y0}^{\X} + B_{0\Y'}^{\X} - B_{00}^{\X}
\end{aligned}
 \displaybreak[3] \\
&\begin{aligned}
\stackrel{\eqref{eq:Conv-harmonic-4}}{=}&
 \sum_{\substack{\mathbb{X}_L(\x)\mathbb{X}_M\mathbb{X}_R = \mathbb{X} \\ \mathbb{Y}_L(\y) \mathbb{Y}_R = \mathbb{Y}}} c(\x,\y ) f_1(\x,\mathbb{X}_M, 0)f_2(\z,\mathbb{X}_L * \mathbb{Y}_L , \y, \mathbb{X}_R * \mathbb{Y}_R, \z') \\
&\hspace{2cm}+ \sum_{\substack{\mathbb{X}_L\mathbb{X}_M(\x)\mathbb{X}_R = \mathbb{X} \\ \mathbb{Y}_L(\y) \mathbb{Y}_R = \mathbb{Y}}} c(\x,\y ) f_1(0,\mathbb{X}_M, \x)f_2(\z,\mathbb{X}_L * \mathbb{Y}_L , \y, \mathbb{X}_R * \mathbb{Y}_R, \z') \\
&+\sum_{\substack{\mathbb{X}_L\mathbb{X}_M \mathbb{X}_R = \mathbb{X} \\ \mathbb{Y}_L(\y)\mathbb{Y}_R = \mathbb{Y}}} f_1(\y,\mathbb{X}_M,\y)f_2(\z,(\mathbb{X}_L * \mathbb{Y}_L),\y,(\mathbb{X}_R * \mathbb{Y}_R),\z') \notag \\
&\hspace{2cm}+\sum_{\substack{\mathbb{X}_L(\x)\mathbb{X}_M \mathbb{X}_R = \mathbb{X} \\ \mathbb{Y}_L(\y)\mathbb{Y}_R = \mathbb{Y}}} c(\y,\x) f_1(\y,\mathbb{X}_M,0)f_2(\z,(\mathbb{X}_L * \mathbb{Y}_L),\y,(\mathbb{X}_R * \mathbb{Y}_R),\z')\\
&\hspace{4cm}+\sum_{\substack{\mathbb{X}_L\mathbb{X}_M(\x) \mathbb{X}_R \\ \mathbb{Y}_L(\y)\mathbb{Y}_R = \mathbb{Y}}} c(\y,\x) f_1(0,\mathbb{X}_M,\y)f_2(\z,(\mathbb{X}_L * \mathbb{Y}_L),\y,(\mathbb{X}_R * \mathbb{Y}_R),\z')
\end{aligned}
 \displaybreak[3] \\
&\begin{aligned}
=&\sum_{\substack{\mathbb{Y}_L(\y)\mathbb{Y}_R = \mathbb{Y} \\ \mathbb{X}_L(\x)\mathbb{X}_M(\x')\mathbb{X}_R = \mathbb{X}}}f_2(\z,\mathbb{X}_L * \mathbb{Y}_L , \y, \mathbb{X}_R * \mathbb{Y}_R, \z')\\
&\hspace{1cm}\times \bigg\{ f_1(\y,\x,\mathbb{X}_M,\x', \y)  + c(\x,\y ) \big(f_1(\x,\mathbb{X}_M,\x', 0) - f_1(\y,\mathbb{X}_M,\x',0)\big)  \\
& \hspace{5cm} + c(\x',\y)\big( f_1(0,\x,\mathbb{X}_M, \x') - f_1(0,\x,\mathbb{X}_M,\y)\big)  \bigg\} 
\end{aligned}\\
&\begin{aligned}
\stackrel{\eqref{eq:strong-parity-revised}}{=}&0
\end{aligned}
\end{align*}
as claimed.
\end{proof}

\begin{lem}\label{eq:Conv-harmonic-parity-y}
For $\Psi_1$ $\in \addmr_0 \cap \Fad \cap \V_{\strprty}$ and $\Psi_2\in \addmr_0$, we have
\begin{align*}
B_{\X0}^{\Y} + B_{0\X'}^{\Y} + B_{\Y0}^{\Y} + B_{0\Y'}^{\Y} - B_{00}^{\Y}=0.
\end{align*}
\end{lem}
\begin{proof}
By symmetry, it suffices to replace $\mathbb{X}$ with $\mathbb{Y}$ in the above lemmas.
\end{proof}

\begin{cor}\label{cor:conclusion}
For $\Psi_1$ $\in \addmr_0 \cap \Fad \cap \V_{\strprty}$ and $\Psi_2\in \addmr_0$, we have
\begin{align*}
\vimo_{d_{\Psi_1(\Psi_2)}}(\z, \mathbb{X} * \mathbb{Y}, \z')
=f_1(\z,\mathbb{X},\z')f_2(\z,\mathbb{Y},\z') + f_1(\z,\mathbb{Y},\z')f_2(\z,\mathbb{X},\z').
\end{align*}

\begin{proof}
This corollary immediately follows from Lemma \ref{lem:vimo-B}, Lemma \ref{eq:Conv-harmonic-z}, Lemma \ref{eq:Conv-harmonic-parity}, and Lemma \ref{eq:Conv-harmonic-parity-y}.
\end{proof}

\end{cor}

\begin{prop}[ = Lemma \ref{lem:essential}]
For $\Psi_1$ $\in \addmr_0 \cap \Fad \cap \V_{\strprty}$ and $\Psi_2\in \addmr_0$, we have
\[
\Delta_*(d_{\Psi_1}(\Psi_2)_{\#}) = d_{\Psi_1}(\Psi_2)_{\#} \otimes 1 + 1\otimes d_{\Psi_1}(\Psi_2)_{\#} + \Psi_{1.\#}\otimes \Psi_{2,\#} + \Psi_{2,\#}\otimes \Psi_{1,\#}.
\]
\end{prop}

\begin{proof}
Since both $\Delta_*$ and $*$ are dual maps relatively, we have
\begin{align*}
\langle \Delta_*(d_{\Psi_1}(\Psi_2)_{\#}) \mid w \otimes 1\rangle = \langle d_{\Psi_1}(\Psi_2)_{\#} \mid w * 1\rangle 
=\langle d_{\Psi_1}(\Psi_2)_{\#} \otimes 1 \mid w \otimes 1\rangle
\end{align*}
and
\begin{align*}
\langle \Delta_*(d_{\Psi_1}(\Psi_2)_{\#}) \mid 1 \otimes w\rangle = \langle d_{\Psi_1}(\Psi_2)_{\#} \mid 1 * w\rangle 
=\langle 1\otimes d_{\Psi_1}(\Psi_2)_{\#}  \mid1\otimes  w \rangle
\end{align*}
for any $w \in Y^*$.

Next, we investigate $\langle \Delta_*(d_{\Psi_1}(\Psi_2)) \mid w_1 \otimes w_2\rangle$ for any non-empty words $w_1$, $w_2 \in Y^*\setminus\{1\}$.
By Corollary \ref{cor:conclusion}, we have
\begin{align*}
&\langle \mathbf{q}^{\vee}( d_{\Psi_1}(\Psi_2)) \mid y_{l+1}((y_{k_{r+s}}\cdots y_{k_{r+1}}) * (y_{k_{r}}\cdots y_{k_{1}})) \rangle \\
&\begin{aligned}
=& \sum_{l_1 + l_2 =l} \langle \mathbf{q}^{\vee}(\Psi_1)\mid y_{l_1+1}y_{k_{r+s}}\cdots y_{k_{r+1}}\rangle \langle \mathbf{q}^{\vee}(\Psi_2)\mid y_{l_2+1}y_{k_{r}}\cdots y_{k_{1}}\rangle \\
&\hspace{1cm}+ \langle \mathbf{q}^{\vee}(\Psi_2)\mid y_{l_1+1}y_{k_{r+s}}\cdots y_{k_{r+1}}\rangle \langle \mathbf{q}^{\vee}(\Psi_1)\mid y_{l_2+1}y_{k_{r}}\cdots y_{k_{1}}\rangle
\end{aligned}
\end{align*}
for $l \in \mathbb{Z}_{\ge 0}$, $r$, $s\in\mathbb{Z}_{>0}$, $k_1,\ldots,k_{r+s} \in \mathbb{Z}_{>0}$.
Thus we have
\begin{align*}
&\langle d_{\Psi_1}(\Psi_2)_{\#} \mid ((y_{k_{r+s}}\cdots y_{k_{r+1}}) * (y_{k_{r}}\cdots y_{k_{1}})) \rangle \\
=&\langle \Psi_{1,\#}\mid y_{k_{r+s}}\cdots y_{k_{r+1}}\rangle \langle \Psi_{2,\#}\mid y_{k_{r}}\cdots y_{k_{1}}\rangle+ \langle \Psi_{2,\#}\mid y_{k_{r+s}}\cdots y_{k_{r+1}}\rangle \langle \Psi_{1,\#}\mid y_{k_{r}}\cdots y_{k_{1}}\rangle.
\end{align*}
Since both $\Delta_*$ and $*$ are dual maps relatively, we have
\begin{align*}
&\langle \Delta_*(d_{\Psi_1}(\Psi_2)_{\#}) \mid y_{k_{r+s}}\cdots y_{k_{r+1}}\otimes y_{k_r}\cdots y_{k_1}\rangle \\
=&\langle d_{\Psi_1}(\Psi_2)_{\#} \mid y_{k_{r+s}}\cdots y_{k_{r+1}}* y_{k_r}\cdots y_{k_1}\rangle \displaybreak[3] \\
=&\langle \Psi_{1,\#}\mid y_{k_{r+s}}\cdots y_{k_{r+1}}\rangle \langle \Psi_{2,\#}\mid y_{k_{r}}\cdots y_{k_{1}}\rangle+ \langle \Psi_{2,\#}\mid y_{k_{r+s}}\cdots y_{k_{r+1}}\rangle \langle \Psi_{1,\#}\mid y_{k_{r}}\cdots y_{k_{1}}\rangle \displaybreak[3]\\
=&\langle \Psi_{1,\#}\otimes \Psi_{2,\#}\mid y_{k_{r+s}}\cdots y_{k_{r+1}} \otimes y_{k_{r}}\cdots y_{k_{1}}\rangle
+ \langle \Psi_{2,\#}\otimes \Psi_{1,\#}\mid y_{k_{r+s}}\cdots y_{k_{r+1}} \otimes y_{k_{r}}\cdots y_{k_{1}}\rangle.
\end{align*}

Hence, it follows that
\begin{align*}
\begin{aligned}
\Delta_*(d_{\Psi_1}(\Psi_2)_{\#}) = d_{\Psi_1}(\Psi_2)_{\#} \otimes 1 + 1\otimes d_{\Psi_1}(\Psi_2)_{\#} + \Psi_{1.\#}\otimes \Psi_{2,\#} + \Psi_{2,\#}\otimes \Psi_{1,\#}
\end{aligned} 
\end{align*}
as claimed.
\end{proof}

\begin{center}

\textbf{Use of AI}

\end{center}

The author used ChatGPT 5.6 Sol for clarification of the exposition, correction of typographical errors, and improvement of proofs in this article. All mathematical content and revisions were carefully checked and verified by the author.

\begin{center}
\textbf{Acknowledgments}
\end{center}
The author would first like to express his deepest gratitude to Professor Minoru Hirose. Without his kindness and encouragement, the author would not have been able to establish the main theorem in this paper.  
He is also grateful to Dr. Nao Komiyama for giving him valuable advice on the theory of commutative power series.  
He would like to thank Dr. Annika Burmester, Prof. Henrik Bachmann, Prof. Hidekazu Furusho, Prof. David Jarossay, Dr. Iu-Iong Ng, Prof. Koji Tasaka, and Dr. Khalef Yaddaden, Shin-taro Yanagida for their helpful comments on his paper.
The author thanks Takeshi Shinohara, Kosei Watanabe, and Ku-Yu Fan for reviewing the Main Theorem.

\bibliographystyle{amsplain}

\bibliography{reference-adjoint}

\end{document}